\documentclass[11pt]{article}
\usepackage{amssymb,amsmath,pifont,amsfonts,amsthm,enumitem}
\usepackage{mathtools}
 \usepackage{lipsum}
\usepackage{multirow,boldline}
\usepackage[top=30mm,right=30mm,bottom=30mm,left=30mm]{geometry}
\usepackage{makecell}
\usepackage{tikz}\usetikzlibrary{shapes,fit}
\usepackage{scalerel}
\usepackage{relsize}
\usepackage{tikz-cd}
\usepackage{rotating}
\usepackage{dynkin-diagrams}
\usepackage{longtable} 
\usepackage[tableposition=above,width=12in]{caption}
\usepackage{arydshln} 
\usetikzlibrary{positioning}

\newcommand{\F}{\mathbb{F}}
\newcommand{\Z}{\mathbb{Z}}
\newcommand{\N}{\mathbb{N}}
\newcommand{\R}{\mathbb{R}}
\newcommand{\C}{\mathbb{C}}

\newcommand{\xmark}{\text{\ding{56}}}

\newcommand\simcal[1]{\stackrel{\widetilde{\textcolor{white}{aaaa}}}{\smash{\mathcal{#1}}\rule{0pt}{0.4ex}}}
\newcommand\shimcal[1]{\stackrel{\widetilde{\textcolor{white}{mmmmm}}}{\smash{\mathcal{#1}}\rule{0pt}{0.5ex}}}

\makeatletter
\newcommand{\imod}[1]{\allowbreak\mkern4mu({\operator@font mod}\,\,#1)}
\makeatother

\theoremstyle{plain}
\newtheorem{thm}{Theorem}
\newtheorem{cjr}{Conjecture}

\newtheorem{cor}{Corollary}

\newtheorem{Theorem}{Theorem}[section]
\newtheorem{Proposition}[Theorem]{Proposition}
\newtheorem{Lemma}[Theorem]{Lemma}
\newtheorem{Corollary}[Theorem]{Corollary}

\theoremstyle{definition}
\newtheorem{defn}{Definition}
\newtheorem{Definition}[Theorem]{Definition}
\newtheorem{Example}[Theorem]{Example}
\newtheorem{Remark}[Theorem]{Remark}

\newtheorem*{maintheorem!}{Theorem \ref{maintheorem!}}
\newtheorem*{classificationthm!}{Theorem \ref{classificationthm!}}
\newtheorem*{myconjecture!}{Conjecture \ref{myconjecture!}}
\newtheorem*{mydefinition!}{Definition \ref{mydefinition!}}
\newtheorem*{existence!}{Theorem \ref{existence!}}
\newtheorem*{maincor1}{Corollary \ref{maincor1}}
\newtheorem*{maincor2}{Corollary \ref{maincor2}}
\newtheorem*{maincor3}{Corollary \ref{maincor3}}

\DeclareMathOperator{\cd}{cd}
\DeclareMathOperator{\Gal}{Gal}

\DeclareMathOperator{\PSp}{PSp}

\DeclareMathOperator{\Spin}{Spin}

\DeclareMathOperator{\Iso}{Iso}

\DeclareMathOperator{\Det}{det}

\DeclareMathOperator{\Stab}{Stab}
\DeclareMathOperator{\Fix}{Fix}
\DeclareMathOperator{\Aut}{Aut}

\DeclareMathOperator{\SL}{SL}

\DeclareMathOperator{\GL}{GL}

\DeclareMathOperator{\GU}{GU}
\DeclareMathOperator{\PSU}{PSU}
\DeclareMathOperator{\AGL}{AGL}
\DeclareMathOperator{\Ad}{Ad}

\DeclareMathOperator{\Char}{char}

\DeclareMathOperator{\Span}{span}

\begin{document}
\title{Maximal connected $k$-subgroups of maximal rank in connected reductive algebraic $k$-groups}
 \author{Damian Sercombe}

\maketitle


\begin{abstract}

\noindent Let $k$ be any field and let $G$ be a connected reductive algebraic $k$-group. Associated to $G$ is an invariant first studied by Satake \cite{Sa0, Sa} and Tits \cite{T} that is called the index of $G$ (a Dynkin diagram along with some additional combinatorial information). Tits \cite{T} showed that the $k$-isogeny class of $G$ is uniquely determined by its index and the $k$-isogeny class of its anisotropic kernel $G_a$. For the cases where $G$ is absolutely simple, Satake \cite{Sa0} and Tits \cite{T} classified all possibilities for the index of $G$. Let $H$ be a connected reductive $k$-subgroup of maximal rank in $G$. We introduce an invariant of the $G(k)$-conjugacy class of $H$ in $G$ called the embedding of indices of $H \subset G$. This consists of the index of $H$ and the index of $G$ along with an embedding map that satisfies certain compatibility conditions. We introduce an equivalence relation called index-conjugacy on the set of $k$-subgroups of $G$, and observe that the $G(k)$-conjugacy class of $H$ in $G$ is determined by its index-conjugacy class and the $G(k)$-conjugacy class of $H_a$ in $G$. We show that the index-conjugacy class of $H$ in $G$ is uniquely determined by its embedding of indices. For the cases where $G$ is absolutely simple of exceptional type and $H$ is maximal connected in $G$, we classify all possibilities for the embedding of indices of $H \subset G$. Finally, we establish some existence results. In particular, we consider which embeddings of indices exist when $k$ has cohomological dimension $1$ (resp. $k=\R$, $k$ is $\mathfrak{p}$-adic).
 \end{abstract}

\section*{Introduction}

\noindent The classification of finite dimensional simple complex Lie algebras into four infinite classical families $A_n$, $B_n$, $C_n$, $D_n$ and five exceptional cases $G_2$, $F_4$, $E_6$, $E_7$, $E_8$ was completed by Killing and Cartan in the late 19th century. 
This was refined by Dynkin in 1947 using the notion of Dynkin diagrams. In the 1950s Borel and Chevalley investigated affine algebraic groups and showed that simple algebraic groups over an algebraically closed field are also classified up to isogeny by their Dynkin diagrams. 
For fields that are not algebraically closed, the situation is more complicated, and was studied by Tits and Satake in the 1960s. All algebraic groups in this paper are considered to be smooth and affine.

\vspace{2mm}\noindent Let $k$ be any field and let $G$ be a connected reductive algebraic $k$-group. Associated to $G$ is an invariant $\mathcal{I}(G)$ called the \textit{index} of $G$ (a Dynkin diagram along with some additional combinatorial information, see $\S \ref{index}$ for a precise definition). This invariant was first studied by Satake \cite{Sa0, Sa} and Tits \cite{T} in the 1960s. Following on from work of Satake \cite{Sa0,Sa}, Tits \cite{T} showed that $G$ is uniquely determined up to $k$-isogeny by its index and the $k$-isogeny class of its anisotropic kernel $G_a$. Moreover, in Table II of \cite{T}, Tits lists all possibilities for the index of an absolutely simple algebraic $k$-group. Combining these results, for any field $k$ we essentially have a classification `up to anisotropy' of connected reductive algebraic $k$-groups. If $k$ is separably closed, $k$ is finite or $k=\R$ then the anisotropic kernel of $G$ is well understood and we have a `complete' classification of connected reductive algebraic $k$-groups. For most other fields, however, a classification of anisotropic connected reductive algebraic $k$-groups remains far out of reach.

\vspace{2mm}\noindent An important problem is to classify the maximal $k$-subgroups $H$ of $G$ up to conjugacy by some element of $G(k)$. The maximal $k$-subgroups of $G$ that are parabolic are well understood, and can be read immediately off of the index of $G$ (see Proposition $21.12$ of \cite{B}). If $k$ is perfect and $H$ is not parabolic then $H$ is reductive by Corollary $3.7$ of \cite{BT} (Borel-Tits). 
For imperfect $k$, the maximal $k$-subgroups of $G$ that are neither parabolic nor reductive have been investigated in \cite{G1}. We are interested in the cases where $H$ is reductive.

\vspace{2mm}\noindent For $G$ absolutely simple, the maximal $k$-subgroups of $G$ have (mostly) been classified when $k$ is separably closed (Borel-de Siebenthal \cite{BD}, Dynkin \cite{D,D1}, Seitz \cite{Se, Se1}, Liebeck-Seitz \cite{LS,LS2,LS4}, Testerman \cite{Te}), when $k$ is finite (Aschbacher \cite{As}, Liebeck-Seitz \cite{LS,LS2,LS3}, Liebeck-Saxl-Seitz \cite{LSS}) and when $k=\R$ (Komrakov \cite{K}, Karpelevich \cite{Ka}, Taufik \cite{Ta}, de Graff-Marrani \cite{DM} and others). However, not much is known for other fields.

\vspace{2mm}\noindent In this paper, for the most part, we continue in the spirit of Tits and Satake and work with an arbitrary field $k$. We consider the problem of classifying connected reductive $k$-subgroups $H$ of maximal rank in any connected reductive algebraic $k$-group $G$ up to $G(k)$-conjugacy. We call this the \textit{classification problem}. The subproblem where $H$ is $k$-anisotropic is called the \textit{anisotropic classification problem} -- it appears to depend intrinsically on the choice of field $k$. We develop techniques to reduce the classification problem to the anisotropic classification problem in full generality (i.e. for any given $k$). In such generality, this seems to be as far as one can go.

\vspace{2mm}\noindent We pay special attention to the cases where $G$ is absolutely simple of exceptional type and $H$ is maximal among connected $k$-subgroups of $G$ (a.k.a \textit{maximal connected}). The cases where $G$ is absolutely simple of classical type will be the subject of a future paper. If $k$ is separably closed, $k$ is finite or $k=\R$, then our results reduce to what is already known in the literature. We also find some new results for the case where $k$ is $\mathfrak{p}$-adic.


\vspace{2mm}\noindent Henceforth we use the following setup. Some of the ensuing definitions are quite technical, and so we provide a more detailed exposition of this setup in Sections \ref{Boreltits}, \ref{combinatoricschap} and the beginning of Section \ref{mainchapter}. 

\vspace{2mm}\noindent Let $k$ be any field. Let $\overline{k}$ be an algebraic closure of $k$ and let $K$ be the separable closure of $k$ in $\overline{k}$. Let $\Gamma=\Gal(K/k)$ be the absolute Galois group of $k$. Let $G$ be a connected reductive (affine) algebraic $k$-group. We will see in Section \ref{Boreltits} that $k$-groups come equipped with a natural action of $\Gamma$.

\vspace{2mm}\noindent In this paragraph we define the index of $G$. Let $S$ be a maximal $k$-split torus of $G$. A result of Grothendieck tells us that there exists a maximal $k$-torus $T$ of $G$ that contains $S$. We endow the character group $X(T)$ with a total order $<$ that is compatible with the restriction map $X(T) \to X(S)$. Let $\Lambda$ be the system of simple roots for $G$ with respect to $T$ that is compatible with $<$. Let $\Lambda_0$ be the subset of $\Lambda$ that vanishes on $S$. Denote $W:=N_G(T)/T$ and $\smash{\tilde{W}}:=(N_G(T) \rtimes \Gamma)/T$. There exists a $\smash{\tilde{W}}$-invariant inner product $(\cdot \hspace{0.5mm},\cdot)$ on $X(T)_{\R}$. There is a natural action $\hat{\iota}:\Gamma \to \GL\!\big(X(T)_{\R}\big)$ that stabilises $\Lambda$ and $\Lambda_0$ and is by isometries, which we call the \textit{Tits action}. The \textit{index} of $G$ (with respect to $T$ and $<$) 
is the quadruple $\mathcal{I}(G)=\big(X(T)_{\R}, \Lambda, \Lambda_0, \hat{\iota}(\Gamma) \big)$. Satake \cite{Sa0} showed that, as a combinatorial object, $\mathcal{I}(G)$ is independent of the choice of $T$ and $<$. We illustrate $\mathcal{I}(G)$ using a modification of a Dynkin diagram called a Tits-Satake diagram.

\vspace{2mm}\noindent One can define a root system purely combinatorially (that is, without reference to an algebraic group). Such an object is called an abstract root system. Similarly, we present a purely combinatorial definition of an index (essentially due to Tits \cite{T}). We will reuse some notation as we move back and forth between algebraic and combinatorial definitions.

\vspace{2mm}\noindent An \textit{abstract index} $\mathcal{I}=(E,\Delta,\Delta_0,\Pi)$ is a quadruple consisting of a finite-dimensional real inner product space $E$, a system of simple roots $\Delta$ in $E$, a subset $\Delta_0$ of $\Delta$ and a subgroup $\Pi$ of the isometry group of $E$ that stabilises both $\Delta$ and $\Delta_0$, such that certain conditions are satisfied (see Definition \ref{abstractindexdefn}). Satake \cite{Sa0} and Tits \cite{T} showed that the index of $G$ is indeed an abstract index.

\vspace{2mm}\noindent We use the following notation: $\langle \Delta \rangle$ is the root system in $E$ that is generated by $\Delta$, $\langle \Delta \rangle^+$ is the associated set of positive roots, $W_{\Delta}$ is the Weyl group of $\Delta$, $E_a$ is the smallest subspace of $E$ that contains $\Delta_0$ and that has an orthogonal complement in $E$ that is fixed pointwise by $\Pi$, $E_{\Delta}$ is the subspace of $E$ that is spanned by $\Delta$, $\overline{E}$ is the orthogonal complement of $E_{\Delta}$ in $E$ and $\overline{E}_a:=\overline{E} \cap E_a$. If $\Delta=\Delta_0$ and $\overline{E}_a$ is trivial then $\mathcal{I}$ is \textit{anisotropic}, otherwise $\mathcal{I}$ is \textit{isotropic}. If $\Delta_0=\varnothing$ and $\Pi$ is trivial then $\mathcal{I}$ is \textit{split}. If $\Delta_0 = \varnothing$ then $\mathcal{I}$ is \textit{quasisplit}. If $\Delta$ is irreducible then $\mathcal{I}$ is \textit{irreducible}. If $\Delta$ is of classical (resp. exceptional) type then $\mathcal{I}$ is of \textit{classical} (resp. \textit{exceptional}) \textit{type}. Let $I$ be a $\Pi$-stable subset of $\Delta$ that contains $\Delta_0$. The abstract index $(E,I, \Delta_0, \Pi)$ is a \textit{subindex} of $\mathcal{I}$. The \textit{minimal subindex} of $\mathcal{I}$ is $\mathcal{I}_m:=(E,\Delta_0, \Delta_0, \Pi)$.

\vspace{2mm}\noindent Let $\smash{'}\mathcal{I}=(\smash{'}\hspace{-0.3mm}E,\smash{'}\hspace{-0.4mm}\Delta,\smash{'}\hspace{-0.4mm}\Delta_0,\smash{'}\Pi)$ be another abstract index. An \textit{isomorphism} from $\mathcal{I}$ to $\smash{'}\mathcal{I}$ is a bijective isometry $\psi:E \to \smash{'}\hspace{-0.3mm}E$ that satisfies $\psi(\Delta)=\smash{'}\hspace{-0.4mm}\Delta$, $\psi(\Delta_0)=\smash{'}\hspace{-0.4mm}\Delta_0$ and $\psi\Pi \psi^{-1}=\smash{'}\Pi$. We denote $\psi(\mathcal{I}):=\smash{'}\mathcal{I}$. In Table II of \cite{T}, Tits classifies all isomorphism classes of irreducible abstract indices. 

\vspace{2mm}\noindent An abstract index $\mathcal{I}$ is \textit{$k$-admissible} if there exists a connected reductive $k$-group $X$ such that the index of $X$ is isomorphic to $\mathcal{I}$. For example, $\dynkin{F}{II}$ is $\R$-admissible but not $\F_q$-admissible for any finite field $\F_q$ (refer to Table \ref{F_4}). In $\S 3$ of \cite{T}, Tits shows that any abstract index is $k$-admissible for some field $k$.

\vspace{2mm}\noindent Let $H$ be a connected reductive proper $k$-subgroup of maximal rank in $G$. We now introduce a (novel) invariant of the $G(k)$-conjugacy class of $H$ in $G$, that we call the embedding of indices of $H \subset G$.

\vspace{2mm}\noindent Let $S_H$ be a maximal $k$-split torus of $H$ and let $T_H$ be a maximal $k$-torus of $H$ that contains $S_H$. Without loss of generality, we can assume that $S_H \subseteq S$ and that $<$ is compatible with the restriction map $X(T) \to X(S_H)$. Consider the Levi $k$-subgroup $L:=C_G(S_H)$ of $G$. There exists $g \in L(K)$ such that $(T_H)^g=T$. Let $\theta:X(T_H)_{\R} \to X(T)_{\R}$ be defined by $\chi \mapsto \chi \cdot g^{-1}$ for $\chi \in X(T_H)$. Denote $\smash{\tilde{W}_H}:=(N_H(T_H) \rtimes \Gamma)/T_H$. There exists a $\smash{\tilde{W}_H}$-invariant inner product $(\cdot \hspace{0.5mm},\cdot)_H$ on $X(T_H)_{\R}$ such that $\theta$ is a bijective isometry.

\vspace{2mm}\noindent Let $\chi_1, \chi_2 \in X(T_H)$. The total order $<$ on $X(T)$ induces a total order $<_{g^{-1}}$ on $X(T_H)$ that is defined by $\chi_1<_{g^{-1}} \chi_2$ if $\chi_1 \cdot g < \chi_2 \cdot g$. Let $\Delta$ be the system of simple roots for $H$ with respect to $T_H$ that is compatible with $<_{g^{-1}}$. Let $\Delta_0$ be the subset of $\Delta$ that vanishes on $S_H$ and let $\hat{\iota}_H:\Gamma \to \GL\!\big(X(T_H)_{\R}\big)$ be the Tits action (stabilising $\Delta$ and $\Delta_0$). The index of $H$ is $\mathcal{I}(H)=\big(X(T_H)_{\R}, \Delta, \Delta_0, \hat{\iota}_H(\Gamma)\big)$. We call the triple $\big(\mathcal{I}(G),\mathcal{I}(H),\theta\big)$ an \textit{embedding of indices} of $H \subset G$. A priori, this triple depends on several choices (i.e. $T$, $T_H$, $<$ and $g$). However, we will see in Theorem \ref{maintheorem!} that $\big(\mathcal{I}(G),\mathcal{I}(H),\theta\big)$ is a combinatorial invariant of $H \subset G$.

\vspace{2mm}\noindent We introduce a (novel) equivalence relation on the set of $k$-subgroups of $G$. Two $k$-subgroups $A$ and $B$ of $G$ are \textit{index-conjugate} if there exists $x \in G(K)$ such that $A^x=B$, $x^{-1}x^{\Gamma} \subset A$ and, for any maximal $k$-split torus $S_A$ of $A$, $(S_A)^x$ is a maximal $k$-split torus of $B$. It is easy to see that index-conjugacy is an equivalence relation on the set of $k$-subgroups of $G$ that is finer than $G(K)$-conjugacy but coarser than $G(k)$-conjugacy. In fact, index-conjugacy restricts to $G(k)$-conjugacy (resp. $G(K)$-conjugacy) on the subset of $k$-split (resp. $k$-anisotropic) $k$-subgroups of $G$. It follows that the $G(k)$-conjugacy class of $H$ in $G$ is determined by its index-conjugacy class and the $G(k)$-conjugacy class of its anisotropic kernel $H_a$ in $G$.

\vspace{2mm}\noindent Inspired by the notion of an abstract index, we introduce a purely combinatorial definition of an embedding of indices. This definition is based on the observation that an embedding of $k$-groups of the same rank induces an embedding of their respective character spaces that satisfies certain conditions on their root systems.

\vspace{2mm}\noindent Let $\mathbb{P}$ denote the set of prime numbers. For $p \in \mathbb{P} \cup \{0\}$, a $p$-closed subsystem of a root system is a slight generalisation of a closed subsystem (see $\S \ref{apsubsystems}$ for the precise definition). 

\begin{defn}\label{mydefinition!} Let $p \in \mathbb{P} \cup \{0\}$. A \textit{$p$-embedding of abstract indices} is a triple that consists of two abstract indices $\mathcal{G}=(F,\Lambda,\Lambda_0, \Xi)$ and $\mathcal{H}=(E,\Delta,\Delta_0,\Pi)$, and a bijective isometry $\theta:E \to F$ that satisfies the following conditions:

\vspace{1.5mm}\noindent $(A.1)$ $\Phi:=\langle \Lambda \rangle$ is $\Pi^{\theta}$-stable, $\theta(\Delta) \subset \Phi^+$, $\Psi:=\langle \theta(\Delta) \rangle$ is a $p$-closed subsystem of $\Phi$, $\Lambda_a:= \Lambda \cap \theta(E_a)$ is a base of $\Phi_a:=\Phi \cap \theta(E_a)$ and ${}_{in}\Lambda_a \subseteq \theta(\Delta_0)$ (where ${}_{in}\Lambda_a$ is the union of irreducible components of $\Lambda_a$ that is maximal with respect to the condition ${}_{in}\Lambda_a \subseteq \langle \theta(\Delta_0) \rangle =:\Psi_0$).

\vspace{1.5mm}\noindent $(A.2)$ For every $\sigma \in \Pi$ there exists a unique $w_{\sigma} \in W_{\Lambda_a}$ such that $w_{\sigma} \sigma^{\theta} \in \Xi$. Moreover, the map $\Pi \to \Xi$ given by $\sigma \mapsto w_{\sigma} \sigma^{\theta}$ is surjective.

\vspace{1.5mm}\noindent $(A.3)$ $\Lambda_a$ is $\Xi$-stable and contains $\Lambda_0$.

\vspace{1.5mm}\noindent $(A.4)$ ${}_{in}\Lambda_a$ is contained in $\Lambda_0$.
\end{defn}

\noindent An \textit{embedding of abstract indices} is a $p$-embedding of abstract indices for some $p \in \mathbb{P} \cup \{0\}$. 

\vspace{2mm}\noindent Henceforth let $p \in \mathbb{P} \cup \{0\}$ and let $(\mathcal{G},\mathcal{H},\theta)$ be a $p$-embedding of abstract indices. We introduce some terminology associated to Definition \ref{mydefinition!}.

\vspace{2mm}\noindent The subindex $\mathcal{L}:=(F,\Lambda_a,\Lambda_0, \Xi)$ of $\mathcal{G}$ is called the \textit{$(\mathcal{H},\theta)$-embedded subindex of $\mathcal{G}$}. We say $(\mathcal{G},\mathcal{H},\theta)$ is \textit{(an)isotropic} if $\mathcal{H}$ is (an)isotropic. We say $(\mathcal{G},\mathcal{H},\theta)$ is \textit{split} (resp. \textit{quasisplit}) if $\mathcal{H}$ and $\mathcal{G}$ are both split (resp. both quasisplit). The map $\theta$ is an \textit{embedding} of $\mathcal{H}$ in $\mathcal{G}$. We often abuse notation and identify $\mathcal{H}$ with its image under $\theta$. 

\vspace{2mm}\noindent Let $V$ be the largest subspace of $F$ that is contained in the span of $\Phi$ and that is perpendicular to all roots in $\Psi$. If the fixed point subspace of $V$ under the action of $\Pi^{\theta}$ is trivial and $\Psi$ is maximal among $p$-closed $\Pi^{\theta}$-stable subsystems of $\Phi$ then $(\mathcal{G},\mathcal{H},\theta)$ is \textit{maximal}. If ${}_{in}\Lambda_a=\Lambda_a$ then $(\mathcal{G},\mathcal{H},\theta)$ is \textit{independent}.

\vspace{2mm}\noindent Let $(\smash{'}\mathcal{G},\smash{'}\mathcal{H},\smash{'}\hspace{-0.1mm}\theta)$ be another $p$-embedding of abstract indices. An \textit{isomorphism} from $(\mathcal{G},\mathcal{H},\theta)$ to $(\smash{'}\mathcal{G},\smash{'}\mathcal{H},\smash{'}\hspace{-0.1mm}\theta)$ is a bijective isometry $\phi:F \to \smash{'}\hspace{-0.5mm}F$ such that $\phi(\mathcal{G})=\smash{'}\mathcal{G}$ and $\phi \theta(\mathcal{H})=\smash{'}\hspace{-0.1mm}\theta(\smash{'}\mathcal{H})$. If such an isomorphism exists then $(\mathcal{G},\mathcal{H},\theta)$ and $(\smash{'}\mathcal{G},\smash{'}\mathcal{H},\smash{'}\hspace{-0.1mm}\theta)$ are \textit{isomorphic}. If $F=\smash{'}\hspace{-0.5mm}F$ and $\phi \in W_{\Lambda}$ then $\phi$ is a \textit{conjugation} from $(\mathcal{H},\theta)$ to $(\smash{'}\mathcal{H},\smash{'}\hspace{-0.1mm}\theta)$ in $\mathcal{G}$. If such a conjugation exists then $(\mathcal{H},\theta)$ and $(\smash{'}\mathcal{H},\smash{'}\hspace{-0.1mm}\theta)$ are \textit{conjugate} in $\mathcal{G}$.

\vspace{2mm}\noindent Given a field $k$ with $p=\Char(k)$, we say that $(\mathcal{G},\mathcal{H},\theta)$ is \textit{$k$-admissible} if there exists a pair of connected reductive $k$-groups $H_1 \subset G_1$ that share a maximal torus such that the embedding of indices of $H_1 \subset G_1$ is isomorphic to $(\mathcal{G},\mathcal{H},\theta)$.

\vspace{2mm}\noindent We can now state the main results of this paper. 

\vspace{2mm}\noindent Our first result gives a combinatorial description of the index-conjugacy classes of connected reductive $k$-subgroups of maximal rank in $G$.

\begin{thm}\label{maintheorem!} Let $k$ be a field and let $p=\Char(k)$. Let $G$ be a connected reductive algebraic $k$-group and let $H$ be a connected reductive $k$-subgroup of maximal rank in $G$. Let $\big(\mathcal{I}(G),\mathcal{I}(H),\theta\big)$ be an embedding of indices of $H \subset G$. Then

\vspace{2mm}\noindent $(i)$ $\big(\mathcal{I}(G),\mathcal{I}(H),\theta\big)$ is a $p$-embedding of abstract indices, and

\vspace{1mm}\noindent $(ii)$ $H$ is maximal connected in $G$ if and only if $\big(\mathcal{I}(G),\mathcal{I}(H),\theta\big)$ is maximal.

\vspace{2mm}\noindent In addition, if $\smash{'}\hspace{-0.4mm}H$ is another connected reductive $k$-subgroup of maximal rank in $G$ and $\big(\mathcal{I}(G),\mathcal{I}(\smash{'}\hspace{-0.4mm}H),\smash{'}\hspace{-0.2mm}\theta\big)$ is an embedding of indices of $\smash{'}\hspace{-0.4mm}H \subset G$, then

\vspace{2mm}\noindent $(iii)$ $H$ is index-conjugate to $\smash{'}\hspace{-0.4mm}H$ in $G$ if and only if $\big(\mathcal{I}(H),\theta\big)$ is conjugate to $\big(\mathcal{I}(\smash{'}\hspace{-0.4mm}H),\smash{'}\hspace{-0.2mm}\theta\big)$ in $\mathcal{I}(G)$.
\end{thm}


\noindent We know from Proposition $2$ of \cite{Sa0} that the isomorphism class of $\mathcal{I}(G)$ is an invariant of $G$. As a consequence of Theorem \ref{maintheorem!}$(i)$ and $(iii)$, we see that the conjugacy class of $\big(\mathcal{I}(H),\theta\big)$ in $\mathcal{I}(G)$ is an invariant of the $G(k)$-conjugacy class of $H$ in $G$.

\vspace{2mm}\noindent It is well known that, over any field $k$, any abstract root system can be realised as the root system of some connected reductive algebraic $k$-group. In other words, any abstract root system is $k$-admissable for all fields $k$. A weaker analogue holds for indices, in that any abstract index is $k$-admissible for some field $k$. The following conjecture is an analogue of those two results. In essence, it is claiming that Definition \ref{mydefinition!} is as restrictive as is possible subject to the constraint that Theorem \ref{maintheorem!}$(i)$ holds.

\begin{cjr}\label{myconjecture!} Let $p \in \mathbb{P} \cup \{0\}$ and let $(\mathcal{G},\mathcal{H},\theta)$ be a $p$-embedding of abstract indices. Then $(\mathcal{G},\mathcal{H},\theta)$ is $k$-admissible for some field $k$ with characteristic $p$. 
\end{cjr}

\noindent The following task arises naturally from Theorem \ref{maintheorem!}. For each $p \in \mathbb{P} \cup \{0\}$ and each irreducible abstract index $\mathcal{G}$, we wish to classify all conjugacy classes of maximal $p$-embeddings of abstract indices in $\mathcal{G}$. Due to the size of this task, in Theorem \ref{classificationthm!} we consider only the cases where $\mathcal{G}$ is of exceptional type. The cases where $\mathcal{G}$ is of classical type will be the subject of a future paper. The tables referred to in Theorem \ref{classificationthm!} can be found in the Appendix.

\begin{thm}\label{classificationthm!} For each $p \in \mathbb{P} \cup \{0\}$ and each (isomorphism class of) abstract index $\mathcal{G}$ of exceptional type, all $\mathcal{G}$-conjugacy classes of isotropic maximal $p$-embeddings of abstract indices are classified in Tables \ref{G_2}, \ref{F_4}, \ref{E_6}, \ref{E_7} and \ref{E_8} for the cases where $\mathcal{G}$ is of type $G_2$, $F_4$, $E_6$, $E_7$ and $E_8$ respectively.
\end{thm}

\noindent Note that in Theorem \ref{classificationthm!} we ignore the anisotropic embeddings. This is because they are easy, since axiom $(A.3)$ of Definition \ref{mydefinition!} becomes vacuous (see Remark \ref{isoaniso}).

\vspace{2mm}\noindent Combining Theorems \ref{maintheorem!} and \ref{classificationthm!} (and the above remark), for any field $k$ and any absolutely simple $k$-group $G$ of exceptional type, we have a classification of all possible index-conjugacy classes of reductive maximal connected $k$-subgroups of maximal rank in $G$. However, not all of these embeddings exist over a given field $k$. So we now turn our attention to the \textit{existence problem}, that is, the problem of determining which (isomorphism classes of) embeddings of indices exist over a given field $k$.

\vspace{2mm}\noindent For the next result, given a $p$-embedding of abstract indices $(\mathcal{G},\mathcal{H},\theta)$, recall that $\mathcal{L}$ denotes the $(\mathcal{H},\theta)$-embedded subindex of $\mathcal{G}$ and $\mathcal{H}_m$ denotes the minimal subindex of $\mathcal{H}$. For $G$ a $k$-group with index $\mathcal{G}$, subindices of $\mathcal{G}$ correspond to Levi $k$-subgroups of $G$ (see Proposition $21.12$ of \cite{B}). So, in the statement of Theorem \ref{existence!}, if $G$ exists then certainly so does $L$.

\begin{thm}\label{existence!} Let $p \in \mathbb{P} \cup \{0\}$ and let $(\mathcal{G},\mathcal{H},\theta)$ be a $p$-embedding of abstract indices. Let $k$ be a field of characteristic $p$. Let $G$ be a connected reductive algebraic $k$-group with index isomorphic to $\mathcal{G}$. Let $L$ be a Levi $k$-subgroup of $G$ with index isomorphic to $\mathcal{L}$. There exists a $k$-subgroup of $G$ with embedding of indices isomorphic to $(\mathcal{G},\mathcal{H},\theta)$ if and only if there exists a $k$-subgroup of $L$ with embedding of indices isomorphic to $(\mathcal{L},\mathcal{H}_m,\theta)$.
\end{thm}

\noindent Combining Theorems \ref{maintheorem!}, \ref{classificationthm!} and \ref{existence!}, we have reduced the classification problem to the anisotropic classification problem for any field $k$ and any absolutely simple $k$-group $G$ of exceptional type. 

\vspace{2mm}\noindent Theorem \ref{existence!} has some interesting consequences.

\begin{cor}\label{maincor1} Let $p \in \mathbb{P} \cup \{0\}$. Let $(\mathcal{G},\mathcal{H},\theta)$ be an independent $p$-embedding of abstract indices. Let $k$ be a field of characteristic $p$ such that $\mathcal{G}$ is $k$-admissible. Then $(\mathcal{G},\mathcal{H},\theta)$ is $k$-admissible.
\end{cor}

\begin{cor}\label{maincor2} Let $p \in \mathbb{P} \cup \{0\}$. Let $(\mathcal{G},\mathcal{H},\theta)$ be a quasisplit $p$-embedding of abstract indices. Let $k$ be a field of characteristic $p$ such that $\mathcal{G}$ is $k$-admissible. Then $(\mathcal{G},\mathcal{H},\theta)$ is $k$-admissible if and only if $\Pi$ is a quotient of the absolute Galois group $\Gamma$ of $k$.
\end{cor}

\noindent For our final result, we consider the existence problem for three specific classes of fields.

\begin{cor}\label{maincor3} Let $(\mathcal{G},\mathcal{H},\theta)$ be an embedding of abstract indices that is listed in Table \ref{G_2}, \ref{F_4}, \ref{E_6}, \ref{E_7} or \ref{E_8}. If there is a $\checkmark$ in column $8$ (resp. $9$, $10$) of the corresponding row then $(\mathcal{G},\mathcal{H},\theta)$ is $k$-admissible for some field $k$ with cohomological dimension $1$ (resp. $k=\R$, $k$ is $\mathfrak{p}$-adic). If there is a $\xmark$ then $(\mathcal{G},\mathcal{H},\theta)$ is not $k$-admissible for such a field.

\vspace{2mm}\noindent [In the cases where we list more than one value for $\Pi$ in a single row of one of the tables, we put a $\checkmark$ in the corresponding column if at least one of these values for $\Pi$ corresponds to a $k$-admissible $(\mathcal{G},\mathcal{H},\theta)$.]
\end{cor}

\noindent There are a few $\mathfrak{p}$-adic cases for which we are unsure of the answer. Such cases are indicated by a $?$ symbol in the relevant entry of the table. Our results in Corollary \ref{maincor3} agree with Table $5.1$ of \cite{LSS} for the cases where $k$ is a finite field and with Tables $4$-$62$ of \cite{K} and Tables $1$-$134$ of \cite{DM} for the case $k=\R$. We are unaware of any existing reference for the $\mathfrak{p}$-adic cases.

\vspace{2mm}\noindent We illustrate our results with the following excerpt from Table \ref{F_4}, for the case where $\Delta=B_4$ and $\Lambda=F_4$.

\begin{table}[!htb]\begin{center}\begin{tabular}{| c | c | c | c | c | c | c | c |}    \hline
\multicolumn{1}{|c|}{} & \multicolumn{1}{c|}{} & \multicolumn{1}{c|}{} & \multicolumn{1}{c|}{} & \multicolumn{1}{c|}{} & \multicolumn{3}{c|}{\small{Special fields}} \\

    $\Delta$ & $\mathcal{H}$ & $\Psi_0$ & $\Phi_a$ & $\mathcal{G}$ & \small{$\hspace{-0.8mm}\!\cd 1 \!\hspace{-1mm}$} & $\R$ & $\hspace{-1mm}\!\!Q_p\!\!\hspace{-1mm}$  \\ \hline \hline

    \multirow{4.2}{*}{$B_4$} & $\dynkin{B}{oooo}$ & $\varnothing$ & $\varnothing$ & $\dynkin{F}{I}$ & $\checkmark$ & $\checkmark$ & $~\checkmark$ \topstrut  \\
& $\dynkin{B}{ooo*}$ & $\smash{\widetilde{A_1}}$ & $\smash{\widetilde{A_1}}$ &&&&  \\ 
& $\dynkin{B}{oo**}$ & $B_2$ & $B_2$ &&&&  \\ 
& $\dynkin{B}{o***}$ & $B_3$ & $B_3$ & $\dynkin{F}{II}$ & \ding{55} & $\checkmark$ & \ding{55} \\ \hline
  \end{tabular}\end{center}\end{table}

\vspace{2mm}\noindent Assume that $G$ is an absolutely simple $k$-group of type $F_4$ and $H$ is a $k$-isotropic reductive maximal connected $k$-subgroup of $G$ of type $B_4$.

\vspace{2mm}\noindent By Theorems \ref{maintheorem!} and \ref{classificationthm!}, there are only two possibilities for the pair of indices $\big(\mathcal{I}(G),\mathcal{I}(H)\big)$. Either $\mathcal{I}(G) \cong \dynkin{F}{I}$ and $\mathcal{I}(H)\cong\dynkin{B}{oooo}$ (that is, both $G$ and $H$ are $k$-split) or $\mathcal{I}(G)\cong\dynkin{F}{II}$ and $\mathcal{I}(H)\cong\dynkin{B}{o***}$. Furthermore, there exists at most one embedding of each type (up to index-conjugacy). Let $S_H$ be a maximal $k$-split torus of $H$ and let $L:=C_G(S_H)$. If $G$ and $H$ are both $k$-split then $S_H=L$, otherwise $L$ has root system $B_3$.

\vspace{2mm}\noindent It follows from Theorem \ref{existence!} that the split embedding exists over any field $k$, and the other embedding exists over any field $k$ for which the index $\dynkin{F}{II}$ is $k$-admissible. For example, using Corollary \ref{maincor3}, if $k=\R$ then both of these types of embeddings exist: $\Spin(5,4) \subset F_{4(4)}$ is the split embedding and $\Spin(8,1) \subset F_{4(-20)}$ is the other one (this notation for real forms of exceptional groups is standard, see for instance \cite{DM}). If the cohomological dimension of $k$ is $1$ (e.g. $k$ is finite) or if $k$ is $\mathfrak{p}$-adic then only the split embedding exists. For example, if $k=\F_q$ for some prime power $q$ then $\Spin_9(q) \subset F_4(q)$ is the split embedding.

\vspace{2mm}\noindent We conclude with the following remark. A $k$-subgroup of $G$ is \textit{regular} if it is normalised by a maximal $k$-torus of $G$. Most of the constructions in this paper can immediately be extended to a regular reductive (not necessarily connected) $k$-subgroup of $G$ by replacing $H$ with $N_G(H)^{\circ}$. With a little more work, we could prove an appropriate generalisation of Theorem \ref{maintheorem!} for regular reductive $k$-subgroups of $G$. This will be considered in a future paper.

\section{Preliminaries}\label{Boreltits}

\noindent In this section we present some theory of a connected reductive algebraic group $G$ over an arbitrary field $k$. Many of the results in this section should be familiar to experts in the field (however, we have found it hard to source some of the proofs). 

\vspace{2mm}\noindent In $\S \ref{*action}$ we describe an action of the absolute Galois group $\Gamma$ of $k$ on the character space $X(T)_{\R}$ of a maximal $k$-torus $T$ of $G$. In $\S \ref{subspstruct}$ we investigate the subspace structure of $X(T)_{\R}$. In $\S \ref{Gammabase}$ we construct a linear order on $X(T)$ that is compatible with the action of $\Gamma$. In $\S \ref{index}$ we introduce an important invariant of $G$ called its index. Finally, in $\S \ref{orientationmaps}$ we investigate a functorial relationship between maximal tori of $G$ and their respective character spaces.

\vspace{2mm}\noindent Let $k$ be a field. Let $p$ be the characteristic of $k$. Let $\overline{k}$ be an algebraic closure of $k$ and let $K$ be the separable closure of $k$ in $\overline{k}$. Let $\Gamma=\Gal(K/k)$ be the absolute Galois group of $k$. 

\vspace{2mm}\noindent One can define a (linear) algebraic $k$-group in several equivalent ways. We follow Borel \cite{B} (AG.$\S 12$, I.$\S 1$) and say that a group $G$ is an \textit{algebraic $k$-group} (a.k.a. \textit{defined over $k$}) if it admits the structure of an affine algebraic $k$-set such that the group operations of multiplication and inversion are defined over $k$. As discussed in AG.$\S 14.3$ of \cite{B}, $G$ comes equipped with a natural action $\star$ of $\Gamma$ on $G(K)$ with fixed point subgroup $G(k)$. A \textit{$k$-morphism} of algebraic $k$-groups is a group homomorphism that is defined over $k$ as a map of algebraic sets.

\begin{Lemma}\label{groupautaction} The $\star$-action of $\Gamma$ on $G(K)$ is by abstract group automorphisms.
\begin{proof} As in $\S 1.1$ of \cite{MT}, we identify $\GL_N$ (for $N \in \N$) with the closed subset $$\big\{ \big(A:=(a_{ij}),y\big) \in \overline{k}{}^{N^2} \times \overline{k}\hspace{0.5mm} \big| \hspace{0.5mm}1 \leq i,j \leq N; \Det(A)y=1 \big\}$$ of affine $(N^2\smash{+}1)$-space. The radical ideal of $\GL_N$ is generated by $\Det(A)y-1$. Observe that multiplication and inversion are $k$-polynomial maps of the $a_{ij}$'s and $y$. Moreover, the determinant operator is a $k$-polynomial function of the $a_{ij}$'s. So $\GL_N$ is defined over $k$ and $\Gamma$ acts on $\GL_N(K)$ by operating on the matrix entries. It is easy to check that this action of $\Gamma$ on $\GL_N(K)$ is by abstract group automorphisms.

\vspace{2mm}\noindent By Proposition $1.10$ of \cite{B}, there exists a $k$-rational representation $\mathcal{R}:G \to \GL_n$ for some integer $n$ (that is, $G$ is $k$-isomorphic to a $k$-subgroup $\mathcal{R}(G)$ of $\GL_n$). 
Then we are done since $\Gamma$ stabilises $\mathcal{R}(G)$ and $\mathcal{R}$ is $\Gamma$-equivariant.
\end{proof}
\end{Lemma}

\noindent Note that the $\star$-action of $\Gamma$ on $G(K)$ is not by automorphisms of algebraic sets (as it stems from a field automorphism, see Theorem $30$ of \cite{St}). 

\vspace{2mm}\noindent In this paper we use the following notation and terminology. For $G$ an algebraic $k$-group, let $G'$ be the derived subgroup of $G$, let $G^{\circ}$ be the connected component of $G$ containing the identity, let $R(G)$ be the radical of $G$ and let $R_u(G)$ be the unipotent radical of $G$. If $R_u(G)$ is trivial then $G$ is \textit{reductive}. If $G$ is connected and $R(G)$ is trivial then $G$ is \textit{semisimple}. If $G$ is semisimple and has no non-trivial connected normal $k$-subgroups then $G$ is \textit{$k$-simple}. We say that $G$ is \textit{absolutely simple} if it is $\overline{k}$-simple. 

\vspace{2mm}\noindent Let $Z$ be a torus that is defined over $k$ (a.k.a. a \textit{$k$-torus}). If $Z$ is $k$-isomorphic to a direct product of finitely many copies of the multiplicative group then $Z$ is \textit{$k$-split}. If $Z$ does not contain a non-trivial $k$-split subtorus then $Z$ is \textit{$k$-anisotropic}. By Proposition $1$ of \cite{Sa}, there exists a unique maximal $k$-split subtorus $Z_s$ of $Z$ and a unique maximal $k$-anisotropic subtorus $Z_a$ of $Z$. Note that $Z=Z_s Z_a$. 
By Theorem $34.3$ of \cite{H}, there exists a finite Galois extension of $k$ over which $Z$ is split. Let $X(Z)$ be the character group of $Z$. For any subtorus $R$ of $Z$, let $X_0(R)$ be the subgroup of $X(Z)$ consisting of all characters that vanish on $R$. 

\vspace{2mm}\noindent Henceforth let $G$ be a connected reductive algebraic $k$-group. 

\vspace{2mm}\noindent Let $S$ be a maximal $k$-split torus of $G$. All maximal $k$-split tori of $G$ are $G(k)$-conjugate to $S$ by Theorem $20.9(ii)$ of \cite{B}. The \textit{$k$-rank} $r_k(G)$ of $G$ is the dimension of $S$. If $r_k(G)=0$ then $G$ is \textit{$k$-anisotropic}. Otherwise, $G$ is \textit{$k$-isotropic}. 
The group $G_a:=C_G(S)'$ is $k$-anisotropic and is called the \textit{(semisimple) anisotropic kernel} of $G$. 

\begin{Theorem}[\textbf{Grothendieck}, \textit{Theorem $34.4$ of \cite{H}}]\label{Grothendieck} There exists a maximal torus of $G$ that is defined over $k$ (a.k.a. a \textit{maximal $k$-torus}). Any $k$-torus of $G$ is contained in a maximal $k$-torus of $G$. 
\end{Theorem} 

\noindent By Theorem \ref{Grothendieck}, there exists a maximal $k$-torus $T$ of $G$ that contains $S$. If $T=S$ then $G$ is \textit{$k$-split}. Let $\rho :X(T) \to X(S)$ be the restriction homomorphism sending $\chi \mapsto \chi|_S$. 
Let $\Phi$ be the root system of $G$ with respect to $T$ and let $\Phi_0:=\Phi \cap X_0(S)$. Let $W:=N_G(T)/T$ and $W_0:=N_{C_G(S)}(T)/T$.
The character space $X(T)_{\R}$ is the real vector space obtained by applying the functor $\otimes _{\Z}\hspace{0.5mm} \R$ to the $\Z$-module $X(T)$. 

\vspace{2mm}\noindent We will need the following useful results.

\begin{Proposition}[\textit{Proposition $20.4$ of \cite{B}}]\label{kparabolicspairing2} Let $S_0$ be a $k$-split torus of $G$ that is not contained in ${Z(G)}^{\circ}$. Then $C_G(S_0)$ is a Levi $k$-subgroup of some parabolic $k$-subgroup of $G$.
\end{Proposition}

\begin{Lemma}\label{mybabyfixed} Let $S_0$ be a $k$-split subtorus of $S$. Let $\Psi$ be the subset of $\Phi$ that vanishes on $S_0$. Then $\Psi$ is the root system of $C_G(S_0)$ with respect to $T$.
\begin{proof} By Proposition \ref{kparabolicspairing2}, $C_G(S_0)=\langle T,U_{\alpha} \hspace{-0.5mm}\mid\hspace{-0.5mm} \alpha \in \Sigma \rangle$ for some parabolic subsystem $\Sigma$ of $\Phi$. Let $\alpha \in \Phi$. Observe that $U_{\alpha} \subseteq C_G(S_0)$ if and only if $S_0 \subseteq \ker(\alpha)$. 
Hence $\Sigma=\Psi$.
\end{proof}
\end{Lemma}

\noindent In particular, Lemma \ref{mybabyfixed} shows that $\Phi_0$ is the root system of $C_G(S)$ with respect to $T$.

\begin{Remark}\label{niceactionremark} We note that all results in Section \ref{Boreltits} hold as long as $S=T_s$. That is, we can relax the assumption that $S$ is maximal among $k$-split tori of $G$.
\end{Remark}

\subsection{The induced action of $\Gamma$ on $X(T)_{\R}$}\label{*action}

\noindent In this section we define and study an induced action of $\Gamma$ on the character space $X(T)_{\R}$. 
We then endow $X(T)_{\R}$ with an inner product that is invariant under the actions of $W$ and $\Gamma$.

\vspace{2mm}\noindent Recall that $\star$ denotes the natural action of $\Gamma$ on $G(K)$. Since $T$ is defined over $k$, this induces an action of $\Gamma$ on $X(T)$ given by $(\chi^{\sigma}) (t):= \sigma \big(\chi (\sigma^{-1} \star t)\big)$ for $\sigma \in \Gamma$, $\chi \in X(T)$ and $t \in T(K)$. It follows easily from Lemma \ref{groupautaction} that this action is by $\Z$-module automorphisms. 
So, by linearly extending, we have a $\R$-linear action $\iota:\Gamma \to \GL\!\big(X(T)_{\R}\big)$, which we refer to as the \textit{induced action of $\Gamma$ on $X(T)_{\R}$}. 

\vspace{2mm}\noindent The absolute Galois group $\Gamma$ is a profinite group, which is finite if and only if $k$ is algebraically closed or is a real closed field (see Corollary VIII.$9.2$ of \cite{L}). However, the following holds.

\begin{Lemma}\label{finiteaction} $\iota(\Gamma)$ is a finite group.
\begin{proof} Recall from Theorem $34.3$ of \cite{H} that there exists a finite Galois extension $L$ of $k$ over which $T$ is split. Let $H=\Gal(K/L)$ and let $\pi:\Gamma \to \Gamma / H$ be the natural projection. Observe that $\iota$ factors through $\pi$, and hence $\iota(\Gamma)$ is finite.
\end{proof}
\end{Lemma}




\noindent We next show that $\iota(\Gamma)$ is well-behaved with respect to the root system $\Phi$.

\begin{Lemma}\label{niceaction} $\Phi$ and $\Phi_0$ are both $\iota(\Gamma)$-stable.
\begin{proof} Let $\sigma \in \Gamma$ and let $\mathfrak{g}$ be the Lie algebra of $G$. By Theorem $3.4$ of \cite{B}, $\mathfrak{g}$ is defined over $k$ (that is, there exists a Lie $k$-algebra $\mathfrak{g}_k \subset \mathfrak{g}$ such that the inclusion map induces an isomorphism $\overline{k} \otimes_k \mathfrak{g}_k \xrightarrow{\smash{\sim}}\mathfrak{g}$). Then $\GL(\mathfrak{g})$ is defined over $k$ as an algebraic group (see Example $1.6(8)$ of \cite{B}). By AG.$\S 14.1$ of \cite{B}, $\Gamma$ acts (semi-linearly) on $K \otimes_k \mathfrak{g}_k =:\mathfrak{g}(K)$ via the first factor with fixed point subset $\mathfrak{g}_k$. This induces an action $\circ$ of $\Gamma$ on $\GL(\mathfrak{g},K)$ given by $(\sigma \circ \phi)(x) = \sigma \big(\phi(\sigma^{-1}(x))\big)$ for $x \in \mathfrak{g}(K)$ and $\phi \in \GL(\mathfrak{g},K)$. The adjoint representation $\Ad:T \to \GL(\mathfrak{g})$ is defined over $k$ by $\S 3.13$ of \cite{B}.

\vspace{2mm}\noindent Let $\mathfrak{g}_{\chi} := \big\{v \in \mathfrak{g} \hspace{0.5mm} \big| \hspace{-0.5mm} \Ad(t)v= \chi(t)v \textnormal{ for all } t \in T \big\}$ for some $\chi \in X(T)$. Note that $\mathfrak{g}_{\chi}$ is defined over $K$ since $\chi$ is defined over $K$. Let $v \in \mathfrak{g}_{\chi}(K)$ and $t \in T(K)$. Denote $w:=\sigma(v)$ and $s:=\sigma^{-1} \star t$. Since $\Ad$ is $\Gamma$-equivariant and $\Gamma$ acts semi-linearly on $\mathfrak{g}$, we have $$\Ad(t)w = {}^{\sigma}\hspace{-1.5mm}\Ad(t)w =\big(\sigma \circ \hspace{-0.2mm}\Ad(s)\big)w= \sigma\big(\hspace{-0.2mm}\!\Ad(s)v\big) = \sigma\big(\chi(s)v\big) = \sigma \big(\chi(s)\big)w={}^{\sigma}\hspace{-0.5mm}\chi(t)w.$$ It follows that $\sigma \big(\mathfrak{g}_{\chi}(K)\big) \subseteq \mathfrak{g}_{{}^\sigma\hspace{-0.5mm}\chi}(K)$. Replacing $\sigma$ with $\sigma^{-1}$ gives us $\sigma\big(\mathfrak{g}_{\chi}(K)\big) =\mathfrak{g}_{{}^\sigma\hspace{-0.5mm}\chi}(K)$. Recall that $\Phi$ is defined to be the set of non-trivial characters $\chi \in X(T)$ that satisfy $\mathfrak{g}_{\chi} \neq 0$. If $\mathfrak{g}_{\chi}(K)=0$ then $\mathfrak{g}_{\chi}=0$. Hence $\Phi$ is $\iota(\Gamma)$-stable. 

\vspace{2mm}\noindent Finally, observe that $X_0(S)$ is a $\Gamma$-invariant subgroup of $X(T)$ (since $S$ is defined over $k$). So $\Phi_0 := \Phi \cap X_0(S)$ is also $\Gamma$-invariant.
\end{proof}
\end{Lemma}

\begin{Corollary}\label{downtown} Let $H$ be a connected reductive $K$-subgroup of $G$ that contains $T$. Let $\Sigma$ be the root system of $H$ with respect to $T$. Then $H$ is defined over $k$ if and only if $\Sigma$ is $\iota(\Gamma)$-stable.
\begin{proof} Let $\alpha \in \Phi$ and let $U_{\alpha}$ be the associated $T$-root subgroup of $G$. By Theorem $34.4(b)$ of \cite{H}, $U_{\alpha}$ is defined over $K$ and the associated isomorphism $u_{\alpha}: \smash{\overline{k}^+} \to U_{\alpha}$ is a $K$-isomorphism. Then, by $\S 2.4(6)$ of \cite{Sa}, for any $\sigma \in \Gamma$ we have \begin{equation}\label{littleaction} \sigma \star U_{\alpha} = U_{\iota(\sigma)(\alpha)}. \end{equation} 
Assume that $\Sigma$ is $\iota(\Gamma)$-stable. Recall that $H=\langle T, U_{\alpha} \hspace{0.5mm}|\hspace{0.5mm} \alpha \in \Sigma \rangle$. Then, using $(\ref{littleaction})$, we see that $H(K)$ is stabilised by the $\star$-action of $\Gamma$. That is, $H$ is defined over $k$. The converse follows immediately from Lemma \ref{niceaction}.
\end{proof}
\end{Corollary}

\noindent Consider the action of $\Gamma$ on $\GL\!\big(X(T)_{\R}\big)$ given by $\psi^\sigma:= \iota(\sigma)\psi \iota(\sigma)^{-1}$ for $\sigma \in \Gamma$ and $\psi \in \GL\!\big(X(T)_{\R}\big)$. In a slight abuse of notation, we also refer to this action as $\iota$. Let $\epsilon:W \hookrightarrow \GL\!\big(X(T)_{\R}\big)$ be the natural embedding. It follows from Lemma \ref{niceaction} that $\epsilon(W)$ is $\iota(\Gamma)$-stable (since a Weyl group is generated by reflections about its roots).

\begin{Lemma}\label{Weylgpstab} The natural embedding $\epsilon:W \hookrightarrow \GL\!\big(X(T)_{\R}\big)$ is $\Gamma$-equivariant.
\begin{proof} Since $T$ and $N_G(T)$ are both defined over $k$, the action $\star$ of $\Gamma$ on $G(K)$ descends to an action $\star$ of $\Gamma$ on the quotient $W=N_G(T)(K)\big/T(K)$. Let $\sigma \in \Gamma$, $\chi \in X(T)_{\R}$ and $w \in W$. Recall from $\S 8.1$ of \cite{MT} that $W$ acts on $X(T)$ by $w(\chi) =\chi \cdot w$. We check that $$(\epsilon(w)^\sigma)(\chi) =\iota(\sigma) \epsilon(w) \iota(\sigma)^{-1} (\chi)= \sigma\sigma^{-1}\chi\sigma w \sigma^{-1} = (\sigma \star w)(\chi).$$ So the natural embedding of $W$ in $\GL\!\big(X(T)_{\R}\big)$ is indeed $\Gamma$-equivariant.
\end{proof}
\end{Lemma}

\noindent Since $T$ is defined over $k$, we can define the \textit{extended Weyl group} $\tilde{W}:=(N_G(T)(K) \rtimes_{\star} \Gamma)/T(K)$ of $G$ with respect to $T$. The natural action of $\tilde{W}$ on $T(K)$ is faithful, and it induces a faithful action of $\tilde{W}$ on $X(T)$ given by $w(\chi) =\chi \cdot w$ for $\chi \in X(T)$ and $w \in \tilde{W}$. Extending linearly, we have a natural embedding $\tilde{\epsilon}:\tilde{W} \hookrightarrow \GL\!\big(X(T)_{\R}\big)$ such that $\tilde{\epsilon}|_{W}=\epsilon$. It follows from Lemma \ref{Weylgpstab} that $\tilde{\epsilon}(\smash{\tilde{W}})=\epsilon(W) \rtimes \iota(\Gamma)$, where $\Gamma$ acts on $\epsilon(W)$ by sending a reflection $w_{\alpha} \mapsto w_{\iota(\sigma)(\alpha)}$ for $\alpha \in \Phi$ and $\sigma \in \Gamma$. Note that $\smash{\tilde{W}}$ is finite by Lemma \ref{finiteaction}. Henceforth, as is the convention with Weyl groups, we consider the embedding $\tilde{\epsilon}$ to be implicit.

\vspace{2mm}\noindent In summary, we have a finite group $\smash{\tilde{W}}=W \rtimes \iota(\Gamma)$ that acts faithfully on $T(K)$ and $X(T)_{\R}$. So there exists a $\smash{\tilde{W}}$-invariant inner product $(\cdot \hspace{0.5mm},\cdot)$ on $X(T)_{\R}$, which is unique up to scalar multiplication on each irreducible $\smash{\tilde{W}}$-submodule of $X(T)_{\R}$. 

\subsection{The subspace structure of $X(T)_{\R}$}\label{subspstruct}

\noindent In this section we relate the subspace structure of $X(T)_{\R}$ to the subtori structure of $T$.

\vspace{2mm}\noindent For any subtorus $Z$ of $T$, recall that $X_0(Z)_{\R}$ is defined to be the subspace of $\R$-linear combinations
of characters of $T$ that vanish on $Z$.

\begin{Lemma}\label{abcdp} Let $T_1$ and $T_2$ be subtori of $T$. Then 

\vspace{2mm}\noindent $(i)$ $X_0(T_1T_2)_{\R} = X_0(T_1)_{\R} \cap X_0(T_2)_{\R}$, and 

\vspace{1mm}\noindent $(ii)$ $X_0((T_1 \cap T_2)^{\circ})_{\R} = X_0(T_1)_{\R} +X_0(T_2)_{\R}$.
\begin{proof} $(i)$ Let $u \in X_0(T_1T_2)_{\R}$. It is immediate that $u \in X_0(T_1)_{\R} \cap X_0(T_2)_{\R}$ since $T_1, T_2 \subseteq T_1T_2$. Conversely, let $v \in X_0(T_1)_{\R} \cap X_0(T_2)_{\R}$ and let $t=t_1t_2$ where $t_1 \in T_1$ and $t_2 \in T_2$. Write $v=\sum_j\lambda_jv_j$ for some $\lambda_j \in \R$ and $v_j \in X(T)$. Observe that $v_j(t)=v_j(t_1)+v_j(t_2)=0$ for each $j$, and so $v \in X_0(T_1T_2)_{\R}$.

\vspace{2mm}\noindent $(ii)$ We first show that the dimensions are equal. For any subtorus $Z$ of $T$, observe that $\dim X_0(Z)_{\R} = \dim T - \dim Z$. Moreover, $\dim T_1T_2 = \dim T_1 + \dim T_2 - \dim (T_1 \cap T_2)^{\circ}$. Combining this with part $(i)$, we have $\dim X_0((T_1 \cap T_2)^{\circ})_{\R} = \dim \!\big(X_0(T_1)_{\R} + X_0(T_2)_{\R}\big)$.

\vspace{2mm}\noindent Now let $x =x_1+x_2$ where $x_1 \in X_0(T_1)_{\R}$ and $x_2 \in X_0(T_2)_{\R}$. Then certainly $x$ is a $\R$-linear combination of characters that vanish on $(T_1 \cap T_2)^{\circ}$. That is, $X_0((T_1 \cap T_2)^{\circ})_{\R} \supseteq X_0(T_1)_{\R} +X_0(T_2)_{\R}$. This completes the proof.
\end{proof}
\end{Lemma}

\noindent We need the following elementary result. For $V$ a vector space, $v \in V$ and $H < \GL(V)$, let $V^H$ be the fixed point subspace of $V$ under $H$ and let $H(v)$ be the $H$-orbit of $v$ in $V$.

\begin{Lemma}\label{easyyyy} Let $V$ be an finite-dimensional inner product space. Let $H$ be a finite subgroup of the isometry group of $V$. Let $C$ be a $H$-submodule of $V$ such that $V=V^H \oplus C$. Then $C=\{v \in V \hspace{0.5mm} |\hspace{0.5mm} \sum_{u \in H(v)} u =0 \}=(V^H)^{\perp}$.
\begin{proof} Let $c \in C$ and $v \in V^H$. Observe that $\sum_{d \in H(c)}d \in C \cap V^H=\{0 \}$. Then $0=(\sum_{d \in H(c)}d,v)=|H(c)|(c,v)$ and so $c \in (V^H)^{\perp}$. Conversely, let $x \in (V^H)^{\perp}$. Then $\sum_{y \in H(x)} y \in (V^H)^{\perp} \cap V^H =\{0\}$. We are done as $C$ and $(V^H)^{\perp}$ are both complements of $V^H$ in $V$.
\end{proof}
\end{Lemma}

\noindent As in $\S \ref{*action}$, let $\iota:\Gamma \to \GL\!\big(X(T)_{\R}\big)$ be the induced action, let $\smash{\tilde{W}}$ be the extended Weyl group of $G$ with respect to $T$ and let $(\cdot \hspace{0.5mm},\cdot)$ be a $\smash{\tilde{W}}$-invariant inner product on $X(T)_{\R}$.

\vspace{2mm}\noindent A subtorus $R$ of $T$ is \textit{standard} if $R=T^Z$ for some subgroup $Z$ of $\smash{\tilde{W}}$. For example, if $R={Z\big(C_G(R)\big)}_s^{\circ}$ then $R$ is standard. 
Similarly, a subgroup $\Omega$ of $\smash{\tilde{W}}$ is \textit{standard} if $\Omega=\Fix_{\smash{\tilde{W}}}(U)$ for some subtorus $U$ of $T$. 

\vspace{2mm}\noindent Let $U$ be a subtorus of $T$. As discussed in $\S \ref{*action}$, the natural action of $\smash{\tilde{W}}$ on $T(K)$ induces an action of $\smash{\tilde{W}}$ on $X(T)$ and hence on $X(T)_{\R}$. Similarly, the action of $\Stab_{\smash{\tilde{W}}}(U)$ on $U(K)$ induces an action of $\Stab_{\smash{\tilde{W}}}(U)$ on $X(U)_{\R}$. A linear embedding $X(U)_{\R} \hookrightarrow X(T)_{\R}$ is a \textit{$\smash{\tilde{W}}$-embedding} if it is $\Stab_{\smash{\tilde{W}}}(U)$-equivariant.



\begin{Lemma}\label{niceactioncor} Let $U$ be a standard subtorus of $T$. There exists a $\smash{\tilde{W}}$-embedding $X(U)_{\R} \hookrightarrow X(T)_{\R}$, where the image of $X(U)_{\R}$ is the orthogonal complement of $X_0(U)_{\R}$ in $X(T)_{\R}$. Any other $\smash{\tilde{W}}$-embedding $X(U)_{\R} \hookrightarrow X(T)_{\R}$ has the same image. In particular, if $U=S$ then the image of $X(U)_{\R}$ is the fixed point subspace of $X(T)_{\R}$ under $\iota(\Gamma)$.
\begin{proof} 

\vspace{2mm}\noindent Consider the short exact sequence of $\Stab_{\smash{\tilde{W}}}(U)$-modules \begin{equation}\label{ses1} 0 \to X_0(U) \to X(T) \to X(U) \to 0.\end{equation} The functor $\otimes_{\Z}\hspace{0.6mm} \R$ is exact (since $\R$ is torsion-free as a $\Z$-module). So applying $\otimes_{\Z}\hspace{0.6mm} \R$ to the short exact sequence $(\ref{ses1})$ gives us a direct sum of $\Stab_{\smash{\tilde{W}}}(U)$-modules \begin{equation}\label{ses2} X(T)_{\R} = X_0(U)_{\R} \oplus X(U)_{\R}.\end{equation} 
In other words, there exists a $\smash{\tilde{W}}$-embedding $e:X(U)_{\R} \hookrightarrow X(T)_{\R}$ such that $e(X(U)_{\R})$ is a linear complement of $X_0(U)_{\R}$. Since $U$ is standard, we have $T^{\Fix_{\smash{\tilde{W}}}(U)}=U$. It follows that $e(X(U)_{\R})$ is the fixed point subspace of $X(T)_{\R}$ under $\Fix_{\smash{\tilde{W}}}(U)$. Then, by Lemma \ref{easyyyy}, the decomposition in $(\ref{ses2})$ is orthogonal and $X_0(U)_{\R}=\{v \in X(T)_{\R} \hspace{0.5mm} |\hspace{0.5mm} \sum_{\sigma \in \Fix_{\smash{\tilde{W}}}(U)} v^{\sigma} =0 \}$. It is easy to see that any $\smash{\tilde{W}}$-embedding $X(U)_{\R} \hookrightarrow X(T)_{\R}$ must arise in this way. 
This proves the first and second assertions.

\vspace{2mm}\noindent Now consider the case where $U=S$. Certainly $S$ is standard. Let $\smash{\tilde{W}_S}$ be the extended Weyl group of $C_G(S)$ with respect to $T$. Under the natural embedding of $\smash{\tilde{W}_S}$ in $\smash{\tilde{W}}$, observe that $\Fix_{\smash{\tilde{W}}}(S)=\smash{\tilde{W}_S} =W_0 \rtimes \iota(\Gamma)$. By Lemma \ref{mybabyfixed}, $W_0$ is the Weyl group of $\Phi_0$. Hence, by the first assertion, $W_0$ acts trivially on $X(S)_{\R}$. This completes the proof.
\end{proof}
\end{Lemma}

\noindent Henceforth, we will take this identification in Lemma \ref{niceactioncor} to be implicit.
\begin{Corollary}\label{niceactioncor1} Let $T_1$ and $T_2$ be standard $k$-subtori of $T$. Under the identification in Lemma \ref{niceactioncor}, we have 

\vspace{2mm}\noindent $(i)$ $X(T_1T_2)_{\R} = X(T_1)_{\R} + X(T_2)_{\R}$, and 

\vspace{1mm}\noindent $(ii)$ $X((T_1 \cap T_2)^{\circ})_{\R} = X(T_1)_{\R} \cap X(T_2)_{\R}$.
\begin{proof} Taking the orthogonal complement of both sides, under the identification in Lemma \ref{niceactioncor}, we see that $(i)$ (resp. $(ii)$) is equivalent to statement $(i)$ (resp. $(ii)$) of Lemma \ref{abcdp}.
\end{proof}
\end{Corollary}

\subsection{A $\Gamma$-order on $X(T)$}\label{Gammabase}

\noindent In this section we construct a total order on $X(T)$ that is compatible with the induced action $\iota:\Gamma \to \GL\!\big(X(T)_{\R}\big)$. 

\vspace{2mm}\noindent There are some subtleties with this construction, as we illustrate with the following lemma. Note that a total order on $X(T)$ uniquely determines a compatible base of $\Phi$. 
We say that $G$ is \textit{$k$-quasisplit} if it contains a Borel subgroup that is defined over $k$. 

\begin{Lemma}\label{quasiiiiiii} There exists a base of $\Phi$ that is $\iota(\Gamma)$-stable if and only if $G$ is $k$-quasisplit.
\begin{proof} Let $\Delta$ be a base of $\Phi$ and consider the associated Borel subgroup $B_{\Delta}:= \big\langle T, U_{\alpha} \hspace{0.5mm}\big|\hspace{0.5mm} \alpha \in \Delta \big\rangle$ of $G$. By Corollary $7.5$ (and $\S 34.2$) of \cite{H}, $B_{\Delta}$ is defined over $K$ since it is generated by connected $K$-groups. It follows from Equation $(\ref{littleaction})$ and $\S 27.3$ of \cite{H} that the map $\Delta \mapsto B_{\Delta}$ is a $\Gamma$-equivariant bijection from the set of bases of $\Phi$ to the set of Borel $K$-subgroups of $G$ that contain $T$. Taking the $\Gamma$-fixed points of this map gives us a bijection from the set of $\iota(\Gamma)$-stable bases of $\Phi$ to the set of Borel $k$-subgroups of $G$ that contain $T$.

\vspace{2mm}\noindent Now let $G$ be $k$-quasisplit and let $B$ be a Borel $k$-subgroup of $G$. Then $B$ contains a maximal $k$-split torus $S_0$ of $G$ and we can decompose $B=R_u(B) \rtimes C_G(S_0)$. Since all maximal $k$-split tori of $G$ are $G(k)$-conjugate, there exists $g \in G(k)$ satisfying $gS_0g^{-1}=S$. Then $gBg^{-1}$ is a Borel $k$-subgroup of $G$ that contains $T$. 
\end{proof}
\end{Lemma}


\noindent In general, however, $G$ is not necessarily $k$-quasisplit. So we proceed as follows to find an appropriate order for $X(T)$. 

\vspace{2mm}\noindent Let $\mathcal{T}:= \{1 \subset T_1 \subset T_2 \subset ... \subset T_r = T\}$ be a chain of $K$-subtori of $T$. A \textit{$\mathcal{T}$-order} on $X(T)$ is a translation-invariant total order $<$ that satisfies the following condition:

\vspace{2mm}\noindent $(*)$ For every $j=1,...,r$ and every positive $x \in X(T)$ that does not vanish on $T_j$, then the restriction of $x$ to $T_j$ is also positive.

\begin{Lemma}\label{funkycompatibleordering} Let $\mathcal{T}$ be a chain of $K$-subtori of $T$. There exists a $\mathcal{T}$-order on $X(T)$.
\begin{proof} Write $\mathcal{T}=\{1 \subset T_1 \subset T_2 \subset ... \subset T_r = T \}$. For each $j=1,...,r$, let $\rho_j:X(T) \to X(T_j)$ be the restriction map. Let $\mathcal{B}=\{b_1,...,b_n\}$ be an ordered basis of $X(T)$ where there exist integers $1 \leq i_1 < i_2 < ... < i_r = n$ such that $\{\rho_j(b_1),...,\rho_j(b_{i_j})\}$ generates $X(T_j)$ for each $j =1,...,r$. Certainly, such a basis exists.

\vspace{2mm}\noindent Let $<$ be the lexicographic order on $X(T)$ with respect to $\mathcal{B}$. It is easy to see that $<$ is a total order on $X(T)$ that is translation-invariant. Let $j \in \{1,...,r\}$ and let $x,y \in X(T)$. By construction, if $x > y$ and $\rho_j(x) \neq \rho_j(y)$ then $\rho_j(x) > \rho_j(y)$. In particular, $<$ satisfies condition $(*)$.
\end{proof}
\end{Lemma}

\noindent A \textit{$\Gamma$-order} on $X(T)$ is a translation-invariant total order $<$ such that the set of positive elements of $X(T) \setminus X_0(S)$ is $\iota(\Gamma)$-stable. The following consequence of Lemma \ref{funkycompatibleordering} is well-known, see for instance $\S 2$ of \cite{Sa} or $\S 21.8$ of \cite{B}.

\begin{Corollary}\label{funkycompatibleorderingcor} There exists a $\Gamma$-order on $X(T)$.
\begin{proof} Let $\mathcal{T}$ be the chain $1 \subset S \subseteq T$. By Lemma \ref{funkycompatibleordering}, there exists a $\mathcal{T}$-order $<$ on $X(T)$. Since $\iota(\Gamma)$ descends to the trivial action on $X(S)$, the set of positive elements of $X(T) \setminus X_0(S)$ is $\iota(\Gamma)$-stable. 
\end{proof}
\end{Corollary}

\noindent A base of $\Phi$ that is compatible with a $\Gamma$-order on $X(T)$ is called a \textit{$\Gamma$-system of simple roots} (a.k.a. a \textit{$\Gamma$-base}) of $\Phi$. Recall that $W_0:=N_{C_G(S)}(T)/T$. Let $W_{\Gamma}:=N_{N_G(S)}(T)/T$.

\begin{Lemma}[\textit{$\S 3$ of \cite{Sa}}]\label{propertiesofgammasystem} Let $\Delta$ be a $\Gamma$-base of $\Phi$. Then

\vspace{2mm}\noindent $(i)$ $\Delta_0:=\Delta \cap X_0(S)$ is a base of $\Phi_0$, 

\vspace{1mm}\noindent $(ii)$ $W_0$ is the Weyl group of $\Phi_0$ and $W_{\Gamma}=\big\{w \in W \hspace{0.5mm}\big|\hspace{0.5mm} w\big(X_0(S)\big)=X_0(S)\big\}$,

\vspace{1mm}\noindent $(iii)$ for each $\sigma \in \Gamma$ there exists a unique $w_{\sigma} \in W_0$ such that $w_{\sigma}\iota(\sigma)(\Delta)=\Delta$, and

\vspace{1mm}\noindent $(iv)$  $W_{\Gamma}$ acts simply transitively on the set of all $\Gamma$-bases of $\Phi$.
\end{Lemma}

\subsection{The index of $G$}\label{index}

\noindent In this section we define and investigate an important invariant of $G$ called its index. We also present Tits' partial classification of connected reductive $k$-groups in terms of their indices.

\vspace{2mm}\noindent Let $\smash{\tilde{W}}$ be the extended Weyl group of $G$ with respect to $T$ and let $(\cdot \hspace{0.5mm},\cdot)$ be a $\smash{\tilde{W}}$-invariant inner product on $X(T)_{\R}$. Let $<$ be a $\Gamma$-order on $X(T)$. Let $\Delta$ be the system of simple roots for $G$ with respect to $T$ that is compatible with $<$. Let $\Delta_0$ be the subset of $\Delta$ that vanishes on $S$. 
Let $\iota:\Gamma \to \GL\!\big(X(T)_{\R}\big)$ be the induced action.

\vspace{2mm}\noindent Let $\sigma \in \Gamma$. Recall from Lemma \ref{niceaction} that $\Phi$ is $\iota(\sigma)$-stable, and so $\iota(\sigma)(\Delta)$ is another choice of base of $\Phi$. Since $W$ acts simply transitively on the set of bases of $\Phi$, there exists a unique $w_{\sigma} \in W$ such that $w_{\sigma}\iota(\sigma)(\Delta)=\Delta$. In fact $w_{\sigma} \in W_0$ by Lemma \ref{propertiesofgammasystem}$(iii)$. The \textit{Tits action} $\hat{\iota}:\Gamma \to \GL\!\big(X(T)_{\R}\big)$ is defined by $\hat{\iota}(\sigma):=w_{\sigma} \iota(\sigma)$ (it is indeed an action, see $\S 2.3$ of \cite{T}). 

\begin{Lemma}\label{Titsactionproperties} The Tits action $\hat{\iota}$ of $\Gamma$ on $X(T)_{\R}$ stabilises $\Delta$ and $\Delta_0$ and is by linear isometries. Moreover, we identify $X(S)_{\R}$ with the largest subspace of $X(T)_{\R}$ that is fixed pointwise by $\hat{\iota}(\Gamma)$ and is perpendicular to all roots in $\Delta_0$.
\begin{proof} The first statement is immediate from Lemmas \ref{niceaction}, \ref{propertiesofgammasystem}$(ii)$ and by choice of $(\cdot \hspace{0.5mm},\cdot)$. Again by Lemma \ref{propertiesofgammasystem}$(ii)$, the actions $\iota$ and $\hat{\iota}$ of $\Gamma$ coincide on the subspace of $X(T)_{\R}$ that is perpendicular to all roots in $\Delta_0$. 
The second statement then follows from Lemma \ref{niceactioncor}.
\end{proof}
\end{Lemma}

\noindent If $\hat{\iota}(\Gamma)$ is trivial then $G$ is of \textit{inner type}. Otherwise, $G$ is of \textit{outer type}. Note that $\iota(\Gamma) \leq W$ if and only if $G$ is of inner type. 

\vspace{2mm}\noindent The \textit{index} $\mathcal{I}(G)$ of $G$ (with respect to $T$ and $<$) is the quadruple $\big(X(T)_{\R}, \Delta , \Delta_0, \hat{\iota}(\Gamma)\big)$. Note that our definition contains slightly more information than Tits' original definition in \cite{T}.

\vspace{2mm}\noindent If $G$ is semisimple then we illustrate $\mathcal{I}(G)$ using a \textit{Tits-Satake diagram}, which is constructed by taking the Dynkin diagram of $G$, blackening each vertex in $\Delta_0$ and linking with a solid gray bar (or somehow indicating) all of the $\hat{\iota}(\Gamma)$-orbits of $\Delta$. 
For arbitrary connected reductive $G$, we illustrate $\mathcal{I}(G)$ using the notation $\mathcal{I}(G') \times \mathcal{T}^d_c$, where $c=\dim\hspace{-0.7mm}\big({Z(G)}_s^{\circ}\big)$, $d=\dim\hspace{-0.7mm}\big({Z(G)}^{\circ}\big)$ and $\mathcal{I}(G')$ is represented by its Tits-Satake diagram. 
In the rare cases where it is not immediately obvious from this illustration, we also explicitly describe the subgroup $\hat{\iota}(\Gamma)$ of $\GL\!\big(X(T)_{\R}\big)$.


\begin{Example} Let $k=\F_q$ for $q$ a prime power and let $G(k)=\GU_3(q)$. Then $r(G)=3$, $r_k(G)=1$, $\Delta=A_2$, $\Delta_0 = \varnothing$ and $\hat{\iota}(\Gamma) \cong \Z_2$. So we write $\mathcal{I}(G)=\begin{tikzpicture}[baseline=-0.3ex]\dynkin[ply=2,rotate=90,fold radius = 2mm]{A}{oo}\end{tikzpicture} \times \mathcal{T}^1_0$. 
\end{Example} 


\noindent The following theorem reduces the problem of classifying connected reductive algebraic $k$-groups to that of classifying their $k$-anisotropic kernels. However, classifying $k$-anisotropic semisimple $k$-groups in full generality is a very difficult problem.

\begin{Theorem}[\textbf{Tits}, \textit{Theorem $2$ of \cite{T}}]\label{Tits} A connected reductive algebraic $k$-group is uniquely determined up to $k$-isomorphism by its $K$-isomorphism class, its index and the $k$-isomorphism class of its anisotropic kernel.
\end{Theorem}

\noindent Inspired by the notion of an abstract root system, in $\S \ref{abstractindex}$ we give a purely combinatorial definition of an index. We will see in Theorem \ref{indeximpo} that, as a combinatorial object, the index of $G$ is independent of the choice of $T$ and $<$ (this is essentially Proposition $2$ of \cite{Sa0}).

\subsection{Orientation maps}\label{orientationmaps}

\noindent In this section we take two maximal $k$-tori of $G$ and we construct a bijective isometry between their respective character spaces. We call this an orientation map, and we investigate some of its properties.

\vspace{2mm}\noindent Recall that $T$ is a maximal $k$-torus of $G$. Let $T_1$ be another maximal $k$-torus of $G$. By Theorem $20.9(ii)$ of \cite{B}, there exists $g \in G(K)$ such that $(T_1)^g=T$. The induced map $X(T_1) \to X(T)$ given by $\chi \mapsto \chi \cdot g^{-1}$ is a $\Z$-module isomorphism. 
This extends uniquely to a $\R$-vector space isomorphism $\theta_g:X(T_1)_{\R} \to X(T)_{\R}$. We call $\theta_g$ the \textit{orientation map induced by $g$}. 

\vspace{2mm}\noindent Recall that $\Phi$ is the root system of $G$ with respect to $T$. Let $\Phi_1$ be the root system of $G$ with respect to $T_1$. The following result shows that the map $g \mapsto \theta_g$ is is compatible with the respective root systems.

\begin{Lemma}\label{propertiesinduced} Let $\alpha \in \Phi$ and let $U_{\alpha} <G$ be the associated root $K$-group with respect to $T$. Then $\theta_g^{-1}(\Phi)=\Phi_1$ and $(U_{\alpha})^{g^{-1}}$ is the root $K$-group associated to $\theta_g^{-1}(\alpha)$ with respect to $T_1$.
\begin{proof} Let $u_{\alpha}:K^+ \to U_{\alpha}$ be the root $K$-isomorphism associated to $\alpha$ with respect to $T$ (it is not necessarily defined over $k$). Let $t \in T$ and $c \in K^+$. Recall from Theorem $8.17(c)$ of \cite{MT} that, up to scalar multiplication, $u_{\alpha}$ is the unique $K$-morphism from $K^+$ to $G$ that is a $K$-isomorphism onto its image and that satisfies the property $u_{\alpha}(c)^t=u_{\alpha}(\alpha(t)c)$. Let $t_1:=t^{g^{-1}}$ and $\beta :=\theta_g^{-1}(\alpha)$. Observe that $(U_{\alpha})^{g^{-1}}$ is a $1$-dimensional unipotent $K$-subgroup of $G$ that is normalised by $T_1$. We check that $$\big(u_{\alpha}(c)^{g^{-1}}\big)^{t_1}=\big(u_{\alpha}(c)^{t}\big)^{g^{-1}}=\big(u_{\alpha}(\alpha(t)c)\big)^{g^{-1}}=\big(u_{\alpha}(\beta(t_1)c)\big)^{g^{-1}}$$ 
and hence $\beta \in \Phi_1$ and $(U_{\alpha})^{g^{-1}}$ is the associated root $K$-group with respect to $T_1$. 
Conversely, if $\alpha_1 \in \Phi_1$ then $\theta_g(\alpha_1) \in \Phi$ by the same argument in reverse.
\end{proof}
\end{Lemma}

\begin{Corollary}\label{propertiesinducedcor} Let $H$ be a connected reductive $k$-subgroup of $G$ that contains $T_1$. Let $\Delta_H$ be a system of simple roots for $H$ with respect to $T_1$. Then $\theta_g(\Delta_H)$ is a system of simple roots for $H^g$ with respect to $T$. 
\begin{proof} Let $\Psi:=\langle \Delta_H\rangle$. Let $V$ denote a $T_1$-root group and $U$ denote a $T$-root group. We have $H=\big\langle T_1, V_{\alpha} \hspace{0.5mm}\big|\hspace{0.5mm} \alpha \in \Psi\big\rangle$. Then $H^g=\big\langle T, U_{\alpha} \hspace{0.5mm}\big|\hspace{0.5mm} \alpha \in \theta_g(\Psi) \big\rangle$ by Lemma \ref{propertiesinduced}. But $\theta_g(\Delta_H)$ is a base of $\theta_g(\Psi)$ since $\theta_g$ is a bijective isometry.
\end{proof}
\end{Corollary} 

\noindent Consider the natural action $\star$ of $\Gamma$ on $G(K)$. Let $\iota:\Gamma \to \GL\!\big(X(T)_{\R}\big)$ and $\iota_1:\Gamma \to \GL\!\big(X(T_1)_{\R}\big)$ be the respective induced actions. Let $T_2$ be another maximal $k$-torus of $G$ and let $h \in G(K)$ such that $(T_2)^h=T_1$. Let $\pi:N_G(T)(K) \to W$ be the natural projection.

\vspace{2mm}\noindent The following result shows that the map $g \mapsto \theta_g$ exhibits functorial-like behaviour.

\begin{Lemma}\label{functorialbehaviour} $\theta_{gh}=\theta_g\theta_h$ and $\theta_g^{-1}=\theta_{g^{-1}}$. If $g \in N_G(T)(K)$ then $\theta_g=\pi(g)^{-1}$. Moreover, $\theta_{\sigma \star g}=\iota(\sigma)\theta_g \iota_1(\sigma)^{-1}$ for every $\sigma \in \Gamma$. 
\begin{proof} Let $\gamma \in X(T_2)$. Then $$\theta_{gh}(\gamma)=\gamma \cdot (gh)^{-1}=(\gamma \cdot h^{-1})\cdot g^{-1}=\theta_h(\gamma) \cdot g^{-1} =\theta_g\theta_h(\gamma).$$ If $h=g^{-1}$ then $\theta_g\theta_h$ is the identity map on $X(T)_{\R}$ and so $\theta_g^{-1}=\theta_{h}$. Now let $\chi \in X(T)$. If $g \in N_G(T)(K)$ (that is, $T_1=T$) then it is immediate that $\theta_g(\chi)=\chi \cdot g^{-1}=\pi(g)^{-1}(\chi)$. 

\vspace{2mm}\noindent Finally, let $\sigma \in \Gamma$ and let $\xi \in X(T_1)$. We check that $$\iota(\sigma)\theta_g \iota_1(\sigma)^{-1}(\xi)= \sigma \sigma^{-1} \xi \sigma g^{-1} \sigma^{-1}=\xi \cdot(\sigma \star g^{-1}) =\theta_{\sigma \star g}(\xi).\qedhere$$
\end{proof}
\end{Lemma}

\noindent Now recall that $W$ is the Weyl group of $G$ with respect to $T$ and $\smash{\tilde{W}}$ is the extended Weyl group of $G$ with respect to $T$. Let $W_1$ be the Weyl group of $G$ with respect to $T_1$ and let $\smash{\tilde{W}}_1$ be the extended Weyl group of $G$ with respect to $T_1$.

\begin{Corollary}\label{propertiesinducedcor1} Under the natural embeddings of $\smash{\tilde{W}}$ in $\GL\!\big(X(T)_{\R}\big)$ and $\smash{\tilde{W}}_1$ in $\GL\!\big(X(T_1)_{\R}\big)$, we have $\theta_gW_1\theta_g^{-1}=W$ and $\theta_g\smash{\tilde{W}}_1\theta_g^{-1}=\smash{\tilde{W}}$.
\begin{proof} The fact that $\theta_gW_1\theta_g^{-1}=W$ is immediate from Lemma \ref{propertiesinduced}, since a Weyl group is generated by reflections about its roots.

\vspace{2mm}\noindent Let $\sigma \in \Gamma$. Observe that $(\sigma \star g)g^{-1} \in N_G(T)$. Denote $w:=\pi\big((\sigma \star g)g^{-1}\big) \in W$. Then applying Lemma \ref{functorialbehaviour} gives us $\theta_g\iota_1(\sigma)\theta_g^{-1}=w\iota(\sigma)$. Recall from $\S \ref{*action}$ that $\smash{\tilde{W}}=W \rtimes \iota(\Gamma)$. Similarly, $\smash{\tilde{W_1}}=W_1 \rtimes \iota_1(\Gamma)$. Combining this with $\theta_gW_1\theta_g^{-1}=W$ tells us that $\theta_g\smash{\tilde{W}}_1\theta_g^{-1} \subseteq \smash{\tilde{W}}$. The same argument in reverse gives us $\theta_g\smash{\tilde{W}}_1\theta_g^{-1} \supseteq \smash{\tilde{W}}$.
\end{proof}
\end{Corollary}

\begin{Corollary}\label{inducedisometry} Let $(\cdot \hspace{0.5mm},\cdot)$ be a $\smash{\tilde{W}}$-invariant inner product on $X(T)_{\R}$. Then there exists a $\smash{\tilde{W}}_1$-invariant inner product on $X(T_1)_{\R}$ such that $\theta_g$ is an isometry.
\begin{proof} Just as we observed for $X(T)_{\R}$, there exists a $\smash{\tilde{W}}_1$-invariant inner product $(\cdot \hspace{0.5mm},\cdot)_1$ on $X(T_1)_{\R}$ that is unique up to scalar multiplication on each irreducible $\smash{\tilde{W}}_1$-submodule of $X(T_1)_{\R}$. Hence, by Corollary \ref{propertiesinducedcor1}, we can choose a normalisation of $(\cdot \hspace{0.5mm},\cdot)_1$ such that $\theta_g$ is an isometry.
\end{proof}
\end{Corollary}

\noindent If $<$ is a total order on $X(T_1)$ and $\chi_1, \chi_2 \in X(T)$, let $<_g$ be the total order on $X(T)$ defined by $\chi_1<_g \chi_2$ if and only if $\chi_1 \cdot g < \chi_2 \cdot g$. 

\begin{Lemma}\label{preservesGammabase} Let $\mathcal{T}$ be a chain of $K$-subtori of $T_1$. Let $<$ be a $\mathcal{T}$-order on $X(T_1)$. Then $<_g$ is a $\mathcal{T}^g$-order on $X(T)$. In particular, if $(S_1)^g = S$ and $<$ is a $\Gamma$-order on $X(T_1)$ then $<_g$ is a $\Gamma$-order on $X(T)$.
\begin{proof} Recall that the map $X(T_1) \to X(T)$ given by $\chi \mapsto \chi \cdot g^{-1}$ is a $\Z$-module isomorphism. The first assertion then follows immediately. Applying the first assertion to the case where $\mathcal{T}=\{1 \subset S_1 \subset T_1\}$ gives us the second assertion.
\end{proof}
\end{Lemma}

\begin{Corollary}\label{gammalinorder} Let $S_1:=(T_1)_s$ and assume that $(S_1)^g \subseteq S$. There exists a $\Gamma$-order $<$ on $X(T_1)$ such that $<_g$ is a $\Gamma$-order on $X(T)$.
\begin{proof}Let $\mathcal{T}$ be the chain of $K$-tori $1 \subset S_1 \subseteq S^{g^{-1}} \subset T_1$. By Lemma \ref{funkycompatibleordering}, there exists a $\mathcal{T}$-order $<$ on $X(T_1)$. Note that $<$ is a $\Gamma$-order by the proof of Corollary \ref{funkycompatibleorderingcor}. By Lemma \ref{preservesGammabase}, $<_g$ is a $\mathcal{T}^g$-order on $X(T)$. Since $S=T_s$ and $\mathcal{T}^g=\{1 \subset (S_1)^g \subseteq S \subset T\}$, again using the argument in Corollary \ref{funkycompatibleorderingcor}, we see that $<_g$ is a $\Gamma$-order on $X(T)$.
\end{proof}
\end{Corollary} 

\section{Combinatorial structures associated to algebraic groups}\label{combinatoricschap}

\noindent In this section we investigate some combinatorial structures associated to reductive algebraic groups.

\vspace{2mm}\noindent In $\S \ref{apsubsystems}$ we classify the almost primitive subsystems of each irreducible abstract root system. In $\S \ref{abstractindex}$ we give a purely combinatorial definition of an index that we call an abstract index. In $\S \ref{embeddingsofabstractindices}$ we introduce the notion of an embedding of abstract indices. Finally, in $\S \ref{charmaximalexceptional}$ we characterise when an embedding of abstract indices is maximal.

\subsection{Almost primitive subsystems}\label{apsubsystems}

\noindent Let $\Phi$ be an abstract root system with Weyl group $W$ and isometry group $\Iso(\Phi)$. Let $p$ either be a prime number or $0$. 

\vspace{2mm}\noindent In this section we introduce the notion of an almost primitive subsystem of $\Phi$ (a generalisation of a maximal subsystem of $\Phi$). We then classify all $W$-conjugacy classes of almost primitive subsystems $\Psi$ of irreducible $\Phi$ and list their associated groups $\Stab_W(\Delta)$ and $\Stab_{\Iso(\Phi)}(\Delta)$ (where $\Delta$ is a base of $\Psi$).

\vspace{2mm}\noindent Associated to $\Phi$ is a set of structure constants $c_{\alpha\beta}^{mn} \in \Z$ for $\alpha, \beta \in \Phi$ and integers $m,n>0$ (see for instance $\S 11.1$ of \cite{MT}). A closed subsystem $\Psi$ of $\Phi$ is \textit{$p$-closed} if, for $\alpha, \beta \in \Psi$ and integers $m, n > 0$, we have $m\alpha+n\beta \in \Psi$ whenever $c_{\alpha\beta}^{mn}$ is non-zero mod $p$.

\vspace{2mm}\noindent Let $\Psi$ be a $p$-closed proper subsystem of $\Phi$ with Weyl group $W_{\Psi}<W$. Let $\Delta$ be a base of $\Psi$. By Proposition $28$ of \cite{C}, there is a natural isomorphism from $\Stab_W(\Delta)$ to the factor group $N_W(W_{\Psi})\big/W_{\Psi}$. Similarly, there is a natural isomorphism from $\Stab_{\Iso(\Phi)}(\Delta)$ to $N_{\Iso(\Phi)}(W_{\Psi})\big/W_{\Psi}$. Hence, by the third isomorphism theorem, $\Stab_{\Iso(\Phi)}(\Delta)$ is an extension of $\Stab_W(\Delta)$ by $N_{\Iso(\Phi)}(W_{\Psi})\big/N_W(W_{\Psi})$.

\vspace{2mm}\noindent For any $\alpha \in \Phi$ and associated reflection $w_{\alpha} \in W$, we have $gw_{\alpha}g^{-1}=w_{g(\alpha)}$ for any $g \in W$ (see the proof of Proposition $2.1.8$ of \cite{C1}). It follows that $\Stab_W(\Psi)=N_W(W_{\Psi})$ and $\Stab_{\Iso(\Phi)}(\Psi)=N_{\Iso(\Phi)}(W_{\Psi})$.

\vspace{2mm}\noindent Following \cite{K} we introduce the following definitions. If $\Psi$ is maximal among $p$-closed subsystems of $\Phi$ that are invariant under $N_W(W_{\Psi})$ (resp. $N_{\Iso(\Phi)}(W_{\Psi})$) then $\Psi$ is \textit{$p$-primitive} (resp. \textit{almost $p$-primitive}) in $\Phi$. If $\Psi$ is $p$-primitive (resp. almost $p$-primitive) in $\Phi$ for some $p \geq 0$ then $\Psi$ is \textit{primitive} (resp. \textit{almost primitive}) in $\Phi$.

\vspace{2mm}\noindent If $\Phi$ is irreducible then the $W$-conjugacy classes of almost primitive subsystems $\Psi$ of $\Phi$ are listed in Table \ref{classlist} (if $\Phi$ is of classical type) and Table \ref{exceptlist} (if $\Phi$ is of exceptional type) along with their respective groups $\Stab_W(\Delta)$. If $\Psi$ is almost $p$-primitive in $\Phi$ only for certain characteristic(s) $p$ then we specify $p$. This information is well-known, and can be found directly by inspection of the root systems $\Phi$. Alternatively, refer to Table $1$ of \cite{K} for the cases where $\Phi$ is of classical type and to $\S 2$ of \cite{LSS} for the cases where $\Phi$ is of exceptional type. If $\Psi$ is almost primitive but not primitive in $\Phi$ then, by inspection of Tables \ref{classlist} and \ref{exceptlist}, there is only one possibility: $\Psi=A_2$ and $\Phi = D_4$. 

\vspace{2mm}\noindent We use the following parameters in Table \ref{classlist} and Lemma \ref{isoautgps}. Let $j,k,l,m$ and $n$ be positive integers that satisfy $2 \leq j \leq n$, $1 \leq k < n/2$ and $lm=n$. 

\begin{table}[!ht]
\centering\caption{Almost primitive subsystems $\Psi$ of classical type $\Phi$}\label{classlist}
\begin{tabular}{c | c c } 
$\Phi$ & $\Psi$ & $\Stab_W(\Delta)$  \\ [0.2ex] \hline 
$A_{n-1}$ & $A_{k-1}A_{n-k-1}$, $(A_{l-1})^m$ & $1$, $S_m$ \topstrut \\ 
$B_n$ ($p \neq 2$) & $B_{n-1}$, $D_jB_{n-j}$ & $\Z_2$, $\Z_2$ \topstrut \\ 
$B_n$ ($p=2$) & $B_kB_{n-k}$, $D_n$, $(B_l)^m$ & $1$, $\Z_2$, $S_m$ \topstrut \\ 
$C_n$ & $C_kC_{n-k}$, $\widetilde{D_n}$ ($p=2$), $(C_l)^m$ & $1$, $\Z_2$, $S_m$ \topstrut \\ 
$D_n$ ($n>4$) & $A_{n-1}$ ($\hspace{0.5mm}(n,2)$ classes), $D_kD_{n-k}$, $(D_l)^m~$ & $~\Z_2$, $\Z_2$, $(\Z_2)^{m-1} \hspace{-0.5mm}\rtimes\hspace{-0.5mm} S_m$ \topstrut \\ 
$D_4$ & $(A_1)^4$, $A_3$ ($3$ classes), $A_2$ & $(\Z_2)^2$, $\Z_2$, $\Z_2$ \topstrut \\ 
\end{tabular}
\end{table}

\begin{table}[!ht]
\centering\caption{Almost primitive subsystems $\Psi$ of exceptional type $\Phi$}\label{exceptlist}
\begin{tabular}{c | c c } 
$\Phi$ & $\Psi$ & $\Stab_W(\Delta)$  \\ [0.2ex] \hline 
$G_2$ & $A_2$, $A_1\widetilde{A_1}$, $\widetilde{A_2}$ ($p=3$) & $\Z_2$, $1$, $\Z_2$ \topstrut \\ 
$F_4$ & $B_4$, $A_2\widetilde{A_2}$, $C_3A_1$ ($p \neq 2$), $D_4$ & $1$, $\Z_2$, $1$, $S_3$ \topstrut \\ 
$\!F_4$ ($p=2$) & $C_4$, $\widetilde{D_4}$ 
 & $1$, $S_3$ 
 \topstrut \\
$E_6$ & $(A_2)^3$, $D_5$, $A_5A_1$, $D_4$, $\varnothing$ & $S_3$, $1$, $1$, $S_3$, $W(E_6)$ \topstrut \\ 
$E_7$ & $A_5A_2$, $D_6A_1$, $D_4(A_1)^3$, $A_7$, $E_6$, $(A_1)^7$, $\varnothing$ & ~$\Z_2$, $1$, $S_3$, $\Z_2$, $\Z_2$, $\GL_3(2)$, $W(E_7)$  \topstrut \\ 
\multirow{2}{*}{$E_8$} & $(A_4)^2$, $A_8$, $E_6A_2$, $D_8$, $(A_2)^4$, & $\Z_4$, $\Z_2$, $\Z_2$, $1$, $\GL_2(3)$,  \topstrut \\ 
& $E_7A_1$, $(D_4)^2$, $(A_1)^8$, $\varnothing$ & $1$, $\Z_2 \hspace{-0.5mm}\times\hspace{-0.5mm} S_3$, $\AGL_3(2)$, $W(E_8)$ \topstrut \\ 
\end{tabular}
\end{table}

\begin{Lemma}\label{isoautgps} Let $\Phi$ be irreducible and let $\Psi$ be almost primitive in $\Phi$ (as in Tables \ref{classlist} and \ref{exceptlist}). Then $\Stab_{\Iso(\Phi)}(\Delta)=\Stab_W(\Delta)$ unless one of the following occurs.

\vspace{2mm}\noindent $(i)$ If $\Phi=D_4$ and $\Psi=(A_1)^4$ (resp. $A_3$, $A_2$) then $\Stab_{\Iso(\Phi)}(\Delta)=S_4$ (resp. $(\Z_2)^2$, $\Z_2 \times S_3$).

\vspace{1mm}\noindent $(ii)$ If $\Phi=D_n$ ($n>4$ is even) and $\Psi=D_kD_{n-k}$ (resp. $(D_l)^m$) then $\Stab_{\Iso(\Phi)}(\Delta)=(\Z_2)^2$ (resp. $(\Z_2)^m \rtimes S_m$).

\vspace{2mm}\noindent $(iii)$ If $\Phi=A_n$ ($n>1$), $D_n$ ($n$ is odd) or $E_6$ then $\Stab_{\Iso(\Phi)}(\Delta)=\Z_2 \times \Stab_W(\Delta)$.

\begin{proof} By careful inspection of the relevant root systems.
\end{proof} 


\end{Lemma}

\noindent Let $\Psi^{\perp}:=\big\{\alpha \in \Phi ~\big|~ ( \alpha, \beta ) =0 \textnormal{ for all } \beta \in \Psi \big\}$ denote the perp of $\Psi$ in $\Phi$.

\begin{Lemma}\label{trivialperp} Let $\Phi$ be irreducible and let $\Psi$ be almost $p$-primitive in $\Phi$. Then $\Psi^{\perp}=\varnothing$ unless $p \neq 2$, $\Phi=B_n$ and $\Psi=B_{n-1}$ for some integer $n \geq 2$.
\begin{proof} Assume that $\Psi^{\perp}$ is non-empty. Then $\Psi \Psi^{\perp}$ is a $N_{\Iso(\Phi)}(W_{\Psi})$-invariant subsystem of $\Phi$ that contains $\Psi$. Since $\Psi$ is almost $p$-primitive, it follows that $\Psi \Psi^{\perp}$ is not $p$-closed in $\Phi$. By inspection of Tables \ref{classlist} and \ref{exceptlist}, this occurs only when $p \neq 2$, $\Phi=B_n$ and $\Psi=B_{n-1}$.
\end{proof}
\end{Lemma}

\subsection{Abstract indices}\label{abstractindex}

\noindent In this section we interpret an index to be an abstract combinatorial object. We formalise this by defining an \textit{abstract index} without reference to any algebraic group. This mirrors the theory of connected reductive algebraic groups over an algebraically closed field, where we consider a root system to be both an invariant of a group as well as an abstract combinatorial object. We present some results of Satake \cite{Sa0} and Tits \cite{T} using this combinatorial language.

\vspace{2mm}\noindent Let $E$ be a finite-dimensional real inner product space. Let $\Delta$ be a system of simple roots in $E$ (we require that each simple reflection is an isometry, but we do not require that $\Delta$ spans $E$). Let $\Delta_0$ be a subset of $\Delta$. Let $\Pi$ be a subgroup of the isometry group of $E$ that stabilises both $\Delta$ and $\Delta_0$.

\vspace{2mm}\noindent We use the following notation. Let $\Iso(E)$ denote the isometry group of $E$. For any subspace $M$ of $E$, $M^{\perp}$ is the orthogonal complement of $M$ in $E$, $\Stab_{\Iso(E)}(M)$ (resp. $\Fix_{\Iso(E)}(M)$) is the subgroup of $\Iso(E)$ that stabilises $M$ (resp. acts trivially on $M$) and 
$\pi_M:\Stab_{\Iso(E)}(M) \to \Stab_{\Iso(E)}(M)\big/\Fix_{\Iso(E)}(M)$ is the natural projection. Let $\langle \Delta \rangle$ (resp. $\langle \Delta_0 \rangle$) be the root system generated by $\Delta$ (resp. $\Delta_0$). Let $\langle \Delta \rangle^+$ be the set of positive roots of $\langle \Delta \rangle$ with respect to $\Delta$. Let $W_{\Delta}$ be the Weyl group of $\Delta$ and let $w_0$ be the longest element of $W_{\Delta}$. Let $E_{\Delta}$ (resp. $E_{\Delta_0}$) be the subspace of $E$ that is spanned by $\Delta$ (resp. $\Delta_0$) and let $\overline{E}:=(E_{\Delta})^{\perp}$. Let $E_s$ be the fixed point subspace of $(E_{\Delta_0})^{\perp}$ under $\Pi$ and let $E_a:=(E_s)^{\perp}$. Denote $\overline{E}_s:=\overline{E} \cap E_s$, $\overline{E}_a:=\overline{E} \cap E_a$ and $\overline{E_a}:=\smash{(E_{\Delta_0})^{\perp}} \cap E_a$. Let $_s\Delta$ be the image of the projection of $\Delta$ onto $E_s$. The orbits of the action of $\Pi$ on $\Delta \setminus \Delta_0$ are called \textit{distinguished}.

\begin{Lemma}\label{indexlemma} $\Delta \cap E_a = \Delta_0$
\begin{proof} By definition, $\Delta \cap E_a \supseteq \Delta_0$. Assume that $\Delta \cap E_a \supset \Delta_0$. That is, $\Delta \cap \overline{E_a}$ is non-trivial. Let $v$ be the sum of all roots in $\Delta \cap \overline{E_a}$ and let $E_v$ be the subspace of $E$ spanned by $v$. Observe that $v$ is non-zero (since $\Delta$ is linearly independent) and is fixed by $\Pi$ (since $\Delta \cap \overline{E_a}$ is $\Pi$-stable). So $E_v$ is a non-trivial subspace of $\overline{E_a}$ that is fixed pointwise by $\Pi$. But this contradicts the definition of $\overline{E_a}$. 
\end{proof}
\end{Lemma}

\begin{Definition}[\textbf{Satake \cite{Sa0}, Tits \cite{T}}]\label{abstractindexdefn} The quadruple $\mathcal{I}:=(E,\Delta,\Delta_0,\Pi )$ is an \textit{abstract index} if the following three properties are satisfied: 

\vspace{2mm}\noindent $(i)$ \textit{(relative root system)} ${}_s\Delta$ is a system of simple roots for a non-reduced abstract root system, 

\vspace{1mm}\noindent $(ii)$ \textit{(self-opposition)} the involution $-w_0$ of $\Delta$ commutes with $\Pi$ and stabilises $\Delta_0$, and

\vspace{1mm}\noindent $(iii)$ \textit{(inductive closure)} if one removes a distinguished orbit $\mathcal{O}$ from $\Delta \setminus \Delta_0$ then the result $(E,\Delta \setminus \mathcal{O},\Delta_0,\Pi )$ is again an abstract index.
\end{Definition}

\noindent By $\S 2.5.1$ of \cite{T}, two elements of $\Delta \setminus \Delta_0$ have the same projection to $E_s$ if and only if they belong to the same distinguished orbit of $\Pi$. So we have a $1-1$ correspondence between $_s\Delta$ and the set of distinguished orbits of $\Pi$. We can compute the isomorphism type of $_s\Delta$ using the method of $\S 2.5$ of \cite{T}.

\vspace{2mm}\noindent Let $\mathcal{I}=(E,\Delta,\Delta_0,\Pi)$ be an abstract index. If $\Delta$ spans $E$ then $\mathcal{I}$ is \textit{semisimple}. If $\Pi$ is trivial then $\mathcal{I}$ is of \textit{inner type}, otherwise $\mathcal{I}$ is of \textit{outer type}. If $\Delta_0=\varnothing$ and $\Pi$ is trivial then $\mathcal{I}$ is \textit{split}. If $\Delta_0 = \varnothing$ then $\mathcal{I}$ is \textit{quasisplit}. 
If $\Delta=\Delta_0$ and $\overline{E}_s$ is trivial then $\mathcal{I}$ is \textit{anisotropic}, otherwise $\mathcal{I}$ is \textit{isotropic}. If $\mathcal{I}$ is semisimple and $\Delta$ is irreducible (resp. classical, exceptional) then $\mathcal{I}$ is \textit{irreducible} (resp. \textit{classical}, \textit{exceptional}).

\vspace{2mm}\noindent The \textit{type} of $\mathcal{I}$ is the isomorphism class of $\Delta$ (as a root system) along with the isomorphism class of $\Pi$. 
If $\mathcal{I}$ is irreducible then the isomorphism class of $\Pi$ is uniquely determined by its order $g:=|\Pi|$. Then, as in Table II of \cite{T}, we use the shorthand that $\mathcal{I}$ is \textit{of type ${}^{g}\hspace{-0.5mm}\Delta$}.

\vspace{2mm}\noindent The abstract index $\mathcal{I}':=\big(E_{\Delta},\Delta,\Delta_0,\pi_{E_{\Delta}}(\Pi)\big)$ is the \textit{semisimple component} of $\mathcal{I}$. The abstract index $\mathcal{I}_a:=(E_a,\Delta_0, \Delta_0, \Pi)$ is the \textit{anisotropic kernel} of $\mathcal{I}$ (note that $\Pi$ stabilises $E_a$ and $\Pi \cap \Fix_{\Iso(E)}(E_a)$ is trivial). The \textit{rank} of $\mathcal{I}$ is $r(\mathcal{I}):=|\Delta|+\dim(\smash{\overline{E}})$, the \textit{relative rank} of $\mathcal{I}$ is $r_s(\mathcal{I}):=\dim(E_s)$ and the \textit{semisimple rank} of $\mathcal{I}$ is $r(\mathcal{I}')=|\Delta|$. Let $I$ be a $\Pi$-stable subset of $\Delta$ that contains $\Delta_0$. The abstract index $(E,I, \Delta_0, \Pi)$ is a \textit{subindex} of $\mathcal{I}$ (i.e. a subindex of $\mathcal{I}$ is obtained by removing distinguished orbits). The \textit{type} of $(E,I,\Delta_0,\Pi)$ is the $W_{\Delta}$-conjugacy class of $I$ in $\langle \Delta \rangle$. In particular, $\mathcal{I}_m:=(E,\Delta_0, \Delta_0, \Pi)$ is the \textit{minimal subindex} of $\mathcal{I}$.

\vspace{2mm}\noindent A $\Pi$-stable system of simple roots (resp. root system) in $E$ is called \textit{$\Pi$-irreducible} if it cannot be written as a non-trivial disjoint union of $\Pi$-stable orthogonal systems of simple roots (resp. root systems). Up to rearrangement, there is a unique decomposition of $\Delta$ as a disjoint union $\Delta=\coprod_{i=1}^l \Delta_i$ of $\Pi$-irreducible subsystems of simple roots $\Delta_i$ for $1 \leq i \leq l$. If $\mathcal{I}$ is semisimple and this decomposition is trivial then $\mathcal{I}$ is \textit{simple}. For each $1 \leq i \leq l$, let $E_i$ be the subspace of $E$ spanned by $\Delta_i$, let $(\Delta_i)_0:=\Delta_0 \cap E_i$ and $\Pi_i :=\pi_{E_i}(\Pi)$. A \textit{simple component} of $\mathcal{I}$ is an abstract index $\mathcal{I}_i:=\big(E_i, \Delta_i, (\Delta_i)_0, \Pi_i\big)$ for some $1 \leq i \leq l$. The \textit{central component} of $\mathcal{I}$ is the abstract index $\overline{\mathcal{I}}:=\big(\overline{E}, \varnothing, \varnothing,\pi_{\overline{E}}(\Pi)\big)$. We use the notation $\mathcal{I}=\mathcal{I}' \times \overline{\mathcal{I}} = \mathcal{I}_1 \times ... \times \mathcal{I}_l \times \overline{\mathcal{I}}$ to represent the decomposition of $\mathcal{I}$ into its components.


\vspace{2mm}\noindent Let $\smash{'}\mathcal{I}=(\smash{'}\hspace{-0.3mm}E,\smash{'}\hspace{-0.4mm}\Delta,\smash{'}\hspace{-0.4mm}\Delta_0,\smash{'}\Pi)$ be another abstract index. An \textit{isomorphism} from $\mathcal{I}$ to $\smash{'}\mathcal{I}$ is a bijective isometry $\psi:E \to \smash{'}\hspace{-0.3mm}E$ that satisfies $\psi(\Delta)=\smash{'}\hspace{-0.4mm}\Delta$, $\psi(\Delta_0)=\smash{'}\hspace{-0.4mm}\Delta_0$ and $\psi\Pi \psi^{-1}=\smash{'}\Pi$. We denote $\psi(\mathcal{I}):=\smash{'}\mathcal{I}$.

\begin{Theorem}[\textbf{Satake, Tits}]\label{Titsclass} Isomorphism classes of irreducible abstract indices are classified in Table II of \cite{T}. 
\end{Theorem}

\noindent In $\S 3.1.2$ of \cite{T}, Tits explains how the problem of classifying (isomorphism classes of) abstract indices reduces to that of classifying irreducible abstract indices.

\vspace{2mm}\noindent Now let $k$ be a field and let $G$ be a connected reductive algebraic $k$-group. Let $T$ be a maximal $k$-torus of $G$ that contains a maximal $k$-split torus of $G$, let $<$ be a $\Gamma$-order on $X(T)$ and let $\mathcal{I}(G)$ be the index of $G$ with respect to $T$ and $<$. 

\begin{Theorem}[\textbf{Satake, Tits}, \textit{Prop. $2$ of \cite{Sa0}, $\S 3.2$ of \cite{T}}]\label{indeximpo} The quadruple $\mathcal{I}(G)$ is an abstract index. Let $\mathcal{I}_1(G)$ be another index of $G$ (say with respect to $T_1$ and $<_1$). Then $\mathcal{I}(G)$ and $\mathcal{I}_1(G)$ are isomorphic as abstract indices.
\end{Theorem}

\noindent Up to rearrangement, we can uniquely decompose $G$ as a commuting product $\smash{\big(\prod_{i=1}^l G_i\big) {Z(G)}^{\circ}}$ where each $G_i$ is $k$-simple and ${Z(G)}^{\circ}$ is a $k$-torus. The $G_i$'s are \textit{$k$-simple components} of $G$.

\begin{Proposition}[\textit{Definition chasing, also see $\S 2$ of \cite{T}}]\label{abstractTits} If $G$ is semisimple (resp. absolutely simple, $k$-simple, classical, exceptional, $k$-(quasi)split, $k$-(an)isotropic, of inner/outer type) then $\mathcal{I}(G)$ is semisimple (resp. irreducible, simple, classical, exceptional, (quasi)split, (an)isotropic, of inner/outer type) as an abstract index. The number of simple components of $G$ and $\mathcal{I}(G)$ are equal and, up to rearrangement, $\mathcal{I}(G_i)=\mathcal{I}(G)_i$ for each $1 \leq i \leq l$. In addition, $X\big({Z(G)}^{\circ}\big)_{\R}=\overline{X(T)_{\R}}$, $X(S)_{\R} =(X(T)_{\R})_s$, $X_0(S)_{\R} =((X(T)_{\R})_a$, $\mathcal{I}(G_a)=\mathcal{I}(G)_a$, $\mathcal{I}(G')=\mathcal{I}(G)'$, $\mathcal{I}\big({Z(G)}^{\circ}\big)=\overline{\mathcal{I}(G)}$, $r(G)=r\big(\mathcal{I}(G)\big)$ and $r_k(G)=r_s\big(\mathcal{I}(G)\big)$.
\end{Proposition}

\noindent For a given field $k$, an abstract index $\mathcal{I}$ is \textit{$k$-admissible} if there exists a connected reductive algebraic $k$-group $X$ such that the index of $X$ is isomorphic to $\mathcal{I}$.

\vspace{2mm}\noindent For example, the abstract index $\dynkin{F}{II}$ is $\R$-admissible but not $\F_q$-admissible for any finite field $\F_q$. In $\S 3.3$ of \cite{T}, Tits gives some necessary conditions for the $k$-admissability of an abstract index when $k$ is finite, $k=\R$, $k$ is $\mathfrak{p}$-adic and $k$ is a number field.

\begin{Theorem}[\textbf{Tits}, \textit{$\S 3$ of \cite{T}}]\label{Titsadmissability} Any abstract index is $k$-admissible for some field $k$. 
\end{Theorem}




\subsection{Embeddings of abstract indices}\label{embeddingsofabstractindices}

\noindent In this section we define a combinatorial structure called an \textit{embedding of abstract indices}. This definition is made without reference to any algebraic group. 

\vspace{2mm}\noindent Let $p$ either be a prime number or $0$. Let $\mathcal{H}=(E,\Delta, \Delta_0, \Pi)$ be an abstract index. As in $\S \ref{abstractindex}$ we use the following notation: $E_s$ is the fixed point subspace of $(E_{\Delta_0})^{\perp}$ under $\Pi$, $E_a:=(E_s)^{\perp}$ and $\overline{E}:=(E_{\Delta})^{\perp}$. Let $\Psi:=\langle \Delta \rangle$ and $\Psi_0:=\langle \Delta_0 \rangle$.

\vspace{2mm}\noindent Let $\Phi$ be a $\Pi$-stable root system in $E$ such that $\Psi$ is a $p$-closed subsystem of $\Phi$. So $\Pi \leq \Stab_{\Iso(\Phi)}(\Delta)$. We introduce some additional notation. Let $W:=W(\Phi)$, $\Phi_s:=\Phi \cap E_s$ and $\Phi_a :=\Phi \cap E_a$. Note that $\Phi_a \cap \Psi = \Psi_0$ by Lemma \ref{indexlemma}. The \textit{independent subsystem} ${}_{in}\Phi_a$ of $\Phi_a$ is the union of irreducible components of $\Phi_a$ that is maximal with respect to the condition ${}_{in}\Phi_a \subseteq \Psi_0$. If ${}_{in}\Phi_a=\Phi_a$ then the pair $(\mathcal{H},\Phi)$ is \textit{independent} (a.k.a. $\mathcal{H}$ \textit{is independent in} $\Phi$).

\vspace{2mm}\noindent Let $V:=E_{\Phi} \cap \overline{E}$ and let $V^{\Pi}$ be the fixed point subspace of $V$ under $\Pi$. The pair $(\mathcal{H},\Phi)$ is \textit{maximal} (a.k.a. $\mathcal{H}$ \textit{is maximal in} $\Phi$) if $V^{\Pi}$ is trivial and $\Psi$ is maximal among $p$-closed $\Pi$-stable (proper) subsystems of $\Phi$. 
In $\S \ref{charmaximalexceptional}$ we characterise when $(\mathcal{H},\Phi)$ is maximal for the cases where $\Phi$ is of exceptional type.

\begin{Lemma}\label{complemma} Let $\Pi \leq W$. Then $W_{\Delta_0} \rtimes \Pi \leq W_{\Phi_a}$ and $r(\Phi_a)=\dim(E)-r_s(\mathcal{H})$.
\begin{proof} By definition, $\Phi_a$ contains $\Delta_0$ and $\Pi$ acts trivially on $E_s$. Then $W_{\Delta_0} \rtimes \Pi \leq W_{\Phi_a}$ since $\Pi \leq W$. Observe that $E_a$ is minimal among all subspaces of $E$ that contain $\Delta_0$ and that have an orthogonal complement in $E$ that is fixed pointwise by $\Pi$. But $E_{\Phi_a}$ contains $\Delta_0$ and $(E_{\Phi_a})^{\perp}$ is fixed pointwise by $\Pi \leq W_{\Phi_a}$. So $E_{\Phi_a}=E_a$ by minimality. Hence $r(\Phi_a) = \dim(E_a)=\dim(E)-r_s(\mathcal{H})$.
\end{proof}
\end{Lemma}

\noindent Let $'\mathcal{H}=(E,\smash{'}\hspace{-0.4mm}\Delta,\smash{'}\hspace{-0.4mm}\Delta_0,\smash{'}\Pi)$ be another abstract index that shares the ambient space $E$ with $\mathcal{H}$ such that $\smash{'}\Pi$ stabilises $\Phi$ and $\langle \smash{'}\hspace{-0.4mm}\Delta \rangle$ is a $p$-closed subsystem of $\Phi$. A \textit{conjugation} from $\mathcal{H}$ to $'\mathcal{H}$ in $\Phi$ is an element $w \in W$ that satisfies $w(\Delta)=\smash{'}\hspace{-0.4mm}\Delta$, $w(\Delta_0)=\smash{'}\hspace{-0.4mm}\Delta_0$ and $w\Pi w^{-1}=\smash{'}\Pi$. If such a conjugation $w \in W$ exists then the abstract indices $\mathcal{H}$ and $'\mathcal{H}$ are \textit{conjugate in $\Phi$}, and we write $w(\mathcal{H})=\smash{'}\mathcal{H}$.

\vspace{2mm}\noindent We use the notation $\smash{'}\Psi:=\langle \smash{'}\hspace{-0.4mm}\Delta \rangle$, $\smash{'}\Psi_0:=\langle \smash{'}\hspace{-0.4mm}\Delta_0 \rangle$, $\smash{'}\hspace{-0.3mm}E_s$ is the fixed point subspace of $(E\hspace{0.1mm}_{\smash{'}\hspace{-0.6mm}\Delta_0})^{\perp}$ under $\smash{'}\Pi$, $\smash{'}\hspace{-0.3mm}E_a:=(\smash{'}\hspace{-0.3mm}E_s)^{\perp}$, $\smash{'}\Phi_a:=\Phi \cap \smash{'}\hspace{-0.3mm}E_a$, $\hspace{2mm}'\!\!\!\!{}_{in}\Phi_a$ is the independent subsystem of $\smash{'}\Phi_a$, $\smash{'}V$ is the largest subspace of $E_{\Phi}$ that is perpendicular to all roots in $\smash{'}\hspace{-0.4mm}\Delta$ and $\smash{'}V^{\smash{'}\Pi}$ is the fixed point subspace of $\smash{'}V$ under $\smash{'}\Pi$.

\begin{Lemma}\label{trickybitty} Let $\mathcal{H}$ be irreducible and let $\Delta_0 \cong \smash{'}\hspace{-0.4mm}\Delta_0$ (as root systems). Assume that there exists $w \in W$ such that $w(\Delta)=\smash{'}\hspace{-0.4mm}\Delta$ and $w\Pi w^{-1}=\smash{'}\Pi$, but that $\Delta_0$ is not $W$-conjugate to $\smash{'}\hspace{-0.4mm}\Delta_0$. Then $\Pi$ is trivial, $\Delta \cong D_n$ and $\Delta_0  \cong (A_{d-1})^r$ for some integers $n,d,r$ such that $n=rd$ and $d\neq 1$ is a power of $2$. If $n \neq 4$ then $\Stab_W(\Delta)=\Fix_W(\Delta)$.
\begin{proof} Observe that $\mathcal{H}$ is isomorphic to $w^{-1}('\mathcal{H})$ since $\Delta_0 \cong \smash{'}\hspace{-0.4mm}\Delta_0$. That is, there exists $\psi \in \Iso(E)$ such that $\psi(\Delta)=\Delta$, $\psi\Pi\psi^{-1}=\Pi$ and $\psi(\Delta_0)=w^{-1}(\smash{'}\hspace{-0.4mm}\Delta_0$). By assumption, $\psi(\Delta_0) \neq \Delta_0$.

\vspace{2mm}\noindent Using Table II of \cite{T}, we observe that $\Delta_0$ is stabilised by $\Stab_{\Iso(\Phi)}(\Delta)$ unless $\Pi$ is trivial, $\Delta \cong D_n$ and $\Delta_0  \cong (A_{d-1})^r$ for some integers $n,d,r$ such that $n=rd$ and $d\neq 1$ is a power of $2$. Then $\psi$ acts on $\Delta$ as a non-trivial graph automorphism. Using $\S 7$ of \cite{C}, we check that $\psi(\Delta_0)$ is not $W_{\Delta}$-conjugate to $\Delta_0$ (since $n$ is even).

\vspace{2mm}\noindent Assume that $n \neq 4$. Then there is a unique non-trivial graph automorphism of $\Delta$. If $\Stab_W(\Delta)\big/\Fix_W(\Delta)$ is non-trivial then $\psi \in W$, which is a contradiction.
\end{proof}
\end{Lemma}

\begin{Lemma}\label{conjlemma} Let $w \in W$ such that $w(\Delta)=\smash{'}\hspace{-0.4mm}\Delta$ and $w\Pi w^{-1}=\smash{'}\Pi$. Then $(\mathcal{H},\Phi)$ is maximal if and only if $('\mathcal{H},\Phi)$ is maximal. If in addition $w(\Delta_0)=\smash{'}\hspace{-0.4mm}\Delta_0$ then $w(E_a)=\smash{'}\hspace{-0.3mm}E_a$, $w(\Phi_a)=\smash{'}\Phi_a$ and $w({}_{in}\Phi_a)=\hspace{2mm}'\!\!\!\!{}_{in}\Phi_a$.
\begin{proof} Let $w \in W$ such that $w(\Delta)=\smash{'}\hspace{-0.4mm}\Delta$ and $w\Pi w^{-1}=\smash{'}\Pi$. Let $\Omega$ be a $p$-closed $\Pi$-stable proper subsystem of $\Phi$ that strictly contains $\Psi$. Then $w(\Omega)$ is a $p$-closed $\smash{'}\Pi$-stable proper subsystem of $\Phi$ that strictly contains $\smash{'}\Psi$. Observe that $w(V^{\Pi})=\smash{'}V^{\smash{'}\Pi}$. The first assertion follows (from this and the reverse argument).

\vspace{2mm}\noindent In addition let $w(\Delta_0)=\smash{'}\hspace{-0.4mm}\Delta_0$. Let $\tilde{W}_{\Delta_0}:=W_{\Delta_0} \rtimes \Pi$, where $\Pi$ acts on $W_{\Delta_0}$ by sending $w_{\alpha} \mapsto w_{\tau(\alpha)}$ for $\alpha \in \Delta_0$ and $\tau \in \Pi$. The fixed point subspace of $E$ under $\tilde{W}_{\Delta_0}$ is $E_s$. Similarly, the fixed point subspace of $E$ under $\tilde{W}_{\smash{'}\hspace{-0.4mm}\Delta_0}:=W_{\smash{'}\hspace{-0.4mm}\Delta_0} \rtimes \smash{'}\Pi$ is $\smash{'}\hspace{-0.3mm}E_s$. Observe that $w\tilde{W}_{\Delta_0}w^{-1}= \tilde{W}_{\smash{'}\hspace{-0.4mm}\Delta_0}$. Hence $w(E_s)=\smash{'}\hspace{-0.3mm}E_s$ and so $w(E_a)=\smash{'}\hspace{-0.3mm}E_a$ and $w(\Phi_a)=\smash{'}\Phi_a$. Then $w({}_{in}\Phi_a)=\hspace{2mm}'\!\!\!\!{}_{in}\Phi_a$ since $w(\Psi_0)=\smash{'}\Psi_0$.
\end{proof}
\end{Lemma}

\noindent We now restate the following definition from the introduction. 


\begin{mydefinition!} Let $p$ either be a prime number or $0$. A \textit{$p$-embedding of abstract indices} is a triple that consists of two abstract indices $\mathcal{G}=(F,\Lambda,\Lambda_0, \Xi)$ and $\mathcal{H}=(E,\Delta,\Delta_0,\Pi)$, and a bijective isometry $\theta:E \to F$ that satisfies the following conditions:

\vspace{1.5mm}\noindent $(A.1)$  $\Phi:=\langle \Lambda \rangle$ is $\Pi^{\theta}$-stable, $\theta(\Delta) \subset \Phi^+$, $\Psi:=\langle \theta(\Delta) \rangle$ is a $p$-closed subsystem of $\Phi$, $\Lambda_a:= \Lambda \cap \theta(E_a)$ is a base of $\Phi_a:=\Phi \cap \theta(E_a)$ and ${}_{in}\Lambda_a \subseteq \theta(\Delta_0)$ (where ${}_{in}\Lambda_a$ is the union of irreducible components of $\Lambda_a$ that is maximal with respect to the condition ${}_{in}\Lambda_a \subseteq \langle \theta(\Delta_0) \rangle =:\Psi_0$).


\vspace{1.5mm}\noindent $(A.2)$ For every $\sigma \in \Pi$ there exists a unique $w_{\sigma} \in W_{\Lambda_a}$ such that $w_{\sigma} \sigma^{\theta} \in \Xi$. Moreover, the map $\Pi \to \Xi$ given by $\sigma \mapsto w_{\sigma} \sigma^{\theta}$ is surjective.

\vspace{1.5mm}\noindent $(A.3)$ $\Lambda_a$ is $\Xi$-stable and contains $\Lambda_0$.

\vspace{1.5mm}\noindent $(A.4)$ ${}_{in}\Lambda_a$ is contained in $\Lambda_0$.
\end{mydefinition!}

\noindent An \textit{embedding of abstract indices} is a $p$-embedding of abstract indices for some $p \geq 0$. 

\vspace{2mm}\noindent So let $p$ either be a prime number or $0$ and let $(\mathcal{G},\mathcal{H},\theta)$ be a $p$-embedding of abstract indices. We introduce some terminology associated to Definition \ref{mydefinition!}.

\vspace{2mm}\noindent The subindex $\mathcal{L}:=(F,\Lambda_a,\Lambda_0, \Xi)$ of $\mathcal{G}$ is called the \textit{$(\mathcal{H},\theta)$-embedded subindex of $\mathcal{G}$}. If $\mathcal{H}$ is isotropic then $\mathcal{G}$ is also isotropic by $(A.3)$. We say that $(\mathcal{G},\mathcal{H},\theta)$ is \textit{(an)isotropic} if $\mathcal{H}$ is (an)isotropic. If $\big(\theta(\mathcal{H}),\Phi\big)$ is maximal (resp. independent) then $(\mathcal{G},\mathcal{H},\theta)$ is \textit{maximal} (resp. \textit{independent}). 
If $\mathcal{H}$ and $\mathcal{G}$ are both split (resp. quasisplit) then $(\mathcal{G},\mathcal{H},\theta)$ is \textit{split} (resp. \textit{quasisplit}). The map $\theta$ is an \textit{embedding} of $\mathcal{H}$ in $\mathcal{G}$. 

\begin{Remark}\label{abusenotation} We often abuse notation and identify $\mathcal{H}$ with its image under $\theta$. One can think of this as analogous to embedding a subgroup within another group.
\end{Remark}

\noindent Given a field $k$ with $p=\Char(k)$, we say that $(\mathcal{G},\mathcal{H},\theta)$ is \textit{$k$-admissible} if there exists a pair of connected reductive $k$-groups $H \subset G$ such that the embedding of indices of $H \subset G$ is isomorphic to $(\mathcal{G},\mathcal{H},\theta)$.

\vspace{2mm}\noindent Let $(\smash{'}\mathcal{G},\smash{'}\mathcal{H},\smash{'}\hspace{-0.1mm}\theta)$ be another $p$-embedding of abstract indices. An \textit{isomorphism} from $(\mathcal{G},\mathcal{H},\theta)$ to $(\smash{'}\mathcal{G},\smash{'}\mathcal{H},\smash{'}\hspace{-0.1mm}\theta)$ is a bijective isometry $\phi:F \to \smash{'}\hspace{-0.5mm}F$ such that $\phi(\mathcal{G})=\smash{'}\mathcal{G}$ and $\phi \theta(\mathcal{H})=\smash{'}\hspace{-0.1mm}\theta (\smash{'}\mathcal{H})$. If such an isomorphism exists then $(\mathcal{G},\mathcal{H},\theta)$ and $(\smash{'}\mathcal{G},\smash{'}\mathcal{H},\smash{'}\hspace{-0.1mm}\theta)$ are \textit{isomorphic}.

\vspace{2mm}\noindent Now assume that $F=\smash{'}\hspace{-0.5mm}F$. A \textit{conjugation} from $(\mathcal{H},\theta)$ to $(\smash{'}\mathcal{H},\smash{'}\hspace{-0.1mm}\theta)$ in $\mathcal{G}$ is an element $w \in W_{\Lambda}$ that satisfies $w\theta(\mathcal{H})= \smash{'}\hspace{-0.1mm}\theta(\smash{'}\mathcal{H})$. If such a conjugation exists then $(\mathcal{H},\theta)$ and $(\smash{'}\mathcal{H},\smash{'}\hspace{-0.1mm}\theta)$ are \textit{conjugate} in $\mathcal{G}$. If $\theta(\mathcal{H})= \smash{'}\hspace{-0.1mm}\theta(\smash{'}\mathcal{H})$ then $(\mathcal{H},\theta)$ and $(\smash{'}\mathcal{H},\smash{'}\hspace{-0.1mm}\theta)$ are \textit{equal} in $\mathcal{G}$.

\begin{Remark}\label{isoaniso} Given an abstract index $\mathcal{G}$ and a value for $p$, it is an interesting and useful exercise to classify all conjugacy classes of maximal $p$-embeddings of abstract indices $(\mathcal{H},\theta)$ in $\mathcal{G}$. We complete this exercise in Theorem \ref{classificationthm!} for each possible $p$ and for each irreducible index $\mathcal{G}$ of exceptional type (except for the cases where $\mathcal{H}$ is anisotropic, which are easy since $\Lambda=\Lambda_a$ and so axiom $(A.3)$ of Definition \ref{mydefinition!} becomes vacuous).
\end{Remark}

\subsection{Maximal embeddings of abstract indices}\label{charmaximalexceptional}

\noindent Let $p$ either be a prime number or $0$. Let $(\mathcal{G}, \mathcal{H}, \theta)$ be a $p$-embedding of abstract indices, say $\mathcal{G}=(F,\Lambda,\Lambda_0, \Xi)$ and $\mathcal{H}=(E,\Delta,\Delta_0,\Pi)$. We follow Remark \ref{abusenotation} and identify $\mathcal{H}$ with its image under $\theta$. This is equivalent to assuming that $E=F$.

\vspace{2mm}\noindent In this section we specialise to the cases where $\Lambda$ is of exceptional type and $r(\Lambda)=\dim(E)$. Our objective is to prove Proposition \ref{pmaximalexceptchar}, which characterises when $(\mathcal{G}, \mathcal{H}, \theta)$ is maximal.

\vspace{2mm}\noindent As in $\S \ref{abstractindex}$ and $\S \ref{embeddingsofabstractindices}$ we use the following notation: $\Phi:=\langle \Lambda \rangle$, $W:=W_{\Lambda}$, $\Psi:=\langle \Delta \rangle$, $\Psi_0:=\langle \Delta_0 \rangle$, $E_{\Lambda}$ is the subspace of $E$ that is spanned by $\Lambda$, $E_s$ is the fixed point subspace of $(E_{\Delta_0})^{\perp}$ under $\Pi$, $E_a:=(E_s)^{\perp}$, $\overline{E}:=(E_{\Delta})^{\perp}$, $V:=E_{\Lambda} \cap \overline{E}$ and $V^{\Pi}$ is the fixed point subspace of $V$ under $\Pi$. 

\begin{Lemma}\label{alssll} Let $\Phi$ be of exceptional type and let $r(\Phi)=\dim(E)$. Let $\Psi$ be maximal among $p$-closed subsystems of $\Phi$. Then $(\mathcal{G}, \mathcal{H}, \theta)$ is not maximal if and only if $\Pi$ is trivial and $(\Phi,\Psi)$ is either $(E_6,D_5)$ or $(E_7,E_6)$. In particular, if $\Psi$ has maximal rank in $\Phi$ then $(\mathcal{G}, \mathcal{H}, \theta)$ is maximal.
\begin{proof} Assume that $(\mathcal{G}, \mathcal{H}, \theta)$ is not maximal. Then $V^{\Pi}$ is non-trivial and so $r(\Psi)<r(\Phi)$ by definition of $V$. By inspection of Table \ref{exceptlist}, if $r(\Psi)<r(\Phi)$ then $(\Phi,\Psi)$ is one of $(E_6,D_4)$, $(E_6,D_5)$ or $(E_7,E_6)$. We can exclude the case $(\Phi,\Psi)=(E_6,D_4)$ since then $\Psi$ is contained in a $p$-closed copy of $D_5$ in $\Phi$ for any choice of $p$.

\vspace{2mm}\noindent Let $(\Phi,\Psi)=(E_6,D_5)$ or $(E_7,E_6)$. Then $V=V^{\Pi}$ has dimension $1$. Recall from Table \ref{exceptlist} that $\Stab_{\Iso(\Phi)}(\Delta) \cong \Z_2$. The non-trivial element of $\Stab_{\Iso(\Phi)}(\Delta)$ acts non-trivially on $V$. Recall from Definition \ref{mydefinition!} that $\Pi \leq \Stab_{\Iso(\Phi)}(\Delta)$. So $\Pi$ must be trivial.

\vspace{2mm}\noindent Conversely, assume that $\Pi$ is trivial and $(\Phi,\Psi)$ is either $(E_6,D_5)$ or $(E_7,E_6)$. Then $V^{\Pi}=V$, which has dimension $1$. Hence $(\mathcal{G}, \mathcal{H}, \theta)$ is not maximal.

\vspace{2mm}\noindent The second assertion of the lemma follows immediately.
\end{proof}
\end{Lemma}

\noindent We use the following notation in Table \ref{primitivemaxactions}. Let $\mathcal{D}i_{n}$ (resp. $\mathcal{SD}i_{n}$) denote the dihedral (resp. semidihedral) group of order $n$. Let $Q_8$ denote the quaternion group. Let $Y:=\GL_3(2)$ and let $P_1(Y)$ and $P_2(Y)$ be non-conjugate maximal parabolic subgroups of $Y$ (see Cases $(i)$ and $(ii)$ of the proof of Lemma \ref{maxactions} for more details).

\begin{Lemma}\label{maxactions} Let $\Phi$ be of exceptional type and let $r(\Phi)=\dim(E)$. Let $\Psi$ be a non-empty almost $p$-primitive subsystem of $\Phi$ that is not maximal among $p$-closed subsystems of $\Phi$. Then $(\Phi,\Psi)$ is one of the pairs listed in the first two columns of Table \ref{primitivemaxactions}. Given such a pair $(\Phi,\Psi)$, $(\mathcal{G}, \mathcal{H}, \theta)$ is maximal if and only if $\Pi$ is one of the groups listed in the last column of Table \ref{primitivemaxactions}.

\begin{table}[!htb]\begin{center}\caption{Subsystems of exceptional $\Phi$ that are almost primitive, but not maximal}\label{primitivemaxactions}\begin{tabular}{| c | c | c | c | c | c |}    \hline
    $\Phi$ & $\Psi$ & $\!\Stab_{\Iso(\Phi)}(\Delta)\!$ & $\Sigma$ & $\!\big(\!\Stab_{\Iso(\Phi)}(\Delta)\big)_{\Sigma}\!$ & possibilities for $\Pi$ \topstrut \\ \hline \hline

\multirow{2.3}{*}{$F_4$} & $D_4$ & $S_3$ & $B_4$ & $\Z_2$ & $\Z_3$, $S_3$ \topstrut \\ \cline{2-6}
& $\!\widetilde{D_4}$ \small{$\!(p \! = \! 2)\!$} & $S_3$  & $C_4$ & $\Z_2$ & $\Z_3$, $S_3$ \topstrut \\ \hline 

$E_6$ & $D_4$ & $\Z_2 \times S_3$ & $D_5$ & $(\Z_2)^2$ & $\Z_3$, $S_3$, $\Z_6$, $\Z_2 \hspace{-0.5mm}\times\hspace{-0.5mm} S_3$ \topstrut \\ \hline 

\multirow{3.4}{*}{$E_7$} & $D_4(A_1)^3$ & $S_3$ & $D_6A_1$ & $\Z_2$  & $\Z_3$, $S_3$ \topstrut \\ \cline{2-6}
& \multirow{2.3}{*}{$(A_1)^7$} & \multirow{2.3}{*}{$\GL_3(2)=Y$} & $D_6A_1$ & $P_1(Y) \cong S_4$ & \multirow{2.3}{*}{$\Z_7$, $\Z_7 \hspace{-0.5mm}\rtimes\hspace{-0.5mm} \Z_3$, $\GL_3(2)$} \topstrut \\ \cdashline{4-5}
&&& $D_4(A_1)^3$ & $P_2(Y) \cong S_4$ & \topstrut \\ \hline 

\multirow{6.1}{*}{$E_8$} & $(A_2)^4$ & $\GL_2(3)$ & $E_6A_2$ & $\mathcal{D}i_{12}$ & \makecell{$\Z_4$, $\mathcal{D}i_8$, $Q_8$, $\mathcal{SD}i_{16}$,  \\ $\Z_8$, $\SL_2(3)$, $\GL_2(3)$}  \topstrut \\ \cline{2-6}
& $(D_4)^2$ & $\Z_2 \times S_3$ & $D_8$ & $(\Z_2)^2$ & $\Z_3$, $S_3$, $\Z_6$, $\Z_2 \hspace{-0.5mm}\times\hspace{-0.5mm} S_3$ \topstrut \\ \cline{2-6}
& \multirow{3.45}{*}{$(A_1)^8$} & \multirow{3.45}{*}{\makecell{$\AGL_3(2)$ \\ $=L \rtimes Y$ \\ $(L \cong (\Z_2)^3)$}} & $E_7A_1$ & $Y \cong \GL_3(2)$ & \multirow{3.45}{*}{\makecell{$Y_1 \cong \GL_3(2)$, \\ $L \hspace{-0.5mm}\rtimes\hspace{-0.5mm} (\Z_7 \hspace{-0.5mm}\rtimes\hspace{-0.5mm} \Z_3)$, \\ $L \hspace{-0.5mm}\rtimes\hspace{-0.5mm} \Z_7$, $\AGL_3(2)$}} \topstrut \\ \cdashline{4-5}
&&& $D_8$ & $L \rtimes P_1(Y)$ & \topstrut \\ \cdashline{4-5}
&&& $(D_4)^2$ & $L \rtimes P_2(Y)$ & \topstrut \\ \hline 
  \end{tabular}\end{center}\end{table}

\begin{proof} We begin by explaining and justifying the entries in Table \ref{primitivemaxactions}. For any given pair $(\Phi,\Psi)$ as in the statement of the lemma, the purpose of constructing Table \ref{primitivemaxactions} is to find all conjugacy classes of subgroups $\Pi$ of $\Stab_{\Iso(\Phi)}(\Delta)$ for which there does not exist a $p$-closed $\Pi$-stable almost primitive subsystem of $\Phi$ that strictly contains $\Psi$.

\vspace{2mm}\noindent For each $\Phi$ of exceptional type, we inspect all $W$-conjugacy classes of non-empty almost primitive subsystems $\Psi$ of $\Phi$ (which are listed in Table \ref{exceptlist}). We exclude those subsystems $\Psi$ that are maximal in $\Phi$. If $\Psi$ is maximal among $p$-closed subsystems of $\Phi$ for some $p \geq 0$ but $\Psi$ is not maximal in $\Phi$ then $p \neq 2$ and $(\Phi,\Psi)=(F_4,C_3A_1)$. Note that $C_3A_1$ is not almost $2$-primitive in $F_4$. So we also exclude the subsystem $C_3A_1$ of $F_4$. We list the resulting subsystems $\Psi$ in the second column (of Table \ref{primitivemaxactions}) and put their associated groups $\Stab_{\Iso(\Phi)}(\Delta)$ in the third column (using Table \ref{exceptlist} and Lemma \ref{isoautgps}). Note that $\smash{\widetilde{D_4}}$ is $p$-closed in $F_4$ if and only if $p=2$. All other subsystems $\Psi$ listed in column $2$ are $p$-closed for any $p \geq 0$.

\vspace{2mm}\noindent For each pair $(\Phi,\Psi)$ in columns $1$ and $2$, we list all $W$-conjugacy classes of almost primitive subsystems of $\Phi$ that strictly contain $\Psi$ in the fourth column. We can choose a representative $\Sigma$ of each of these conjugacy classes and a base $\Delta_{\Sigma}$ of $\Sigma$ such that $\Stab_{\Iso(\Phi)}(\Delta_{\Sigma})$ stabilises $\Psi$. Let $W_{\Sigma}$ be the Weyl group of $\Sigma$. Note that $\Psi^{\perp}=\Sigma^{\perp}=\varnothing$ by Lemma \ref{trivialperp}.

\vspace{2mm}\noindent For each triple $(\Phi, \Psi, \Sigma)$ listed in columns $1$, $2$ and $4$, we wish to compute the subgroup $\big(\!\Stab_{\Iso(\Phi)}(\Delta)\big)_{\Sigma}$ of $\Stab_{\Iso(\Phi)}(\Delta)$ that stabilises $\Sigma$. We observe that $\big(\!\Stab_{\Iso(\Phi)}(\Delta)\big)_{\Sigma}=\Stab_{N_{\Iso(\Phi)}(W_{\Sigma})}(\Delta)$. By a similar argument to Proposition $28$ of \cite{C}, there is an isomorphism from $\big(\!\Stab_{\Iso(\Phi)}(\Delta)\big)_{\Sigma}$ to $N_{N_{\Iso(\Phi)}(W_{\Sigma})}(W_{\Psi})\big/W_{\Psi}$. Hence, by the third isomorphism theorem, $\big(\!\Stab_{\Iso(\Phi)}(\Delta)\big)_{\Sigma}$ is an extension of $\Stab_{W_{\Sigma}}(\Delta)$ by $N_{N_{\Iso(\Phi)}(W_{\Sigma})}(W_{\Psi}) \big/N_{W_{\Sigma}}(W_{\Psi})$. Observe that $N_{N_{\Iso(\Phi)}(W_{\Sigma})}(W_{\Psi}) \geq \Stab_{\Iso(\Phi)}(\Delta_{\Sigma})$ since $\Stab_{\Iso(\Phi)}(\Delta_{\Sigma})$ stabilises $\Psi$ and that $\Stab_{\Iso(\Phi)}(\Delta_{\Sigma}) \cap W_{\Sigma}= 1$. So, using the second isomorphism theorem, we have  $$N_{N_{\Iso(\Phi)}(W_{\Sigma})}(W_{\Psi}) \big/N_{W_{\Sigma}}(W_{\Psi}) \geq \big(\!\Stab_{\Iso(\Phi)}(\Delta_{\Sigma}) N_{W_{\Sigma}}(W_{\Psi})\big) \big/N_{W_{\Sigma}}(W_{\Psi}) \cong \Stab_{\Iso(\Phi)}(\Delta_{\Sigma}).$$ Recall from Proposition $28$ of \cite{C} that $N_{\Iso(\Phi)}(W_{\Sigma}) \big/ W_{\Sigma} \cong \Stab_{\Iso(\Phi)}(\Delta_{\Sigma})$. It follows that $N_{N_{\Iso(\Phi)}(W_{\Sigma})}(W_{\Psi}) \big/N_{W_{\Sigma}}(W_{\Psi}) \cong \Stab_{\Iso(\Phi)}(\Delta_{\Sigma})$. In summary, we have shown that $$\big(\!\Stab_{\Iso(\Phi)}(\Delta)\big)_{\Sigma} = \Stab_{W_{\Sigma}}(\Delta) \hspace{0.5mm}.\hspace{0.5mm} \Stab_{\Iso(\Phi)}(\Delta_{\Sigma}).$$

\noindent The groups $\Stab_{W_{\Sigma}}(\Delta)$ and $\Stab_{\Iso(\Phi)}(\Delta_{\Sigma})$ can be found using Tables \ref{classlist}, \ref{exceptlist} and Lemma \ref{isoautgps}. It is then easy to find the conjugacy class of the subgroup $\big(\!\Stab_{\Iso(\Phi)}(\Delta)\big)_{\Sigma}$ of $\Stab_{\Iso(\Phi)}(\Delta)$ for all but the following difficult cases, which we consider individually.

\vspace{2mm}\noindent \underline{$(i)$}: Let $\Phi=E_7$ and $\Psi=(A_1)^7$. Here $\Stab_{\Iso(\Phi)}(\Delta) \cong \GL_3(2)$ and $\Sigma=D_6A_1$ or $D_4(A_1)^3$.

\vspace{2mm}\noindent Let $Y:=\Stab_{\Iso(\Phi)}(\Delta)$. There are two conjugacy classes $P_1(Y)$ and $P_2(Y)$ of order $24$ subgroups of $Y$, both of which are parabolic and isomorphic to $S_4$. Consider the action of $Y<S_7$ on the factors of $(A_1)^7$. Then $P_1(Y)$ fixes an $A_1$ whilst $P_2(Y)$ does not fix any. If $\Sigma=D_6A_1$ then $\Stab_W(\Delta_{\Sigma}) \cong 1$, $\Stab_{W_{\Sigma}}(\Delta) \cong S_4$ and hence $\big(\!\Stab_{\Iso(\Phi)}(\Delta)\big)_{\Sigma}=\Stab_{W_{\Sigma}}(\Delta)=P_1(Y)$ since it fixes an $A_1$. If $\Sigma=D_4(A_1)^3$ then $\Stab_W(\Delta_{\Sigma}) \cong S_3$ acts transitively on the factors of $(D_4)^{\perp}=(A_1)^3$ and $\Stab_{W_{\Sigma}}(\Delta) \cong (\Z_2)^2$ acts transitively on the factors of an $(A_1)^4$ in the $D_4$. Hence $\big(\!\Stab_{\Iso(\Phi)}(\Delta)\big)_{\Sigma}=(\Z_2)^2 \hspace{0.5mm}.\hspace{0.5mm} S_3=P_2(Y)$ since it does not fix an $A_1$.

\vspace{2mm}\noindent \underline{$(ii)$}: Let $\Phi=E_8$ and $\Psi=(A_1)^8$. Here $\Stab_{\Iso(\Phi)}(\Delta) \cong \AGL_3(2)$ and $\Sigma=E_7A_1$, $D_8$ or $(D_4)^2$.

\vspace{2mm}\noindent Let $X:=\Stab_{\Iso(\Phi)}(\Delta)$ and let $L = (\F_2)^3$ be the affine space on which $X$ acts. Then $X=L \rtimes Y$, where $Y=\GL_3(2)$ acts linearly on $L$ and $L \cong (\Z_2)^3$ as a group. There are two conjugacy classes of complements to $L$ in $X$. That is, there exists a maximal subgroup $Y_1$ of $X$ that is isomorphic to but not conjugate to $Y$. There are two conjugacy classes $L \rtimes P_1(Y)$ and $L \rtimes P_2(Y)$ of order $192$ subgroups of $X$, both of which are parabolic and isomorphic to each other.

\vspace{2mm}\noindent Consider the action of $X<S_8$ on the factors of $(A_1)^8$. Then $P_1(Y)$ fixes two $A_1$'s whilst $P_2(Y)$ fixes one $A_1$. In fact $L \rtimes P_2(Y)$ does not contain a $S_4$ subgroup that fixes two $A_1$'s. If $\Sigma=D_8$ then $\Stab_W(\Delta_{\Sigma}) \cong 1$ and $\Stab_{W_{\Sigma}}(\Delta) \cong (\Z_2)^3 \rtimes S_4$ contains a $S_4$ subgroup that fixes two $A_1$'s. 
Hence $\big(\!\Stab_{\Iso(\Phi)}(\Delta)\big)_{\Sigma}=\Stab_{W_{\Sigma}}(\Delta)=L \rtimes P_1(Y)$. If $\Sigma=(D_4)^2$ then $\Stab_{W_{\Sigma}}(\Delta) \cong (\Z_2)^4$ and $\Stab_W(\Delta_{\Sigma}) \cong \Z_2 \times S_3$. The $S_3$ in $\Stab_W(\Delta_{\Sigma})$ fixes precisely two $A_1$'s, but no non-trivial subgroup of $\Stab_{W_{\Sigma}}(\Delta)$ fixes both of these $A_1$'s. Hence $\big(\!\Stab_{\Iso(\Phi)}(\Delta)\big)_{\Sigma}=(\Z_2)^4 \hspace{0.5mm}.\hspace{0.5mm} (\Z_2 \times S_3)=L \rtimes P_2(Y)$ since it does not contain a $S_4$ subgroup that fixes two $A_1$'s.

\vspace{2mm}\noindent \underline{$(iii)$}: Let $\Phi=E_8$ and $\Psi=(A_2)^4$. Here $\Stab_{\Iso(\Phi)}(\Delta) \cong \GL_2(3)$ and $\Sigma=E_6A_2$.

\vspace{2mm}\noindent We have $\Stab_{W_{\Sigma}}(\Delta) \cong S_3$ and $\Stab_W(\Delta_{\Sigma}) \cong \Z_2$. All order $12$ subgroups of $\GL_2(3)$ are conjugate and are isomorphic to $\mathcal{D}i_{12}$. Hence $\big(\!\Stab_{\Iso(\Phi)}(\Delta)\big)_{\Sigma}=S_3 \hspace{0.5mm}.\hspace{0.5mm} \Z_2 \cong \mathcal{D}i_{12}$.

\vspace{2mm}\noindent At this point we have justified the entries in the fifth column of Table \ref{primitivemaxactions}.

\vspace{2mm}\noindent For each group $\Stab_{\Iso(\Phi)}(\Delta)$ in column $3$, we list -- in the last column -- all conjugacy classes of subgroups of $\Stab_{\Iso(\Phi)}(\Delta)$ that are not contained in any of the associated subgroups $\big(\!\Stab_{\Iso(\Phi)}(\Delta)\big)_{\Sigma}$ in column $5$. To do this, we construct the poset of conjugacy classes of subgroups for each group $\Stab_{\Iso(\Phi)}(\Delta)$ (this is easy to do unless $\Stab_{\Iso(\Phi)}(\Delta)$ is one of $\GL_2(3)$, $\GL_3(2)$ or $\AGL_3(2)$, in which case we use GAP). This completes our explanation of Table \ref{primitivemaxactions}.

\vspace{2mm}\noindent We now prove the second assertion of the lemma.

\vspace{2mm}\noindent By construction of Table \ref{primitivemaxactions}, there exists a $p$-closed $\Pi$-stable almost primitive subsystem $\Sigma$ of $\Phi$ that strictly contains $\Psi$ if and only if $\Pi$ is not one of the subgroups of $\Stab_{\Iso(\Phi)}(\Delta)$ listed in the last column of Table \ref{primitivemaxactions}. If $\Psi$ is strictly contained in a $\Pi$-stable proper subsystem of $\Phi$ then $\Psi$ is strictly contained in a $\Pi$-stable almost primitive subsystem of $\Phi$. So $\Psi$ is maximal among $p$-closed $\Pi$-stable subsystems of $\Phi$ if and only if $\Pi$ is one of the groups listed in the last column of Table \ref{primitivemaxactions}. 

\vspace{2mm}\noindent Let $(\mathcal{G}, \mathcal{H}, \theta)$ be associated to one of the triples $(\Phi, \Psi, \Pi )$ listed in columns $1$, $2$ and $6$ of Table \ref{primitivemaxactions}. It remains to show that $V^{\Pi}$ is trivial.


\vspace{2mm}\noindent If $r(\Psi)=r(\Phi)$ then $V$ is trivial by definition. We inspect Table \ref{primitivemaxactions}. The only case that satisfies $r(\Psi)<r(\Phi)$ is when $\Phi=E_6$ and $\Psi=D_4$. All possibilities for $\Pi \leq \Stab_{\Iso(\Phi)}(\Delta) \cong \Z_2 \times S_3$ that are listed in column $6$ contain an element $\tau$ of order $3$. Observe that $\tau$ acts on the $2$-dimensional subspace $V$ of $E$ as a rotation by an angle of $2\pi/3$. So there does not exist a non-trivial subspace of $V$ that is fixed pointwise by $\Pi$.
\end{proof}
\end{Lemma}

\noindent The following proposition characterises when $(\mathcal{G}, \mathcal{H}, \theta)$ is maximal.

\begin{Proposition}\label{pmaximalexceptchar} Let $\Phi$ be of exceptional type and let $r(\Phi)=\dim(E)$. Then $(\mathcal{G}, \mathcal{H}, \theta)$ is maximal if and only if one of the following occurs: \begin{itemize}	
\vspace{-1mm}\item $\Psi$ is maximal among $p$-closed subsystems of $\Phi$ and $\Psi$ has maximal rank in $\Phi$, 
\vspace{-2mm}\item $(\Phi,\Psi) =(E_6,D_5)$ or $(E_7,E_6)$ and $\Pi$ is non-trivial,
\vspace{-2mm}\item $(\Phi,\Psi,\Pi)$ is one of the triples listed in columns $1$, $2$ and $6$ of Table \ref{primitivemaxactions}, or
\vspace{-2mm}\item $\Psi=\varnothing$, $\Phi=E_6$, $E_7$ or $E_8$ and no non-empty proper closed subsystem of $\Phi$ is $\Pi$-stable.
\end{itemize} \vspace{-1mm} In particular, if $(\mathcal{G}, \mathcal{H}, \theta)$ is maximal then $\Psi$ is an almost $p$-primitive subsystem of $\Phi$. If in addition $(\mathcal{G}, \mathcal{H}, \theta)$ is isotropic then $\Psi$ is non-empty.
\begin{proof} Let $(\mathcal{G}, \mathcal{H}, \theta)$ be maximal. Assume (for a contradiction) that $\Omega$ is a $p$-closed $N_{\Iso(\Phi)}(W_{\Psi})$-stable proper subsystem of $\Phi$ that strictly contains $\Psi$. Observe that $\Omega$ is $W_{\Psi}$-stable and so $\Omega$ is $N_{\Iso(\Phi)}(W_{\Psi})\big/W_{\Psi}$-stable. It follows that $\Omega$ is $\Stab_{\Iso(\Phi)}(\Delta)$-stable (using the natural isomorphism $\Stab_{\Iso(\Phi)}(\Delta) \smash{\xrightarrow{\sim}} N_{\Iso(\Phi)}(W_{\Psi})\big/W_{\Psi}$ described in $\S \ref{apsubsystems}$). But then $(\mathcal{G}, \mathcal{H}, \theta)$ cannot be maximal since $\Pi \leq \Stab_{\Iso(\Phi)}(\Delta)$. We have our contradiction. Hence $\Psi$ is an almost $p$-primitive subsystem of $\Phi$.

\vspace{2mm}\noindent Henceforth let $\Psi=\varnothing$. Recall from Table \ref{exceptlist} that $\Psi$ is almost $p$-primitive in $\Phi$ if and only if $\Phi=E_6$, $E_7$ or $E_8$ (for any choice of $p$). If $\Phi=E_6$, $E_7$ or $E_8$ then all closed subsystems of $\Phi$ are $p$-closed since $\Phi$ is simply laced. Observe that $V=E$ (since $\Psi=\varnothing$ and $r(\Phi)=\dim(E)$). If $E^{\Pi}$ is non-trivial then, since $\Phi$ is irreducible, either $E^{\Pi} \cap \Phi$ or $(E^{\Pi})^{\perp} \cap \Phi$ (or both) is a non-empty $\Pi$-stable proper closed subsystem of $\Phi$. 
In summary, $(\mathcal{G}, \mathcal{H}, \theta)$ is maximal if and only if $\Phi=E_6$, $E_7$ or $E_8$ and no non-empty proper closed subsystem of $\Phi$ is $\Pi$-stable. Finally, if $(\mathcal{G}, \mathcal{H}, \theta)$ is isotropic then $E^{\Pi}$ is non-trivial and so $(\mathcal{G}, \mathcal{H}, \theta)$ cannot be maximal.

\vspace{2mm}\noindent The result then follows from Lemmas \ref{alssll} and \ref{maxactions}.
\end{proof}
\end{Proposition}

\section{Proof of Theorem \ref{maintheorem!}}\label{mainchapter}

\noindent The objective of this section is to prove Theorem \ref{maintheorem!}. We prove parts $(i)$, $(ii)$ and $(iii)$ of Theorem \ref{maintheorem!} in $\S \ref{indexHG}$, $\S \ref{maxconnectedness}$ and $\S \ref{part(iii)}$ respectively. 

\vspace{2mm}\noindent We begin with some setup.

\vspace{2mm}\noindent Let $k$ be any field and let $p =\Char(k)$. Let $\overline{k}$ be an algebraic closure of $k$ and let $K$ be the separable closure of $k$ in $\overline{k}$. Let $\Gamma=\Gal(K/k)$ be the absolute Galois group of $k$.

\vspace{2mm}\noindent Let $G$ be a connected reductive (algebraic) $k$-group. Let $S$ be a maximal $k$-split torus of $G$ and let $T$ be a maximal $k$-torus of $G$ that contains $S$. There is a natural action of $\Gamma$ on $G(K)$ (in $\S \ref{Boreltits}$ we refer to this action as $\star$). This induces an action $\iota:\Gamma \to \GL\!\big(X(T)_{\R}\big)$ given by $\iota(\sigma)(\chi):= \sigma \chi \sigma^{-1} $ for $\sigma \in \Gamma$ and $\chi \in X(T)$. We endow $X(T)_{\R}$ with an inner product $(\cdot \hspace{0.5mm},\cdot)$ that is invariant under the extended Weyl group of $G$ with respect to $T$. 

\vspace{2mm}\noindent Let $H$ be a connected reductive proper $k$-subgroup of maximal rank in $G$. Let $S_H$ be a maximal $k$-split torus of $H$ and let $T_H$ be a maximal $k$-torus of $H$ that contains $S_H$. Define $\iota_H:\Gamma \to \GL\!\big(X(T_H)_{\R}\big)$ by $\iota_H(\sigma)(\gamma):= \sigma \gamma \sigma^{-1} $ for $\sigma \in \Gamma$ and $\gamma \in X(T_H)$. 

\vspace{2mm}\noindent Theorem $20.9(ii)$ of \cite{B} states that all maximal $k$-split tori of $G$ are $G(k)$-conjugate to $S$. So there exists $g_1 \in G(k)$ such that $(S_H)^{g_1} \subseteq S$. Note that $(S_H)^{g_1}$ is a $k$-split torus. Consider the subgroup $L:=C_G\big((S_H)^{g_1}\big)$ of $G$. By Proposition \ref{kparabolicspairing2}, either $L=G$ or $L$ is a Levi $k$-subgroup of some parabolic $k$-subgroup of $G$. By Theorem $20.9(ii)$ of \cite{B}, there exists $g_2 \in L(K)$ such that $(T_H)^{g_2g_1}=T$. Denote $g:=g_2g_1$. Let $\theta_g:X(T_H)_{\R} \to X(T)_{\R}$ be defined by $\chi \mapsto \chi \cdot g^{-1}$. Using Corollary \ref{inducedisometry}, we endow $X(T_H)_{\R}$ with an inner product $(\cdot \hspace{0.5mm},\cdot)_H$ that is invariant under the extended Weyl group of $H$ with respect to $T_H$ and such that $\theta_g$ is a bijective isometry. For simplicity of notation, henceforth denote $\theta:=\theta_g$. 

\vspace{2mm}\noindent By Corollary \ref{gammalinorder}, there exists a $\Gamma$-order $<$ on $X(T)$ such that $<_{g^{-1}}$ is a $\Gamma$-order on $X(T_H)$. Let $\Lambda$ be the system of simple roots for $G$ with respect to $T$ that is compatible with $<$. Let $\Delta$ be the system of simple roots for $H$ with respect to $T_H$ that is compatible with $<_{g^{-1}}$. 

\vspace{2mm}\noindent We use the following notation. Let $\Lambda_0$ be the subset of $\Lambda$ that vanishes on $S$. Let $\Delta_0$ be the subset of $\Delta$ that vanishes on $S_H$. Let $\Phi:=\langle \Lambda \rangle$, $\Psi:=\langle \theta(\Delta) \rangle$ and $\Psi_0:=\langle \theta(\Delta_0) \rangle$. Let $W_{\Lambda}:=N_G(T)/T$, $W_{\Delta}:=N_H(T_H)/T_H$, $W_{\Lambda_0}:=N_{C_G(S)}(T)/T$ and $W_{\Delta_0}:=N_{C_H(S_H)}(T_H)/T_H$. It follows from Lemma \ref{propertiesofgammasystem}$(ii)$ that the subscript of each of these Weyl groups is indeed the associated root system.

\vspace{2mm}\noindent Let $\sigma \in \Gamma$. By Lemma \ref{propertiesofgammasystem}$(iii)$, there exists a unique $w_{\sigma} \in W_{\Lambda_0}$ such that $\big(w_{\sigma}\iota(\sigma)\big) (\Lambda)=\Lambda$. The Tits action $\hat{\iota}:\Gamma \to \GL\!\big(X(T)_{\R}\big)$ is defined by $\hat{\iota}(\sigma):=w_{\sigma}\iota(\sigma)$. The index of $G$ is the quadruple $\mathcal{I}(G)=\big(X(T)_{\R}, \Lambda, \Lambda_0,\hat{\iota}(\Gamma)\big)$. Similarly, there exists a unique $v_{\sigma} \in W_{\Delta_0}$ such that $\big(v_{\sigma}\iota_H(\sigma)\big) (\Delta)=\Delta$. The Tits action $\hat{\iota}_H:\Gamma \to \GL\!\big(X(T_H)_{\R}\big)$ is defined by $\hat{\iota}_H(\sigma):=v_{\sigma}\iota_H(\sigma)$. The index of $H$ is $\mathcal{I}(H)=\big(X(T_H)_{\R}, \Delta, \Delta_0, \hat{\iota}_H(\Gamma)\big)$.

\vspace{2mm}\noindent The \textit{embedding of indices} of $H \subset G$ (with respect to $T$, $T_H$, $<$ and $g$) is the triple $\big(\mathcal{I}(G),\mathcal{I}(H),\theta\big)$.

\begin{Remark}\label{partialindex} Alternatively, instead of defining the index of $G$ with respect to a maximal $k$-split torus, we could have proceeded as in Remark $2.7.2(d)$ of \cite{T} and defined the \textit{partial index} of $G$ with respect to $T_H$. However, this would lose a lot of detail in the classification.
\end{Remark} 

\subsection{Proof of Theorem \ref{maintheorem!}$(i)$}\label{indexHG}

\noindent \textbf{Theorem 1$(i)$.} \textit{The triple $\big(\mathcal{I}(G),\mathcal{I}(H),\theta\big)$ is a $p$-embedding of abstract indices.}

\begin{proof} We know that $\mathcal{I}(G)$ and $\mathcal{I}(H)$ are both abstract indices by Theorem \ref{indeximpo}, and that $\theta$ is a bijective isometry. It remains to check that $\big(\mathcal{I}(G),\mathcal{I}(H),\theta\big)$ satisfies all of the axioms of Definition \ref{mydefinition!}.

\vspace{2mm}\noindent We recall the following combinatorial notation associated to Definition \ref{mydefinition!}. Let $(X(T_H)_{\R})_a$ be the smallest subspace of $X(T_H)_{\R}$ that contains $\Delta_0$ and has an orthogonal complement that is fixed pointwise by $\hat{\iota}_H(\Gamma)$. Let $\Phi^+$ be the set of positive roots of $\Phi$ with respect to $\Lambda$. Let $\Lambda_a:= \Lambda \cap \theta\big((X(T_H)_{\R})_a\big)$, $\Phi_a:=\Phi \cap \theta\big((X(T_H)_{\R})_a\big)$ and let ${}_{in}\Lambda_a$ be the union of irreducible components of $\Lambda_a$ that is maximal with respect to the condition ${}_{in}\Lambda_a \subseteq \Psi_0$.

\vspace{2mm}\noindent We begin with the following result. 

\begin{Lemma}\label{prelim} $\Lambda_a$ is a system of simple roots for $L$ with respect to $T$. Moreover, $\Phi_a=\langle \Lambda_a \rangle$.
\begin{proof} Recall from Proposition \ref{abstractTits} that $(X(T_H)_{\R})_a$ is the subspace of $X(T_H)_{\R}$ that is generated by all characters that vanish on $S_H$. So, since $\theta$ is a bijective isometry, $\theta^{-1}(\Lambda_a)$ (resp. $\theta^{-1}(\Phi_a)$) is the subset of $\theta^{-1}(\Lambda)$ (res. $\theta^{-1}(\Phi)$) that vanishes on $S_H$. 

\vspace{2mm}\noindent By Corollary \ref{propertiesinducedcor} and Lemma \ref{preservesGammabase}, $\theta^{-1}(\Lambda)$ is a $\Gamma$-system of simple roots for $G$ with respect to $T_H$. Then, by Lemma \ref{mybabyfixed}, $\theta^{-1}(\Phi_a)$ is the root system of $C_G(S_H)$ with respect to $T_H$. Combining this with Lemma \ref{propertiesofgammasystem}$(i)$ and Remark \ref{niceactionremark} tells us that $\theta^{-1}(\Lambda_a)$ is a system of simple roots for $C_G(S_H)$ with respect to $T_H$. Once again applying Corollary \ref{propertiesinducedcor} gives us the result.
\end{proof}
\end{Lemma}

\noindent We check axiom $(A.1)$ of Definition \ref{mydefinition!} in the following lemma. 

\begin{Lemma}\label{A1proof} $\Phi$ is $\hat{\iota}_H(\Gamma)^{\theta}$-stable, $\Psi$ is a $p$-closed subsystem of $\Phi$, $\theta(\Delta) \subset \Phi^+$, $\Phi_a=\langle \Lambda_a \rangle$ and ${}_{in}\Lambda_a \subseteq \theta(\Delta_0)$.
\begin{proof} By Lemma \ref{propertiesinduced}, $\theta^{-1}(\Phi)$ is the root system of $G$ with respect to $T_H$. Observe that $\theta^{-1}(\Phi)$ is $\iota_H(\Gamma)$-stable by Lemma \ref{niceaction} and Remark \ref{niceactionremark}. Moreover, $W_{\Delta_0}<\smash{(W_{\Lambda})^{\theta^{-1}}}$ by Corollary \ref{propertiesinducedcor1}. So $\theta^{-1}(\Phi)$ is also $\hat{\iota}_H(\Gamma)$-stable since $v_{\sigma} \in W_{\Delta_0}$. Hence $\Phi$ is $\hat{\iota}_H(\Gamma)^{\theta}$-stable.

\vspace{2mm}\noindent By Corollary \ref{propertiesinducedcor}, $\Psi$ is the root system of $H^g$ with respect to $T$. Then, by Corollary 13.7 of \cite{MT}, $\Psi$ is a $p$-closed subsystem of $\Phi$. By construction $\theta(\Delta)$ is compatible with $<$, and so $\theta(\Delta) \subset \Phi^+$. We know from Lemma \ref{prelim} that $\Phi_a=\langle \Lambda_a \rangle$. It remains to show that ${}_{in}\Lambda_a \subseteq \theta(\Delta_0)$.

\vspace{2mm}\noindent Up to rearrangement, we can uniquely decompose $L$ as a commuting product $\big(\prod_{i=1}^m L_i\big) {Z(L)}^{\circ}$ where each $L_i$ is $k$-simple and ${Z(L)}^{\circ}$ is a $k$-torus. 
Let ${}_{in}L$ be the product of all $k$-simple components $L_i$ of $L$ such that $L_i \subseteq C_H(S_H)^{g_1}$ (if none exist then say ${}_{in}L:=1$).

\vspace{2mm}\noindent Recall from Lemmas \ref{mybabyfixed} and \ref{propertiesofgammasystem}$(i)$ that $\Delta_0$ is a system of simple roots for $C_H(S_H)$ with respect to $T_H$. Each $(L_i)^{g_1^{-1}}$ is a $k$-simple component of $C_H(S_H)$ since $g_1 \in G(k)$. So there exists a subset ${}_{in}\Delta_0$ of $\Delta_0$ that is a system of simple roots for $({}_{in}L)^{g_1^{-1}}T_H$ with respect to $T_H$. Then, by Corollary \ref{propertiesinducedcor}, $\theta(_{in}\Delta_0)$ is a system of simple roots for $({}_{in}L)T$ with respect to $T$.

\vspace{2mm}\noindent Let $i \in \{1,2,...,m\}$. It follows from Lemma \ref{prelim} that there exists a subset $(\Lambda_a)_i$ of $\Lambda_a$ that is a system of simple roots for $L_iT$ with respect to $T$. Observe that $(\Lambda_a)_i \subseteq \Psi_0$ if and only if $L_i \subseteq C_H(S_H)^{g_1}$. So ${}_{in}\Lambda_a$ is another system of simple roots for $({}_{in}L)T$ with respect to $T$. Since ${}_{in}\Lambda_a$ and $\theta({}_{in}\Delta_0)$ are both compatible with $<$, we have ${}_{in}\Lambda_a=\theta({}_{in}\Delta_0)$ by uniqueness.
 \end{proof}
\end{Lemma}

\noindent \noindent Let $W_{\Lambda_a}:=N_L(T)/T$. Note that $W_{\Lambda_a}$ is indeed the Weyl group of $\Lambda_a$ by Lemma \ref{prelim}. We now check axiom $(A.2)$ of Definition \ref{mydefinition!}. 

\begin{Lemma}\label{typeofG} For every $\sigma \in \Gamma$ there exists a unique $x_{\sigma} \in W_{\Lambda_a}$ such that $\hat{\iota}(\sigma)=x_{\sigma}\hat{\iota}_H(\sigma)^{\theta}$.
\begin{proof} Let $\sigma \in \Gamma$. We claim that $w_{\sigma}$ and $(v_{\sigma})^{\theta}$ are both in $W_{\Lambda_a}$.

\vspace{2mm}\noindent It follows from the proof of Lemma \ref{prelim} that $\Lambda_a$ is the subset of $\Lambda$ that vanishes on $(S_H)^{g_1}$. Then $\Lambda_0 \subseteq \Lambda_a$ since $(S_H)^{g_1} \subseteq S$. Hence $W_{\Lambda_0} \subseteq W_{\Lambda_a}$ and so $w_{\sigma} \in W_{\Lambda_a}$. Now observe that $\Psi_0 \subseteq \Phi_a$ using Corollary \ref{propertiesinducedcor} and Lemma \ref{prelim}. Then $(W_{\Delta_0})^{\theta} \leq W_{\Lambda_a}$ since $\theta$ is a bijective isometry. So indeed $(v_{\sigma})^{\theta} \in W_{\Lambda_a}$.

\vspace{2mm}\noindent Recall from Lemma \ref{A1proof} that $\Phi$ is $\hat{\iota}_H(\sigma)^{\theta}$-stable. Since $W_{\Lambda}$ acts simply transitively on the set of bases of $\Phi$, there exists a unique $x_{\sigma} \in W_{\Lambda}$ such that $\big(x_{\sigma}\hat{\iota}_H(\sigma)^{\theta}\big)(\Lambda)=\Lambda$. Let $n_{\sigma}:=g^{\sigma} g^{-1}$. Note that $n_{\sigma} \in L(K)$ since $g=g_2g_1$ where $g_1 \in G(k)$ and $g_2 \in L(K)$. Moreover, $n_{\sigma}$ normalises $T$ since $T_H$ and $T$ are both defined over $k$. Let $u_{\sigma}$ be the image of $n_{\sigma}$ in $W_{\Lambda_a}$ under the natural projection $N_L(T)(K) \to W_{\Lambda_a}$. Then $u_{\sigma}= \iota_H(\sigma)^{\theta}\iota(\sigma)^{-1}$ by Lemma \ref{functorialbehaviour}. Hence $$ \hat{\iota}(\sigma)=w_{\sigma}u_{\sigma}^{-1}\iota_H(\sigma)^{\theta} = w_{\sigma}u_{\sigma}^{-1}v_{\sigma}^{-\theta}\hat{\iota}_H(\sigma)^{\theta}.$$ Since $\hat{\iota}(\sigma)$ stabilises $\Lambda$, we have $x_{\sigma}=w_{\sigma}u_{\sigma}^{-1}v_{\sigma}^{-\theta}$ by uniqueness.
\end{proof}
\end{Lemma}

\noindent It follows immediately from Lemma \ref{typeofG} that axiom $(A.2)$ of Definition \ref{mydefinition!} is satisfied. Another consequence of Lemma \ref{typeofG} is that the index $\mathcal{I}(G)$ is of inner type if and only if $\hat{\iota}_H(\Gamma)^{\theta} \leq W_{\Lambda}$.

\vspace{2mm}\noindent In the following lemma we check axiom $(A.3)$ of Definition \ref{mydefinition!}.

\begin{Lemma}\label{A3proof} $\Lambda_a$ is a $\hat{\iota}(\Gamma)$-stable subset of $\Lambda$ that contains $\Lambda_0$.
\begin{proof} Recall from the proof of Lemma \ref{typeofG} that $\Lambda_a \supseteq \Lambda_0$. By Lemma \ref{prelim}, $\Phi_a$ is the root system of $L$ with respect to $T$. Since $L$ is defined over $k$, we see that $\Phi_a$ is $\iota(\Gamma)$-stable by Lemma \ref{niceaction}. Then $\Phi_a$ is $\hat{\iota}(\Gamma)$-stable since $w_{\sigma} \in W_{\Lambda_0} \subseteq W_{\Lambda_a}$. Hence $\Lambda_a = \Phi_a \cap \Lambda$ is $\hat{\iota}(\Gamma)$-stable. 
\end{proof}
\end{Lemma}

\noindent Finally, we check axiom $(A.4)$ of Definition \ref{mydefinition!}.

\begin{Lemma}\label{removeiha} ${}_{in}\Lambda_a$ is contained in $\Lambda_0$.
\begin{proof} Recall from the proof of Lemma \ref{A1proof} that ${}_{in}\Lambda_a$ is a system of simple roots for $({}_{in}L)T$ with respect to $T$. Observe that ${}_{in}L$ is $k$-anisotropic since $({}_{in}L)^{g_1^{-1}} \subseteq C_H(S_H)'$ and $g_1 \in G(k)$. Hence $S \cap {}_{in}L=1$ and so ${}_{in}\Lambda_a \subseteq \Lambda_0$ by Lemma \ref{mybabyfixed}.
\end{proof}
\end{Lemma}

\noindent By combining Lemmas \ref{A1proof}, \ref{typeofG}, \ref{A3proof} and \ref{removeiha} we have proved part $(i)$ of Theorem \ref{maintheorem!}. 
\end{proof}

\begin{Remark}\label{indexofL} Recall from $\S \ref{embeddingsofabstractindices}$ that $\mathcal{L}:=\big(X(T)_{\R}, \Lambda_a, \Lambda_0,\hat{\iota}(\Gamma)\big)$ is a combinatorially defined object called the $(\mathcal{I}(H),\theta)$-embedded subindex of $\mathcal{I}(G)$. Let $\mathcal{M}:=\big(X(T_H)_{\R}, \Delta_0, \Delta_0,\hat{\iota}_H(\Gamma)\big)$ be the (combinatorially defined) minimal subindex of $\mathcal{I}(H)$. It is worth noting that $(\mathcal{L},\mathcal{M},\theta)$ is the embedding of indices of $C_H(S_H) \subset L$ with respect to $T$, $T_H$, $<$ and $g$. That is, the embedding of indices of $C_H(S_H) \subset L$ is a combinatorial invariant associated to $\big(\mathcal{I}(G),\mathcal{I}(H),\theta\big)$.
\end{Remark}

\subsection{Proof of Theorem \ref{maintheorem!}$(ii)$}\label{maxconnectedness}

\noindent In this section we prove part $(ii)$ of Theorem \ref{maintheorem!}. 

\vspace{2mm}\noindent We first characterise when $H$ is maximal among connected reductive $k$-subgroups of $G$.

\begin{Lemma}\label{reductionlem} $H$ is maximal among connected reductive $k$-subgroups of $G$ if and only if $\Psi$ is maximal among $p$-closed $\hat{\iota}_H(\Gamma)^{\theta}$-stable subsystems of $\Phi$.
\begin{proof} We will repeatedly use the fact that $\theta$ is a bijective isometry. Recall from Theorem \ref{maintheorem!}$(i)$ that $\Psi$ is a $p$-closed subsystem of $\Phi$. By definition of the index $\mathcal{I}(H)$, $\Delta$ is $\hat{\iota}_H(\Gamma)$-stable. So $\Psi$ is $\hat{\iota}_H(\Gamma)^{\theta}$-stable.

\vspace{2mm}\noindent Let $M$ be a connected reductive proper $k$-subgroup of $G$ that strictly contains $H$. Let $\Phi_M$ be the root system of $M$ with respect to $T_H$. By Lemma \ref{propertiesinduced}, $\theta^{-1}(\Phi)$ is the root system of $G$ with respect to $T_H$. Then $\Phi_M$ is a $p$-closed proper subsystem of $\theta^{-1}(\Phi)$ that strictly contains $\langle \Delta \rangle$ by Corollary $13.7$ of \cite{MT}. Since $M$ is defined over $k$, by Lemma \ref{niceaction} and Remark \ref{niceactionremark}, we see that $\Phi_M$ is $\iota_H(\Gamma)$-stable. Hence $\Phi_M$ is $\hat{\iota}_H(\Gamma)$-stable (by construction of the Tits action $\hat{\iota}_H$, since $\Phi_M$ contains $\Delta_0$).

\vspace{2mm}\noindent Conversely, let $\Sigma$ be a $p$-closed $\hat{\iota}_H(\Gamma)^{\theta}$-stable proper subsystem of $\Phi$ that strictly contains $\Psi$. Then $\theta^{-1}(\Sigma)$ is a $p$-closed $\hat{\iota}_H(\Gamma)$-stable proper subsystem of $\theta^{-1}(\Phi)$ that strictly contains $\langle \Delta \rangle$. Note that $\theta^{-1}(\Sigma)$ must also be $\iota_H(\Gamma)$-stable since it contains $\Delta_0$. Hence, by Corollary \ref{downtown}, $M_{\theta^{-1}(\Sigma)}:=\big\langle T_H, U_{\alpha} \hspace{0.5mm}\big|\hspace{0.5mm} \alpha \in \theta^{-1}(\Sigma) \big\rangle$ is a connected reductive proper $k$-subgroup of $G$ that strictly contains $H$.
\end{proof}
\end{Lemma}

\noindent Let $V$ be the largest subspace of $X(T)_{\R}$ that is contained in the span of $\Phi$ and is perpendicular to all roots in $\Psi$. Let $V^{\Gamma}$ denote the fixed point subspace of $V$ under the action of $\hat{\iota}_H(\Gamma)^{\theta}$. Consider the $k$-subtorus $R:=(Z(H) \cap G')^{\circ}$ of $T_H$. 

\begin{Lemma}\label{reductionlem1} $R_s$ is trivial if and only if $V^{\Gamma}$ is trivial.
\begin{proof} Consider the character space $X(T_H)_{\R}$ and the identification used in Lemma \ref{niceactioncor}. Recall from Lemma \ref{propertiesinduced} that $\theta^{-1}(\Phi)$ is the root system of $G$ with respect to $T_H$. By Lemma \ref{mybabyfixed}, $X_0(Z(H)^{\circ})_{\R}$ (resp. $X_0(Z(G)^{\circ})_{\R}$) is the subspace of $X(T_H)_{\R}$ that is spanned by $\Delta$ (resp. $\theta^{-1}(\Phi)$). 

\vspace{2mm}\noindent Observe that $T_H=(T_H\cap G')Z(G)^{\circ}$ is a commuting product. Then $X(T_H)_{\R} = X(T_H\cap G')_{\R} \oplus X(Z(G)^{\circ})_{\R}$ by Corollary \ref{niceactioncor1}. This decomposition must be orthogonal since $W_H$ centralises $Z(G)^{\circ}$, where $W_H$ denotes the Weyl group of $G$ with respect to $T_H$. Hence $\theta^{-1}(V)=X(Z(H)^{\circ})_{\R} \cap X(T_H \cap G')_{\R}$. Then $\theta^{-1}(V^{\Gamma})=\theta^{-1}(V) \cap X(S_H)_{\R}$ by Lemma \ref{Titsactionproperties}.

\vspace{2mm}\noindent By definition, $R_s=\big(S_H \cap Z(H) \cap G' \big)^{\circ}$. Then, by Corollary \ref{niceactioncor1}$(ii)$, we have $$X(R_s)_{\R}=X(S_H)_{\R} \cap X(Z(H)^{\circ})_{\R} \cap X(T_H \cap G')_{\R}=\theta^{-1}(V^{\Gamma}).$$ Hence $R_s$ is trivial if and only if $V^{\Gamma}$ is trivial.
\end{proof}
\end{Lemma}

\noindent Recall that $\big(\mathcal{I}(G),\mathcal{I}(H),\theta\big)$ is \textit{maximal} if $V^{\Gamma}$ is trivial and $\Psi$ is maximal among $p$-closed $\hat{\iota}_H(\Gamma)^{\theta}$-stable subsystems of $\Phi$. We are now able to prove part $(ii)$ of Theorem \ref{maintheorem!}.

\vspace{2mm}\noindent \textbf{Theorem 1$(ii)$.} \textit{$H$ is maximal connected in $G$ if and only if $\big(\mathcal{I}(G),\mathcal{I}(H),\theta\big)$ is maximal.}

\begin{proof} Let $\big(\mathcal{I}(G),\mathcal{I}(H),\theta\big)$ be maximal. Assume (for a contradiction) that $M$ is a maximal connected proper $k$-subgroup of $G$ that strictly contains $H$. Then $M$ is not reductive by Lemma \ref{reductionlem}. Since $M$ has maximal rank in $G$, $M$ must be parabolic by Corollary $3.3$ of \cite{BT} (Borel-Tits). 

\vspace{2mm}\noindent We construct a Levi $k$-subgroup of $M$ that contains $H$. 
Let $L_0=(L_0)'{Z(L_0)}^{\circ}$ be a Levi $k$-subgroup of $M$. Consider the associated quotient map $\pi:M= R_u(M) \rtimes L_0 \to L_0$ and its restriction $\pi|_H:H \to L_0$. Recall that $T_H$ is a maximal $k$-torus of $H$. Then $\pi(T_H)$ is a maximal $k$-torus of $L_0$ and so ${Z(L_0)}^{\circ} \leq \pi(T_H)$. Then $\pi|_H^{-1}\big({Z(L_0)}^{\circ}\big)$ is a central $k$-torus of $H$ 
and so, by Proposition \ref{kparabolicspairing2}, $L_1:=C_G\big(\pi|_H^{-1}{\big({Z(L_0)}^{\circ}\big)}{}_s\big)$ is a Levi $k$-subgroup of $M$ that contains $H$. Note that $L_1 = C_G\big({Z(L_1)}^{\circ}_s\big)$ by Proposition $20.6(i)$ of \cite{B}. Hence ${Z(L_1)}^{\circ} \cap G'$ is $k$-isotropic (as otherwise $L_1 = G$). 

\vspace{2mm}\noindent By construction, we have ${Z(L_1)}^{\circ} \cap G' \subseteq R \subseteq T_H$. So $R$ is $k$-isotropic. Hence $V^{\Gamma}$ is non-trivial by Lemma \ref{reductionlem1}. This is a contradiction.

\vspace{2mm}\noindent Conversely, let $H$ be maximal connected in $G$. Then $\Psi$ is maximal among $p$-closed $\hat{\iota}_H(\Gamma)^{\theta}$-stable subsystems of $\Phi$ by Lemma \ref{reductionlem}. Assume (for a contradiction) that $V^{\Gamma}$ is non-trivial. Then $R$ is $k$-isotropic by Lemma \ref{reductionlem1}. By Proposition \ref{kparabolicspairing2}, $C_G(R_s)$ is a Levi $k$-subgroup of some parabolic $k$-subgroup $P$ of $G$. That is, $G \supset P \supset C_G(R_s) \supseteq H$. All parabolic subgroups are connected. So we have a contradiction.
\end{proof}

\noindent This completes the proof of part $(ii)$ of Theorem \ref{maintheorem!}.

\subsection{Proof of Theorem \ref{maintheorem!}$(iii)$}\label{part(iii)}

\noindent In this section we prove part $(iii)$ of Theorem \ref{maintheorem!}.

\vspace{2mm}\noindent Let $\smash{'}\hspace{-0.4mm}H$ be a connected reductive $k$-subgroup of maximal rank in $G$. We use an identical setup as we did for $H$.

\vspace{2mm}\noindent Let $\smash{'}\hspace{-0.4mm}S_H$ be a maximal $k$-split torus of $\smash{'}\hspace{-0.4mm}H$ and let $\smash{'}T_H$ be a maximal $k$-torus of $\smash{'}\hspace{-0.4mm}H$ that contains $\smash{'}\hspace{-0.4mm}S_H$. Define $\smash{'}\iota_H:\Gamma \to \GL\!\big(\smash{'}X(T_H)_{\R}\big)$ by $\smash{'}\iota_H(\sigma)(\gamma):= \sigma \gamma \sigma^{-1} $ for $\sigma \in \Gamma$ and $\gamma \in X(\smash{'}T_H)$. There exists $h_1 \in G(k)$ such that $(\smash{'}\hspace{-0.4mm}S_H)^{h_1} \subseteq S$. Denote $\smash{'}\hspace{-0.4mm}L:=C_G\big((\smash{'}\hspace{-0.4mm}S_H)^{h_1}\big)$. There exists $h_2 \in \smash{'}\hspace{-0.4mm}L(K)$ such that $(\smash{'}T_H)^{h_2h_1}=T$. Denote $h:=h_2h_1$. Consider the induced map $\theta_h:X(\smash{'}T_H)_{\R} \to X(T)_{\R}$ defined by $\chi \mapsto \chi \cdot h^{-1}$. Using Corollary \ref{inducedisometry}, we endow $X(\smash{'}T_H)_{\R}$ with an inner product $'(\cdot \hspace{0.5mm},\cdot)_H$ that is invariant under the extended Weyl group of $\smash{'}\hspace{-0.4mm}H$ with respect to $\smash{'}T_H$ and such that $\theta_h$ is a bijective isometry. Henceforth, for simplicity, denote $\smash{'}\hspace{-0.3mm}\theta:=\theta_h$.

\begin{Remark}\label{compatibleity} It is not immediately obvious that we can adjust our choice of $h$ to ensure that $<_{h^{-1}}$ is a $\Gamma$-order on $X(\smash{'}T_H)$. However, this does not matter for the following reason. By Corollary \ref{gammalinorder}, there exists a $\Gamma$-order $\smash{'}\hspace{-1.8mm}<$ on $X(T)$ such that $\smash{'}\hspace{-1.8mm}<_{h^{-1}}$ is a $\Gamma$-order on $X(\smash{'}T_H)$. Let $\smash{'}\mathcal{I}(G)$ be the index of $G$ with respect to $T$ and $\smash{'}\hspace{-1.8mm}<$. Let $W_{\Gamma}:=N_{N_G(S)}(T)/T$. It follows from Lemma \ref{propertiesofgammasystem} that there exists $w \in W_{\Gamma}$ such that $w\big(\smash{'}\mathcal{I}(G)\big)=\mathcal{I}(G)$. 
So, in the proof of Theorem \ref{maintheorem!}$(iii)$, we simply adjust our conjugating element by $w$.
\end{Remark}

\noindent By Remark \ref{compatibleity}, we can assume without loss of generality that $<_{h^{-1}}$ is a $\Gamma$-order on $X(\smash{'}T_H)$. 

\vspace{2mm}\noindent Let $\smash{'}\hspace{-0.4mm}\Delta$ be the system of simple roots for $\smash{'}\hspace{-0.4mm}H$ with respect to $\smash{'}T_H$ that is compatible with $<_{h^{-1}}$. Let $\smash{'}\hspace{-0.4mm}\Delta_0$ be the subset of $\smash{'}\hspace{-0.4mm}\Delta$ that vanishes on $\smash{'}\hspace{-0.4mm}S_H$. Let $W_{\smash{'}\hspace{-0.4mm}\Delta}:=N_{\smash{'}\hspace{-0.6mm}H}(\smash{'}T_H)/\smash{'}T_H$ and $W_{\smash{'}\hspace{-0.4mm}\Delta_0}:=N_{C_{\smash{'}\hspace{-0.6mm}H}(\smash{'}\hspace{-0.4mm}S_H)}(\smash{'}T_H)/\smash{'}T_H$. There exists a unique $\smash{'}\hspace{-0.3mm}v_{\sigma} \in W_{\smash{'}\hspace{-0.4mm}\Delta_0}$ such that $\big(\smash{'}\hspace{-0.3mm}v_{\sigma}\smash{'}\hspace{-0.2mm}\iota_H(\sigma)\big) (\smash{'}\hspace{-0.4mm}\Delta)=\smash{'}\hspace{-0.4mm}\Delta$. The Tits action $\smash{'}\hat{\iota}_H:\Gamma \to \GL\!\big(\smash{'}X(T_H)_{\R}\big)$ is defined by $\smash{'}\hat{\iota}_H(\sigma):=\smash{'}\hspace{-0.3mm}v_{\sigma}\smash{'}\hspace{-0.2mm}\iota_H(\sigma)$ for each $\sigma \in \Gamma$. Let $\mathcal{I}(\smash{'}\hspace{-0.4mm}H)= \big(X(\smash{'}T_H)_{\R}, \smash{'}\hspace{-0.4mm}\Delta, \smash{'}\hspace{-0.4mm}\Delta_0, \smash{'}\hat{\iota}_H(\Gamma) \big)$.

\vspace{2mm}\noindent Recall that $H$ is \textit{index-conjugate} to $\smash{'}\hspace{-0.4mm}H$ in $G$ if there exists $x \in G(K)$ such that $H^x=\smash{'}\hspace{-0.4mm}H$, $(S_H)^x$ is a maximal $k$-split torus of $\smash{'}\hspace{-0.4mm}H$ and $x^{-1}x^{\Gamma} \subset H$. 

\vspace{2mm}\noindent We now restate and prove part $(iii)$ of Theorem \ref{maintheorem!}.

\vspace{2mm}\noindent \textbf{Theorem 1$(iii)$.} \textit{$H$ is index-conjugate to $\smash{'}\hspace{-0.4mm}H$ in $G$ if and only if $\big(\mathcal{I}(H),\theta \big)$ is conjugate to $\big(\mathcal{I}(\smash{'}\hspace{-0.4mm}H),\smash{'}\hspace{-0.3mm}\theta\big)$ in $\mathcal{I}(G)$.}

\begin{proof} Let $H$ be index-conjugate to $\smash{'}\hspace{-0.4mm}H$ in $G$. That is, there exists $x \in G(K)$ such that $H^x=\smash{'}\hspace{-0.4mm}H$, $(S_H)^x$ is a maximal $k$-split torus of $\smash{'}\hspace{-0.4mm}H$ and $x^{-1}x^{\Gamma} \subset H$. 

\vspace{2mm}\noindent By Theorem $20.9(ii)$ of \cite{B}, there exists $x_1 \in H(k)$ such that $(S_H)^{xx_1}=\smash{'}\hspace{-0.4mm}S_H$. Note that $C_H(S_H)$ is $K$-split (Theorem \ref{Grothendieck}). So, again by Theorem $20.9(ii)$ of \cite{B}, there exists $x_2 \in C_H(S_H)(K)$ such that $(T_H)^{xx_1x_2}=\smash{'}T_H$. So $c:=xx_1x_2 \in H(K)$ satisfies $(S_H)^c=\smash{'}\hspace{-0.4mm}S_H$ and $(T_H)^c=\smash{'}T_H$.

\vspace{2mm}\noindent By Corollary \ref{propertiesinducedcor} and Lemma \ref{preservesGammabase}, $\theta_c^{-1}(\smash{'}\hspace{-0.4mm}\Delta)$ is a $\Gamma$-system of simple roots for $H$ with respect to $T_H$. Denote $(W_{\Delta})_{\Gamma}:=N_{N_H(S_H)}(T_H)/T_H$. By Lemma \ref{propertiesofgammasystem}$(iv)$, there exists a unique $w_0 \in (W_{\Delta})_{\Gamma}$ such that $w_0(\Delta)=\theta_c^{-1}(\smash{'}\hspace{-0.4mm}\Delta)$. Let $d \in N_{N_H(S_H)}(T_H)(K)$ be a preimage of $w_0$ under the natural projection $N_{N_H(S_H)}(T_H)(K) \to (W_{\Delta})_{\Gamma}$ and let $y:=cd^{-1} \in H(K)$. Then $\theta_y(\Delta)=\smash{'}\hspace{-0.4mm}\Delta$ by Lemma \ref{functorialbehaviour}.

\vspace{2mm}\noindent Observe that $C_H(S_H)^y=C_{\smash{'}\hspace{-0.4mm}H}(\smash{'}\hspace{-0.4mm}S_H)$ since $H^y=\smash{'}\hspace{-0.4mm}H$ and $(S_H)^y=\smash{'}\hspace{-0.4mm}S_H$. Recall from Lemmas \ref{mybabyfixed} and \ref{propertiesofgammasystem}$(i)$ that $\Delta_0$ is a system of simple roots for $C_H(S_H)$ with respect to $T_H$. Similarly, $\smash{'}\hspace{-0.4mm}\Delta_0$ is a system of simple roots for $C_{\smash{'}\hspace{-0.4mm}H}(\smash{'}\hspace{-0.4mm}S_H)$ with respect to $\smash{'}T_H$. Combining this with Corollary \ref{propertiesinducedcor} and $\theta_y(\Delta)=\smash{'}\hspace{-0.4mm}\Delta$ gives us $\theta_y(\Delta_0)=\smash{'}\hspace{-0.4mm}\Delta_0$.

\vspace{2mm}\noindent Let $\sigma \in \Gamma$. Recall from Lemma \ref{niceaction} that $\langle \Delta \rangle$ is stabilised by $\iota_H(\sigma)$. Similarly, $\langle\smash{'}\hspace{-0.4mm}\Delta_0 \rangle$ is stabilised by $\smash{'}\iota_H(\sigma)$. Since $W_{\Delta}$ acts simply transitively on the set of bases of $\langle \Delta \rangle$, there is a unique element $\kappa_{\sigma} \in W_{\Delta}$ such that $\big(\kappa_{\sigma}\theta_y^{-1} \smash{'}\iota_H(\sigma)\theta_y \big)(\Delta) = \iota_H(\sigma) (\Delta)$. Note that $\theta_yW_{\Delta}\theta_y^{-1}=W_{\smash{'}\hspace{-0.4mm}\Delta}$ since $\theta_y(\Delta)=\smash{'}\hspace{-0.4mm}\Delta$. So $\kappa_{\sigma}=v_{\sigma}^{-1}\theta_y^{-1}\smash{'}\hspace{-0.3mm}v_{\sigma}\theta_y$ by uniqueness. Observe that $y^{-1}y^{\sigma} \in N_H(T_H)(K)$. Hence $\theta_y^{-1} \smash{'}\iota_H(\sigma)\theta_y \iota_H(\sigma)^{-1} \in W_{\Delta}$ by Lemma \ref{functorialbehaviour}. But then $\kappa_{\sigma}\theta_y^{-1} \smash{'}\iota_H(\sigma)\theta_y \iota_H(\sigma)^{-1}=1$ as only the trivial element of $W_{\Delta}$ stabilises $\Delta$. Rearranging, we have $\theta_y \hat{\iota}_H(\sigma)\theta_y^{-1}=\smash{'}\hat{\iota}_H(\sigma)$.

\vspace{2mm}\noindent We have thus far shown that $\theta_y\big(\mathcal{I}(H)\big)=\mathcal{I}(\smash{'}\hspace{-0.4mm}H)$. Let $w_1 \in W_{\Lambda}$ be the image of $gy^{-1}h^{-1} \in N_G(T)(K)$ under the natural projection $N_G(T)(K) \to W_{\Lambda}$. By Lemma \ref{functorialbehaviour}, we identify $w_1=\smash{'}\hspace{-0.3mm}\theta\theta_y\theta^{-1}$. Then $w_1\theta\big(\mathcal{I}(H)\big)=\smash{'}\hspace{-0.3mm}\theta\theta_y\big(\mathcal{I}(H)\big)=\smash{'}\hspace{-0.3mm}\theta\big(\mathcal{I}(\smash{'}\hspace{-0.4mm}H)\big)$. This completes one direction of the proof.

\vspace{2mm}\noindent Conversely, let $\big(\mathcal{I}(H),\theta \big)$ be conjugate to $\big(\mathcal{I}(\smash{'}\hspace{-0.4mm}H),\smash{'}\hspace{-0.3mm}\theta\big)$ in $\mathcal{I}(G)$. That is, there exists $w \in W_{\Lambda}$ such that $w\theta\big(\mathcal{I}(H)\big)= \smash{'}\hspace{-0.3mm}\theta\big(\mathcal{I}(\smash{'}\hspace{-0.4mm}H)\big)$. Since index-conjugacy is coarser than $G(k)$-conjugacy, it suffices to show that $H^{g_1}$ is index-conjugate to $\smash{'}\hspace{-0.4mm}H^{h_1}$ in $G$. Equivalently, without loss of generality, we can assume that $g_1$ and $h_1$ are trivial. 

\vspace{2mm}\noindent Consider the index $\mathcal{I}(H)$. Recall from Proposition \ref{abstractTits} that $(X(T_H)_{\R})_a$ is the subspace of $X(T_H)_{\R}$ that is generated by all characters that vanish on $S_H$. Since $g$ centralises $S_H$, $\theta\smash{\big((X(T_H)_{\R})_a\big)}$ is the subspace of $X(T)_{\R}$ that is generated by all characters that vanish on $S_H$. Similarly, $\smash{'}\hspace{-0.3mm}\theta\smash{\big((X(\smash{'}T_H)_{\R})_a\big)}$ is the subspace of $X(T)_{\R}$ that is generated by all characters that vanish on $\smash{'}\hspace{-0.4mm}S_H$. We know that $w\theta\smash{\big((X(T_H)_{\R})_a\big)}=\smash{'}\hspace{-0.3mm}\theta\smash{\big((X(\smash{'}T_H)_{\R})_a\big)}$ by Lemma \ref{conjlemma}. 

\vspace{2mm}\noindent Let $n \in N_G(T)(K)$ be a preimage of $w$ under the natural projection $N_G(T)(K) \to W_{\Lambda}$. Let $t \in \smash{'}\hspace{-0.4mm}S_H$ and let $\chi \in X(T)$ vanish on and only on $S_H$. Certainly such a $\chi$ exists. Then $0=w(\chi)(t)=\chi(t^n)$ and so $t^n \in S_H$. The reverse argument shows that $s^{n^{-1}} \in \smash{'}\hspace{-0.4mm}S_H$ for $s \in S_H$. Hence $(S_H)^{n^{-1}}=\smash{'}\hspace{-0.4mm}S_H$. 

\vspace{2mm}\noindent Let $z:=h^{-1}n^{-1}g \in G(K)$. Then $(S_H)^z=\smash{'}\hspace{-0.4mm}S_H$ and $(T_H)^z=\smash{'}T_H$. Recall that $\theta_z:X(T_H)_{\R} \to X(\smash{'}T_H)_{\R}$ is defined by $\chi \mapsto \chi \cdot z^{-1}$. By Lemma \ref{functorialbehaviour}, we identify $\theta_z=\smash{'}\hspace{-0.3mm}\theta^{-1} w\theta$. So $\theta_z\big(\mathcal{I}(H)\big)=\mathcal{I}(\smash{'}\hspace{-0.4mm}H)$. In particular, $\theta_z(\Delta)=\smash{'}\hspace{-0.4mm}\Delta$ and hence $H^z=\smash{'}\hspace{-0.4mm}H$ by Corollary \ref{propertiesinducedcor}.

\vspace{2mm}\noindent Let $\sigma \in \Gamma$. Observe that $\theta_zW_{\Delta_0}\theta_z^{-1}=W_{\smash{'}\hspace{-0.4mm}\Delta_0}$ since $\theta_z(\Delta_0)=\smash{'}\hspace{-0.4mm}\Delta_0$. Let $u_{\sigma}:=\theta_z^{-1}\hspace{0.1mm}\smash{'}\hspace{-0.3mm}v_{\sigma}^{-1} \theta_z v_{\sigma} \in W_{\Delta_0}$. By assumption, $\theta_z\hat{\iota}_H(\sigma)\theta_z^{-1} = \smash{'}\hat{\iota}_H(\sigma)$. Rearranging this equation gives us \begin{equation}\label{asfasf} u_{\sigma}=\theta_z^{-1} \smash{'}\iota_H(\sigma) \theta_z \iota_H(\sigma)^{-1}.\end{equation} Let $n_{\sigma}$ be a preimage of $u_{\sigma}$ under the natural projection $N_{C_H(S_H)}(T_H) \to W_{\Delta_0}$. Note that $n_{\sigma}z^{-1}z^{\sigma} \in N_G(T_H)$. Let $W_H$ be the Weyl group of $G$ with respect to $T_H$. Applying Lemma \ref{functorialbehaviour} to Equation $(\ref{asfasf})$ tells us that the image of $n_{\sigma}z^{-1}z^{\sigma}$ in $W_H$ under the natural projection is the trivial element. So $n_{\sigma}z^{-1}z^{\sigma} \in T_H$ and hence $z^{-1}z^{\sigma} \in H$. We are done since any maximal $k$-split torus of $H$ is $H(k)$-conjugate to $S_H$.
\end{proof}

\noindent This completes the proof of part $(iii)$ and hence the entirety of Theorem \ref{maintheorem!}.

\section{Proof of Theorem \ref{classificationthm!}}\label{restrictions}

\noindent In this section we describe the procedure which we use to prove Theorem \ref{classificationthm!}. This involves constructing Tables \ref{G_2}, \ref{F_4}, \ref{E_6}, \ref{E_7} and \ref{E_8} for the cases $G_2$, $F_4$, $E_6$, $E_7$ and $E_8$ respectively. We present these tables and the necessary computations in the Appendix.

\vspace{2mm}\noindent We first recall the statement of Theorem \ref{classificationthm!}. Let $\mathbb{P}$ denote the set of prime numbers. 

\begin{classificationthm!} \textit{For each $p \in \mathbb{P} \cup \{0\}$ and each (isomorphism class of) abstract index $\mathcal{G}$ of exceptional type, all $\mathcal{G}$-conjugacy classes of isotropic maximal $p$-embeddings of abstract indices are classified in Tables \ref{G_2}, \ref{F_4}, \ref{E_6}, \ref{E_7} and \ref{E_8} for the cases where $\mathcal{G}$ is of type $G_2$, $F_4$, $E_6$, $E_7$ and $E_8$ respectively.}
\end{classificationthm!}

\noindent Given a $p$-embedding of abstract indices $(\mathcal{G},\mathcal{H},\theta)$, we follow Remark \ref{abusenotation} and use the notation of $\S \ref{embeddingsofabstractindices}$. That is, $\mathcal{G}=(E,\Lambda,\Lambda_0, \Xi)$, $\mathcal{H}=(E,\Delta,\Delta_0,\Pi)$, $\Phi:=\langle \Lambda \rangle$, $\Phi_0:=\langle \Lambda_0 \rangle$, $W:=W_{\Lambda}$, $\Psi:=\langle \Delta \rangle$, $\Psi_0:=\langle \Delta_0 \rangle$, $E_s$ is the fixed point subspace of $\smash{(E_{\Delta_0})^{\perp}}$ under $\Pi$, $r_s(\mathcal{H}):=\dim (E_s)$, $\Phi_s := \Phi \cap E_s$, $E_a:=(E_s)^{\perp}$, $\Phi_a := \Phi \cap E_a$, $\overline{E_a}:=\smash{(E_{\Delta_0})^{\perp}} \cap E_a$, ${}_{in}\Phi_a$ is the union of irreducible components of $\Phi_a$ that is maximal with respect to ${}_{in}\Phi_a \subseteq \Psi_0$ and $\Iso(\Phi)$ is the isometry group of $E$ that stabilises $\Phi$. By assumption, $\dim E = r(\Phi)$.

\begin{proof} Let $p \in \mathbb{P} \cup \{0\}$ and let $\Phi$ be an irreducible root system of exceptional type. We do not require that $p$ and $\Phi$ are fixed, rather we consider all of their possible values. For every abstract index $\mathcal{G}$ of type $\Phi$, we describe an algorithm which we use to classify all $\mathcal{G}$-conjugacy classes of isotropic maximal $p$-embeddings of abstract indices $(\mathcal{H},\theta)$. We present our results in Table \ref{G_2} (resp. \ref{F_4}, \ref{E_6}, \ref{E_7}, \ref{E_8}) for the case where $\Phi=G_2$ (resp. $F_4$, $E_6$, $E_7$, $E_8$).

\vspace{2mm}\noindent By Proposition \ref{pmaximalexceptchar}, $\Delta$ is a non-empty almost primitive subsystem (of simple roots) of $\Phi$. As such, in columns $1$ and $2$ respectively (of Tables \ref{G_2}, \ref{F_4}, \ref{E_6}, \ref{E_7} and \ref{E_8}), we use Table \ref{exceptlist} and Lemma \ref{isoautgps} to list all $W$-conjugacy classes of non-empty almost primitive subsystems $\Delta$ of $\Phi$ along with their associated groups $\Stab_{\Iso(\Phi)}(\Delta)$.

\vspace{2mm}\noindent The property of $(\mathcal{H},\Phi)$ being maximal depends only on the $W$-conjugacy class of $\Delta$ in $\Phi$ and the conjugacy class of $\Pi$ in $\Stab_{\Iso(\Phi)}(\Delta)$. So, for each $\Delta$ in column $1$, we use Proposition \ref{pmaximalexceptchar} to list -- in column $3$ -- all conjugacy classes of subgroups $\Pi$ of $\Stab_{\Iso(\Phi)}(\Delta)$ for which $(\mathcal{H},\Phi)$ is maximal. We then use Tits' classification of isomorphism classes of abstract indices in Table II and $\S 3.1.2$ of \cite{T} to list -- in column $4$ -- all conjugacy classes in $\Phi$ of abstract indices $\mathcal{H}$ that are associated to the pair $(\Delta,\Pi)$. There is a subtlety here (see Lemma \ref{trickybitty}) in that certain abstract indices of type ${}^1\hspace{-0.5mm}D_n$ may be oriented in two different ways in $\Phi$ (or three if $n=4$). For example, for the case $\Phi=E_8$ and $\Delta=D_8$, the abstract indices $\begin{tikzpicture}[scale=.175,baseline=0.3ex]
    \foreach \y in {0,...,5}
    \draw[thin,xshift=\y cm] (\y cm,0) ++(.3 cm, 0) -- +(14 mm,0);
    \draw[thin] (0 cm,0) circle (3 mm);
    \draw[thin] (2 cm,0) circle (3 mm);
    \draw[thin,fill=black] (4 cm,0) circle (3 mm);
    \draw[thin] (6 cm,0) circle (3 mm);
    \draw[thin,fill=black] (8 cm,0) circle (3 mm);
    \draw[thin] (10 cm,0) circle (3 mm);
    \draw[thin,fill=black] (12 cm,0) circle (3 mm);
    \draw[thin,fill=black] (2 cm,2 cm) circle (3 mm);
    \draw[thin,fill=black] (2 cm, 3mm) -- +(0, 1.4 cm);
  \end{tikzpicture}$ and $\begin{tikzpicture}[scale=.175,baseline=0.3ex]
    \foreach \y in {0,...,5}
    \draw[thin,xshift=\y cm] (\y cm,0) ++(.3 cm, 0) -- +(14 mm,0);
    \draw[thin,fill=black] (0 cm,0) circle (3 mm);
    \draw[thin] (2 cm,0) circle (3 mm);
    \draw[thin,fill=black] (4 cm,0) circle (3 mm);
    \draw[thin] (6 cm,0) circle (3 mm);
    \draw[thin,fill=black] (8 cm,0) circle (3 mm);
    \draw[thin] (10 cm,0) circle (3 mm);
    \draw[thin,fill=black] (12 cm,0) circle (3 mm);
    \draw[thin] (2 cm,2 cm) circle (3 mm);
    \draw[thin,fill=black] (2 cm, 3mm) -- +(0, 1.4 cm);
  \end{tikzpicture}$ are isomorphic but not conjugate in $\Phi$. This is because their respective anisotropic kernels are of type $\smash{\big((A_1)^4\big)'}$ and $\smash{\big((A_1)^4\big)''}$ in $\Phi$. If instead $\Phi=D_9$ and $\Delta=D_8$ then the aforementioned abstract indices are conjugate in $\Phi$ since $\Stab_W(\Delta)\cong \Z_2$. This consideration occurs only five times in our tables, twice in Table \ref{E_7} when $\Delta=D_6A_1$ and three times in Table \ref{E_8} when $\Delta=D_8$.

\vspace{2mm}\noindent For each abstract index $\mathcal{H}$ in column $4$, we find the $W$-conjugacy class of $\Psi_0$ in $\Phi$ and put it in the fifth column. This is easy to do.

\vspace{2mm}\noindent By Lemma \ref{conjlemma}, the $W$-conjugacy class of $\Phi_a$ in $\Phi$ is independent of the choice of representative of the conjugacy class of $\mathcal{H}$ in $\Phi$. As such, for each abstract index $\mathcal{H}$ in column $4$, we compute the $W$-conjugacy class of $\Phi_a$ in $\Phi$ and put it in column $6$. This requires some work, and we perform these computations in Appendices $A$, $B$, $C$, $D$ and $E$ for the cases where $\Phi = G_2$, $F_4$, $E_6$, $E_7$ and $E_8$ respectively.

\vspace{2mm}\noindent Next we inspect all (isomorphism classes of) abstract indices $\mathcal{G}$ of type $\Phi$ (which are listed in Table II of \cite{T}). For each abstract index $\mathcal{H}$ in column $4$, we use conditions $(A.2)$, $(A.3)$ and $(A.4)$ of Definition \ref{mydefinition!} to exclude certain abstract indices of type $\Phi$.

\vspace{2mm}\noindent Using $(A.2)$, if $\Pi \leq \Stab_W(\Delta)$ then we exclude those abstract indices $\mathcal{G}$ of outer type. Otherwise, we exclude those abstract indices $\mathcal{G}$ of inner type (note that this condition is vacuous unless $\Phi=E_6$). Using $(A.3)$, we exclude those abstract indices $\mathcal{G}$ that do not contain a subindex of type $\Phi_a$. 

\vspace{2mm}\noindent Next we compute the $W$-conjugacy class of ${}_{in}\Phi_a$ and $\Phi \cap \overline{E_a}$ in $\Phi$ (once again, by Lemma \ref{conjlemma}, these are independent of the choice of representative of the conjugacy class of $\mathcal{H}$ in $\Phi$). Using $(A.4)$, we exclude those abstract indices $\mathcal{G}$ where $\Phi_0$ does not contain a subsystem of type ${}_{in}\Phi_a$. Again using $(A.4)$, we exclude $\mathcal{G}$ unless $\Phi_a$ contains a subsystem of type ${}_{in}\Phi_a \big(\Phi \cap \overline{E_a}\big)$ such that the ${}_{in}\Phi_a$ component is contained in $\Phi_0$. Note that these two conditions are vacuous if ${}_{in}\Phi_a$ is trivial. The second condition is used only once in our tables, for the case where $\mathcal{H}=\smash{\begin{tikzpicture}[baseline=-0.6ex]\dynkin[ply=2,fold radius = 2mm]{A}{ooo*ooo}\end{tikzpicture}}$ and $\Phi=E_7$. We list the resulting abstract indices $\mathcal{G}$ in column $7$ (if there are none we leave it blank).

\vspace{2mm}\noindent At this point (given the computations in the Appendix) we have justified the entries in the first seven columns of Tables \ref{G_2}, \ref{F_4}, \ref{E_6}, \ref{E_7} and \ref{E_8}. For each pair of indices $(\mathcal{H},\mathcal{G})$ listed in columns $4$ and $7$, the embedding $\theta$ of $\mathcal{H}$ in $\mathcal{G}$ has been uniquely determined (up to $\mathcal{G}$-conjugacy) by the conjugacy class information present in columns $1-4$. This procedure, applied and presented in the Appendix, completes the proof of Theorem \ref{classificationthm!}.
\end{proof}

\section{Existence results}\label{partialconverse}

\noindent In this section we prove Theorem \ref{existence!} and Corollaries \ref{maincor1}, \ref{maincor2} and \ref{maincor3}. We also reintroduce Conjecture \ref{myconjecture!}.

\vspace{2mm}\noindent Given a field $k$ with characteristic $p$, recall that a $p$-embedding of abstract indices $(\mathcal{G},\mathcal{H},\theta)$ is \textit{$k$-admissible} if there exists a pair of connected reductive algebraic $k$-groups $H \subset G$ such that the embedding of indices of $H \subset G$ is isomorphic to $(\mathcal{G},\mathcal{H},\theta)$. In this section we consider the \textit{existence problem}, that is, the problem of determining whether or not $(\mathcal{G},\mathcal{H},\theta)$ is $k$-admissible.

\vspace{2mm}\noindent We begin by restating Conjecture \ref{myconjecture!}, which claims that Definition \ref{mydefinition!} is as restrictive as is possible whilst ensuring that Theorem \ref{maintheorem!}$(i)$ holds. Let $\mathbb{P}$ denote the set of prime numbers.

\begin{myconjecture!} \textit{Let $p \in \mathbb{P} \cup \{0\}$ and let $(\mathcal{G},\mathcal{H},\theta)$ be a $p$-embedding of abstract indices. Then $(\mathcal{G},\mathcal{H},\theta)$ is $k$-admissible for some field $k$ with characteristic $p$.}
\end{myconjecture!}

\noindent Henceforth, given a $p$-embedding of abstract indices $(\mathcal{G},\mathcal{H},\theta)$, we use the notation of $\S \ref{embeddingsofabstractindices}$. That is, $\mathcal{G}=(F,\Lambda,\Lambda_0, \Xi)$, $\mathcal{H}=(E,\Delta,\Delta_0,\Pi)$, $\Phi:=\langle \Lambda \rangle$, $\Psi:=\langle \theta(\Delta) \rangle$, $\Psi_0:=\langle \theta(\Delta_0) \rangle$, $E_s$ is the fixed point subspace of $\smash{(E_{\Delta_0})^{\perp}}$ under $\Pi$, $E_a:=(E_s)^{\perp}$, $\Lambda_a := \Lambda \cap \theta(E_a)$, $\Phi_a := \Phi \cap \theta(E_a)$, ${}_{in}\Lambda_a$ is the union of irreducible components of $\Lambda_a$ that is maximal with respect to ${}_{in}\Lambda_a \subseteq \Psi_0$, $\mathcal{H}_m:=(E,\Delta_0, \Delta_0, \Pi)$ and $\mathcal{L}:=(F,\Lambda_a,\Lambda_0, \Xi)$.

\vspace{2mm}\noindent We now state and prove Theorem \ref{existence!}. Note that, in the statement of Theorem \ref{existence!}, if $G$ exists then so does $L$ (since $\mathcal{L}$ is a subindex of $\mathcal{G}$, see Proposition $21.12$ of \cite{B}).

\begin{existence!} \textit{Let $p \in \mathbb{P} \cup \{0\}$ and let $(\mathcal{G},\mathcal{H},\theta)$ be a $p$-embedding of abstract indices. Let $k$ be a field of characteristic $p$. Let $G$ be a connected reductive algebraic $k$-group with index isomorphic to $\mathcal{G}$. Let $L$ be a Levi $k$-subgroup of $G$ with index isomorphic to $\mathcal{L}$. There exists a $k$-subgroup of $G$ with embedding of indices isomorphic to $(\mathcal{G},\mathcal{H},\theta)$ if and only if there exists a $k$-subgroup of $L$ with embedding of indices isomorphic to $(\mathcal{L},\mathcal{H}_m,\theta)$.}

\begin{proof} We use the standard setup. Let $S$ be a maximal $k$-split torus of $L$ and let $T$ be a maximal $k$-torus of $L$ that contains $S$. Note that $L=C_G\big({Z(L)}^{\circ}_s\big)$ by Proposition $20.6$ of \cite{B}, and so $S$ is also a maximal $k$-split torus of $G$. Let $<$ be a $\Gamma$-order on $X(T)$. 
Let $\mathcal{I}(G)$ be the index of $G$ with respect to $T$ and $<$. By assumption, $\mathcal{I}(G)$ is isomorphic to $\mathcal{G}$. For ease of notation, we identify $\mathcal{I}(G)=\mathcal{G}$. Similarly, let $\mathcal{L}$ be the index of $L$ with respect to $T$ and $<$.

\vspace{2mm}\noindent Let $A$ be a connected reductive $k$-subgroup of maximal rank in $L$. Let $S_A$ be a maximal $k$-split torus of $A$ and let $T_A$ be a maximal $k$-torus of $A$ that contains $S_A$. Without loss of generality, as discussed in the proof of Theorem \ref{maintheorem!}, we can assume that $S_A \subseteq S$ and that $<$ is compatible with the restriction map $X(T) \to X(S_A)$. Let $g \in C_L(S_A)(K)$ such that $(T_A)^g=T$. By construction, $<_{g^{-1}}$ is a $\Gamma$-order on $X(T_A)$. Let $\mathcal{I}(A)$ be the index of $A$ with respect to $T_A$ and $<_{g^{-1}}$. Let $\theta_g:X(T_A)_{\R} \to X(T)_{\R}$ be defined by $\chi \mapsto \chi \cdot g^{-1}$ for $\chi \in X(T_A)$.

\vspace{2mm}\noindent Assume that the embedding of indices of $A \subset L$ is isomorphic to $(\mathcal{L},\mathcal{H}_m,\theta)$. For ease of notation, we identify $\mathcal{I}(A)=\mathcal{H}_m$ and $\theta_g=\theta$. Observe that $A'$ is $k$-anisotropic by Proposition \ref{abstractTits}, since $(\mathcal{H}_m)'$ is anisotropic. That is, $S_A={Z(A)}^{\circ}_s$. 

\vspace{2mm}\noindent Let $H:=\big\langle T_A, U_{\alpha} \hspace{0.5mm}\big|\hspace{0.5mm} \alpha \in \pm \Delta \big\rangle$, where each $U_{\alpha}$ denotes a $T_A$ root group of $G$. Recall from $(A.1)$ of Definition \ref{mydefinition!} that $\Psi$ is a $p$-closed subsystem of $\Phi$. Then, by Lemma \ref{propertiesinduced}, $H^g$ (and hence $H$) is a connected reductive $K$-subgroup of $G$. Moreover, recall from $(A.1)$ of Definition \ref{mydefinition!} that $\theta(\Delta) \subset \Phi^+$, and so $\Delta$ is the unique system of simple roots for $H$ with respect to $T_A$ that is compatible with $<_{g^{-1}}$.

\vspace{2mm}\noindent Let $\Gamma$ be the absolute Galois group of $k$. Consider the induced action $\iota_A:\Gamma \to \GL\!\big(X(T_A)_{\R}\big)$ (as defined in $\S \ref{*action}$). For every $\sigma \in \Gamma$ there exists a unique $w_{\sigma} \in W_{\Delta_0}$ such that $w_{\sigma}\iota_A(\sigma)$ stabilises $\Delta_0$. By assumption, $\Pi=\{w_{\sigma} \iota(\sigma)\hspace{0.5mm}|\hspace{0.5mm} \sigma \in \Gamma\}$ and $\Pi$ stabilises $\Delta$. Recall from $(A.1)$ of Definition \ref{mydefinition!} that $\Phi$ is $\Pi^{\theta}$-stable. Then $\langle \Delta \rangle$ and $\theta^{-1}(\Phi)$ are both $\iota_A(\Gamma)$-stable since they contain $\Delta_0$. Hence $H$ is defined over $k$ by Corollary \ref{downtown} and Remark \ref{niceactionremark}.

\vspace{2mm}\noindent By Lemmas \ref{niceactioncor} and \ref{Titsactionproperties}, $E_a$ is the subspace of $E$ that is spanned by all characters that vanish on $S_A$. Then, by Lemmas \ref{mybabyfixed} and \ref{propertiesofgammasystem}$(i)$ (and Remark \ref{niceactionremark}), $\Delta \cap E_a$ is a system of simple roots for $C_H(S_A)$ with respect to $T_A$. But $\Delta \cap E_a= \Delta_0$ by Lemma \ref{indexlemma}, 
and hence $C_H(S_A)=A$. Combining this with $S_A={Z(A)}^{\circ}_s$, we deduce that $S_A$ is a maximal $k$-split torus of $H$. 

\vspace{2mm}\noindent In summary, we have shown that $H$ is a connected reductive $k$-subgroup of $G$, $T_A$ is a maximal $k$-torus of $H$ containing a maximal $k$-split torus of $H$, $<_{g^{-1}}$ is a $\Gamma$-order on $X(T_A)$ and, finally, that $\mathcal{H}$ is the index of $H$ with respect to $T_A$ and $<_{g^{-1}}$. Certainly $\theta$ is a suitable orientation map, since $g$ was chosen to centralise $S_A$.

\vspace{2mm}\noindent We have already proved the converse in the process of proving Theorem \ref{maintheorem!}$(i)$, see Remark \ref{indexofL}.
\end{proof}
\end{existence!}

\noindent Henceforth, given a $p$-embedding of abstract indices $(\mathcal{G},\mathcal{H},\theta)$, we follow Remark \ref{abusenotation} and identify $\mathcal{H}$ with its image under $\theta$. In particular, we say that $E=F$.

\begin{maincor1} \textit{Let $p \in \mathbb{P} \cup \{0\}$. Let $(\mathcal{G},\mathcal{H},\theta)$ be an independent $p$-embedding of abstract indices. Let $k$ be a field of characteristic $p$ such that $\mathcal{G}$ is $k$-admissible. Then $(\mathcal{G},\mathcal{H},\theta)$ is $k$-admissible.}
\begin{proof} We have ${}_{in}\Lambda_a=\Lambda_a$ by definition of independence. Then $\Delta_0=\Lambda_a$ by Lemma \ref{indexlemma} and $(A.1)$ of Definition \ref{mydefinition!}. Moreover, $\Lambda_a=\Lambda_0$ by $(A.3)$ and $(A.4)$ of Definition \ref{mydefinition!}. The only element of $W_{\Delta_0}$ that stabilises $\Delta_0$ is the trivial element. Hence $\Pi=\Xi$ by $(A.2)$ of Definition \ref{mydefinition!}, since $\Xi$ and $\Pi$ both stabilise $\Delta_0$. In summary, $\mathcal{L}=\mathcal{H}_m$.

\vspace{2mm}\noindent Since $\mathcal{G}$ is $k$-admissible, there exists a connected reductive $k$-group $G$ with index isomorphic to $\mathcal{G}$. Let $L$ be a Levi $k$-subgroup of $G$ with index isomorphic to $\mathcal{L}$ (certainly this exists). Since $\mathcal{L}=\mathcal{H}_m$, the embedding of indices of $L$ in itself is isomorphic to $(\mathcal{L},\mathcal{H}_m,\theta)$. Hence, by Theorem \ref{existence!}, $G$ contains a $k$-subgroup $H$ such that the embedding of indices of $H \subset G$ is isomorphic to $(\mathcal{G},\mathcal{H},\theta)$.
\end{proof}
\end{maincor1}

\begin{maincor2} \textit{Let $p \in \mathbb{P} \cup \{0\}$. Let $(\mathcal{G},\mathcal{H},\theta)$ be a quasisplit $p$-embedding of abstract indices. Let $k$ be a field of characteristic $p$ such that $\mathcal{G}$ is $k$-admissible. Then $(\mathcal{G},\mathcal{H},\theta)$ is $k$-admissible if and only if $\Pi$ is a quotient of the absolute Galois group $\Gamma$ of $k$.}
\begin{proof} Assume that $\Pi$ is a quotient of $\Gamma$. Recall that we identify $\mathcal{H}$ with its image under $\theta$. Let $\mathcal{G}=(F,\Lambda,\Lambda_0, \Xi)$ and $\mathcal{H}=(E,\Delta,\Delta_0,\Pi)$. Note that $\Lambda_0$ and $\Delta_0$ are both empty by definition of quasisplit. 

\vspace{2mm}\noindent Let $\Iso(E)$ be the isometry group of $E$. Note that $\Stab_{\Iso(E)}(\Phi_a) = W_{\Lambda_a} \rtimes \Stab_{\Iso(E)}(\Lambda_a)$, where $\Stab_{\Iso(E)}(\Lambda_a)$ acts on $W_{\Lambda_a}$ in the obvious way (see for instance $(A.8)$ of \cite{H}). Let $\rho:\Stab_{\Iso(E)}(\Phi_a) \to \Stab_{\Iso(E)}(\Lambda_a)$ be the associated projection. It follows from $(A.1)$ of Definition \ref{mydefinition!} that $\Pi \subseteq \Stab_{\Iso(E)}(\Phi_a)$. Then $\rho(\Pi)=\Xi$ by $(A.2)$ of Definition \ref{mydefinition!}. 

\vspace{2mm}\noindent Let $V:=\Span(\Lambda_a)$ and let $V^{\perp}$ be the orthogonal complement of $V$ in $E$. Consider the natural projections $\pi:\Stab_{\Iso(E)}(V) \to\Stab_{\Iso(E)}(V)/ \Fix_{\Iso(E)}(V)$ and $\pi^{\perp}:\Stab_{\Iso(E)}(V^{\perp}) \to \Stab_{\Iso(E)}(V^{\perp})/ \Fix_{\Iso(E)}(V^{\perp})$.

\vspace{2mm}\noindent Since $\mathcal{G}$ is $k$-admissible, there exists a connected reductive $k$-group $G$ with index isomorphic to $\mathcal{G}$. Let $L$ be a $k$-subgroup of $G$ with index isomorphic to $\mathcal{L}$ (this must exist). Recall from Proposition \ref{abstractTits} that $L$ is $k$-quasisplit and the index of $L'$ is isomorphic to $(V,\Lambda_a,\Lambda_0, \pi(\Xi))=:\mathcal{L}'$. By assumption, $\pi(\Pi)$ is a quotient of $\Gamma$. Observe that $\rho \pi(\Pi)=\pi(\Xi)$ since $\rho(\Pi)=\Xi$ and $\pi$ commutes with $\rho$. Then, by Theorem $1.1$ of \cite{G}, there exists a maximal $k$-torus $T_0$ of $L'$ with index isomorphic to $(V,\varnothing,\varnothing, \pi(\Pi))$.

\vspace{2mm}\noindent Let $T:=T_0{Z(L)}^{\circ}$. Recall from Proposition \ref{abstractTits} that the index of ${Z(L)}^{\circ}$ is isomorphic to $(V^{\perp},\varnothing,\varnothing, \pi^{\perp}(\Xi))=:\overline{\mathcal{L}}$. Observe that $\pi^{\perp}(\Pi)=\pi^{\perp}(\Xi)$ since $\rho(\Pi)=\Xi$. Since $\Pi$ is uniquely defined by $\pi(\Pi)$ and $\pi^{\perp}(\Pi)$, we conclude that the embedding of indices of $T \subset L$ is isomorphic to $(\mathcal{L},\mathcal{H}_m,\theta)$. Hence, by Theorem \ref{existence!}, $G$ contains a $k$-subgroup $H$ such that the embedding of indices of $H \subset G$ is isomorphic to $(\mathcal{G},\mathcal{H},\theta)$.

\vspace{2mm}\noindent The converse follows immediately from the definition of the Tits action.
\end{proof}
\end{maincor2}

\noindent It remains to prove Corollary \ref{maincor3}. For $k$ a field, let $\cd(k)$ denote the cohomological dimension of the absolute Galois group of $k$. Recall that an \textit{embedding of abstract indices} is a $p$-embedding of abstract indices for some $p \in \mathbb{P} \cup \{0\}$.

\begin{maincor3} \textit{Let $(\mathcal{G},\mathcal{H},\theta)$ be an embedding of abstract indices that is listed in Table \ref{G_2}, \ref{F_4}, \ref{E_6}, \ref{E_7} or \ref{E_8}. If there is a $\checkmark$ in column $8$ (resp. $9$, $10$) of the corresponding row then $(\mathcal{G},\mathcal{H},\theta)$ is $k$-admissible for some field $k$ with $\cd(k)=1$ (resp. $k=\R$, $k$ is $\mathfrak{p}$-adic). If there is a $\xmark$ then $(\mathcal{G},\mathcal{H},\theta)$ is not $k$-admissible for such a field.}

\vspace{2mm}\noindent \textit{[In the cases where we list more than one value for $\Pi$ in a single row of one of the tables, we put a $\checkmark$ in the corresponding column if at least one of these values for $\Pi$ corresponds to a $k$-admissible $(\mathcal{G},\mathcal{H},\theta)$.]}
\begin{proof} We describe an algorithm which we use to determine whether or not $(\mathcal{G},\mathcal{H},\theta)$ is $k$-admissible for some field $k$ with $\cd(k)=1$ (resp. $k=\R$, $k$ is $\mathfrak{p}$-adic).

\vspace{2mm}\noindent Given $p \in \mathbb{P} \cup \{0\}$, we say that $(\mathcal{G},\mathcal{H},\theta)$ is \textit{$p$-valid} if it is a $p$-embedding of abstract indices.

\vspace{2mm}\noindent We first consider the cases where $(\mathcal{G},\mathcal{H},\theta)$ is not $0$-valid. We put a $\xmark$ in columns $9$ and $10$ (of the corresponding row) since $\Char(k)=0$ if $k=\R$ or $k$ is $\mathfrak{p}$-adic. If $\cd(k) =1$ then $\Lambda_0$ and $\Delta_0$ are both empty by $\S 3.3.1$ of \cite{T}. So we put a $\xmark$ in column $8$ unless $\Lambda_0$ and $\Delta_0$ are both empty. There are five remaining cases, only one of which has a non-cyclic $\Pi$. We ignore this case. 

\vspace{2mm}\noindent For the four remaining cases, the first two are $3$-valid with $\mathcal{G}=\dynkin{G}{I}$ and the other two are $2$-valid with $\mathcal{G}=\dynkin{F}{I}$. For the first two cases (resp. other two cases), let $k$ be the finite field of order $3$ (resp. $2$). The absolute Galois group $\Gamma$ of $k$ is isomorphic to $\smash{\hat{\Z}}$, so $\cd(k)=1$ and any cyclic group is a quotient of $\Gamma$. The abstract index $\mathcal{G}$ is certainly $k$-admissable as it is split. Hence, by Corollary \ref{maincor2}, we put a $\checkmark$ in column $8$.

\vspace{2mm}\noindent Henceforth we consider the cases where $(\mathcal{G},\mathcal{H},\theta)$ is $0$-valid. 

\vspace{2mm}\noindent \underline{\textit{Column $8$}}: $\cd(k) =1$.

\vspace{2mm}\noindent If $\cd(k) =1$ then again $\Lambda_0$ and $\Delta_0$ are both empty by $\S 3.3.1$ of \cite{T}. Consider the field $k=\C(t)$ of rational functions over $\C$ and its absolute Galois group $\Gamma$. The characteristic of $k$ is $0$. By Corollary $3.4.4$ and Theorem $3.4.8$ of \cite{Sz}, $\cd(k)=1$ and any finite group is a quotient of $\Gamma$. The index $\smash{\begin{tikzpicture}[baseline=-0.6ex]\dynkin[ply=2,fold radius=2mm]{E}{oooooo}\end{tikzpicture}}$ is $k$-admissible (since $k$ admits a quadratic extension) and all split indices are $k$-admissible. 
So, by Corollary \ref{maincor2}, we put a $\checkmark$ in column $8$ if $\Lambda_0$ and $\Delta_0$ are both empty and a $\xmark$ otherwise. 

\vspace{2mm}\noindent We cross-check the entries in column $8$ using Table $5.1$ of \cite{LSS}.

\vspace{2mm}\noindent \underline{\textit{Column $9$}}: $k=\R$.

\vspace{2mm}\noindent We use Table II of \cite{T} to check if $\mathcal{G}$ is $k$-admissible. If not, we put a $\xmark$ in column $9$. If so, then we proceed as follows.

\vspace{2mm}\noindent If $(\mathcal{G},\mathcal{H},\theta)$ is independent then we put a $\checkmark$ in column $9$ by Corollary \ref{maincor1}. If $\Lambda_0$ and $\Delta_0$ are both empty and $(\mathcal{G},\mathcal{H},\theta)$ is not independent then, by Corollary \ref{maincor2}, we put a $\checkmark$ (resp. $\xmark$) in column $9$ if $\Pi=1$ or $\Z_2$ (resp. $|\Pi|>2$). 

\vspace{2mm}\noindent If $\Lambda_0$ or $\Delta_0$ is non-empty and $(\mathcal{G},\mathcal{H},\theta)$ is not independent then we take a pair of complex groups $H_{\C} \subset G_{\C}$ that share a maximal torus $T_{\C}$ and that have root systems of type $\Delta$ and $\Lambda$ respectively. We classify conjugacy classes of holomorphic involutions of $G_{\C}$ that stabilise $H_{\C}$ and $T_{\C}$. By $\S 3$ of \cite{K}, this is equivalent to classifying pairs of $\R$-groups $H \subset G$ with respective complexifications $H_{\C} \subset G_{\C}$. 
We do not include these computations in this paper since the case $k = \R$ has already been done by Komrakov \cite{K} and de Graff-Marrani \cite{DM}.

\vspace{2mm}\noindent We cross-check the entries in column $9$ using Table $2$ of \cite{K}. 

\vspace{2mm}\noindent \underline{\textit{Column $10$}}: $k$ is $\mathfrak{p}$-adic.

\vspace{2mm}\noindent We first observe, using II.$7$ of \cite{GMS}, that the absolute Galois group of any $\mathfrak{p}$-adic field has $\smash{\hat{\Z}}$ as a quotient. In addition, observe that the index $\smash{\begin{tikzpicture}[baseline=-0.6ex]\dynkin[ply=2,fold radius=2mm]{E}{oooooo}\end{tikzpicture}}$ is $k$-admissible for $k$ any $\mathfrak{p}$-adic field (since $k$ admits a quadratic extension).

\vspace{2mm}\noindent We use Table II of \cite{T} to check if $\mathcal{G}$ is $k$-admissible for some $\mathfrak{p}$-adic field $k$. If not, we put a $\xmark$ in column $10$. If so, then we proceed as follows.

\vspace{2mm}\noindent If $(\mathcal{G},\mathcal{H},\theta)$ is independent then we put a $\checkmark$ in column $10$ by Corollary \ref{maincor1}. If $\Lambda_0$ and $\Delta_0$ are both empty, $\Pi$ is cyclic and $(\mathcal{G},\mathcal{H},\theta)$ is not independent then, by Corollary \ref{maincor2}, we put a $\checkmark$ in column $10$ (since all finite cyclic groups are quotients of $\smash{\hat{\Z}}$).

\vspace{2mm}\noindent For all remaining cases, we are unsure of whether or not $(\mathcal{G},\mathcal{H},\theta)$ is $k$-admissible for some $\mathfrak{p}$-adic field $k$. So we put a $?$ in column $10$.

\vspace{2mm}\noindent This procedure, presented in the Appendix, completes the proof of Corollary \ref{maincor3}.
\end{proof}
\end{maincor3}

\section*{Appendix}\label{appendix}

\noindent In the appendix we complete the statement and proof of Theorem \ref{classificationthm!}. We perform the calculations needed to construct Tables \ref{G_2}, \ref{F_4}, \ref{E_6}, \ref{E_7} and \ref{E_8} and then we present the tables. In particular, for each abstract index $\mathcal{H}$ listed in the fourth column of our tables, we use the results of $\S \ref{embeddingsofabstractindices}$ to compute the $W$-conjugacy class of $\Phi_a$ in $\Phi$ (the entry in the sixth column of our tables). 

\vspace{2mm}\noindent We consider the case where $\Phi$ is $G_2$, $F_4$, $E_6$, $E_7$ and $E_8$ in Appendices $A$, $B$, $C$, $D$ and $E$ respectively. 

\vspace{2mm}\noindent Recall the setup from $\S \ref{restrictions}$. Let $p$ either be a prime number or $0$. Let $(\mathcal{G}, \mathcal{H}, \theta)$ be a $p$-embedding of abstract indices, say $\mathcal{G}=(F,\Lambda,\Lambda_0, \Xi)$ and $\mathcal{H}=(E,\Delta,\Delta_0,\Pi)$, where $\Lambda$ is irreducible of exceptional type and $r(\Lambda)=\dim(F)$. We follow Remark \ref{abusenotation} and identify $E=F$. We denote $\Phi:=\langle \Lambda \rangle$, $W:=W_{\Lambda}$, $\Psi:=\langle \Delta \rangle$, $\Psi_0:=\langle \Delta_0 \rangle$, $E_s$ is the fixed point subspace of $\smash{(E_{\Delta_0})^{\perp}}$ under $\Pi$, $r_s(\mathcal{H}):=\dim (E_s)$, $\Phi_s := \Phi \cap E_s$, $E_a:=(E_s)^{\perp}$ and $\Phi_a := \Phi \cap E_a$. Let $\Iso(\Phi)$ be the isometry group of $E$ that stabilises $\Phi$. 

\vspace{2mm}\noindent We introduce some additional notation. Let $\alpha_0$ be the highest root of $\Psi$. If $\Psi$ is reducible then let $\Psi^i$ be the $i$'th irreducible component of $\Psi$ with respect to the ordering used in $\mathcal{H}$. For example, if $\mathcal{H}=\smash{\begin{tikzpicture}[baseline=-0.2ex]\dynkin[ply=2,fold radius = 2mm]{E}{*o****}\end{tikzpicture} \hspace{-0.3mm}\times\hspace{-0.5mm} \begin{tikzpicture}[baseline=-0.2ex]\dynkin[ply=2,fold radius = 2mm]{A}{oo}\end{tikzpicture}}$ then $\Psi^1$ (resp. $\Psi^2$) is the irreducible component of $\Psi$ of type $E_6$ (resp. $A_2$). Let $\Delta^i := \Delta \cap \Psi^i$, $\Delta^i_0 := \Delta_0 \cap \Psi^i$ and let $\alpha_0^i$ be the highest root of $\Psi^i$. For $\Sigma$ a subsystem of $\Phi$, let $\Aut_W(W_{\Sigma})$ be the group of graph automorphisms of the Dynkin diagram of $\Sigma$ that are induced by elements of $W$.

\vspace{2mm}\noindent If there exists more than one $W$-conjugacy class of a given isomorphism class of subsystems of $\Phi$, we distinguish them using Carter's notation in $\S 7,8$ and Tables $7$-$11$ of \cite{C}. For example, if $\Phi=E_8$, we denote a parabolic (resp. non-parabolic) subsystem of $\Phi$ of type $(A_3)^2$ by $\smash{{\big((A_3)^2\big)}''}$ (resp. $\smash{{\big((A_3)^2\big)}'}$). If $\Phi$ is not simply-laced then we use a $\widetilde{\textcolor{white}{m}}$ to denote a subsystem of $\Phi$ that consists only of short roots.

\subsection*{Appendix A. $\Phi=G_2 = \dynkin[scale=2,labels={\alpha_1,\alpha_2}]{G}{2}$}\label{G2}

\noindent Let $\Phi=G_2$. We use the following \textit{standard} parametrisation in $\R^2$ of $\Phi$. Let $\alpha_1=2e_2-e_1-e_3$ and $\alpha_2=e_1-e_2$. The highest root of $\Phi$ is $e_1+e_2-2e_3$. The Weyl group of $\Phi$ is $W \cong \mathcal{D}i_{12}$ and the longest element of $W$ is $-1$.

\vspace{2mm}\noindent By Table II of \cite{T}, the only isotropic index of type $\Phi$ is $\dynkin{G}{I}$. This index is of inner type since $\Iso(\Phi)=W$.

\vspace{2mm}\noindent Recall from Table \ref{exceptlist} that $\Delta$ is one of $A_2$, $\smash{\widetilde{A_2}}$ (if $p=3$) or $\smash{A_1\widetilde{A_1}}$. In the fourth column of Table \ref{G_2} we list all indices $\mathcal{H}$ of type $\Delta$ that are maximal in $\Phi$. For each of those indices $\mathcal{H}$, we compute the $W$-conjugacy class of $\Phi_a$ in $\Phi$.

\begin{table}[!htb]\begin{center}\caption{Isotropic maximal embeddings of abstract indices of type $G_2$}\label{G_2}\begin{tabular}{| c | c | c | c | c | c | c | c | c | c | c |}    \hline
\multicolumn{1}{|c|}{} & \multicolumn{1}{c|}{} & \multicolumn{1}{c|}{} & \multicolumn{1}{c|}{} & \multicolumn{1}{c|}{} & \multicolumn{1}{c|}{} &  \multicolumn{1}{c|}{} & \multicolumn{3}{c|}{\small{Special fields}} \\

    $\Delta$ & \small{$\!\Stab_W(\Delta)\!$} & $\Pi$ & $\mathcal{H}$ & $\Psi_0$ & $\Phi_a$ & $\mathcal{G}$ & \small{$\hspace{-0.8mm}\!\cd 1 \!\hspace{-1mm}$} & $\R$ & $\!\!Q_{\mathfrak{p}}\!\!$  \\ \hline \hline 

\multirow{2.3}{*}{$A_2$} & \multirow{2.3}{*}{$\Z_2$} & $1$ & $\dynkin{A}{oo}$ & $\varnothing$ & $\varnothing$ & $\dynkin{G}{I}$ & $\checkmark$ & $\checkmark$ & $~\checkmark$ \topstrut \\ \cline{3-11} 
& & $\Z_2$ & $\begin{tikzpicture}[baseline=-0.3ex]\dynkin[ply=2,rotate=90,fold radius = 2mm]{A}{oo}\end{tikzpicture}$ & $\widetilde{A_1}$ & $\widetilde{A_1}$ & $\dynkin{G}{I}$ & $\checkmark$ & $\checkmark$ & $~\checkmark$ \topstrut \\ \hline 

$\multirow{2.3}{*}{\makecell{$\widetilde{A_2}$ \\ \small{$\!(p \! = \! 3)\!$}}}$ & \multirow{2.3}{*}{$\Z_2$} & $1$ & $\widetilde{\smash{\dynkin{A}{oo}}}$ & $\varnothing$ & $\varnothing$ & $\dynkin{G}{I}$ & $\checkmark$ & $\xmark$ & $~\xmark$  \topstrut \\ \cline{3-11} \noalign{\vskip 0.1mm}  
 & & $\Z_2$ & $\simcal{\begin{tikzpicture}[baseline=-0.3ex]\dynkin[ply=2,rotate=90,fold radius = 2mm]{A}{oo}\end{tikzpicture}}$ & $A_1$ & $A_1$ & $\dynkin{G}{I}$ & $\checkmark$ & $\xmark$ & $~\xmark$  \topstrut \\ \hline 

\multirow{3.45}{*}{$\!A_1\widetilde{A_1}\!$} & \multirow{3.45}{*}{$1$} & \multirow{3.45}{*}{$1$} & $\dynkin{A}{o} \hspace{-0.3mm}\times\hspace{-0.3mm} \widetilde{\smash{\dynkin{A}{o}}}$ & $\varnothing$ & $\varnothing$ & $\dynkin{G}{I}$ & $\checkmark$ & $\checkmark$ & $~\checkmark$ \topstrut \\ \cdashline{4-11} 
& & & $\dynkin{A}{o} \hspace{-0.3mm}\times\hspace{-0.3mm} \widetilde{\smash{\dynkin{A}{*}}}$ & $\widetilde{A_1}$ & $\widetilde{A_1}$ &&&& \topstrut \\ \cdashline{4-11} 
& & & $\dynkin{A}{*} \hspace{-0.3mm}\times\hspace{-0.3mm} \widetilde{\smash{\dynkin{A}{o}}}$ & $A_1$ & $A_1$ &&&& \topstrut  \\ \hline \noalign{\vskip 0.1mm} 
  \end{tabular}\end{center}\end{table}

\vspace{2mm}\noindent If $\mathcal{H}=\begin{tikzpicture}[baseline=-0.3ex]\dynkin[ply=2,rotate=90,fold radius = 2mm]{A}{oo}\end{tikzpicture}$ then $\Phi_s=\big\langle  \alpha_0 \big\rangle \cong A_1$ and $r(\Phi_s)=r_s(\mathcal{H}) =1$. Hence $\Phi_a =(\Phi_s)^{\perp} \cong \smash{\widetilde{A_1}}$ since $\Phi_s$ spans $E_s$.

\vspace{2mm}\noindent If $\mathcal{H}=\simcal{\begin{tikzpicture}[baseline=-0.3ex]\dynkin[ply=2,rotate=90,fold radius = 2mm]{A}{oo}\end{tikzpicture}}$ then $\Phi_s=\big\langle  \alpha_0 \big\rangle \cong \smash{\widetilde{A_1}}$ and $r(\Phi_s)=r_s(\mathcal{H}) =1$. Hence $\Phi_a =(\Phi_s)^{\perp} \cong A_1$.

 \vspace{2mm}\noindent For all remaining indices $\mathcal{H}$ listed in the fourth column of Table \ref{F_4}, we observe that $\Psi_0$ is a parabolic subsystem of $\Phi$ and we use Lemma \ref{complemma} to check that $r(\Phi_a)=r(\Psi_0)$. So $\Phi_a=\Psi_0$ since $\Phi_a$ is a parabolic subsystem of $\Phi$ that contains $\Psi_0$. That is, $\mathcal{H}$ is independent in $\Phi$.

\subsection*{Appendix B. $\Phi=F_4 = \dynkin[scale=2,labels={\alpha_1,\alpha_...,\alpha_4}]{F}{4}$}\label{F4}

\noindent Let $\Phi=F_4$. We use the following \textit{standard} parametrisation in $\R^4$ of $\Phi$. Let $\alpha_1=e_1-e_2$, $\alpha_2=e_2-e_3$, $\alpha_3=e_3$ and $\alpha_4=\big(-\sum_{i=1}^4 e_i\big)/2$. The highest root of $\Phi$ is $e_1-e_4$. The Weyl group of $\Phi$ is $W \cong (\Z_2^3 \rtimes S_4) \rtimes S_3$ and the longest element of $W$ is $-1$.

\vspace{2mm}\noindent By Table II of \cite{T}, the isotropic indices of type $\Phi$ are $\dynkin{F}{I}$ and $\dynkin{F}{II}$. Both of these indices are of inner type since $\Iso(\Phi)=W$.

\begin{table}[!htb]\begin{center}\caption{Isotropic maximal embeddings of abstract indices of type $F_4$}\label{F_4}\hspace*{0pt}\begin{tabular}{| c | c | c | c | c | c | c | c | c | c | c |}    \hline
\multicolumn{1}{|c|}{} & \multicolumn{1}{c|}{} & \multicolumn{1}{c|}{} & \multicolumn{1}{c|}{} & \multicolumn{1}{c|}{} & \multicolumn{1}{c|}{} & \multicolumn{1}{c|}{} & \multicolumn{3}{c|}{\small{Special fields}} \\

    $\Delta$ & \small{$\!\Stab_W(\Delta)\!$} & $\Pi$ & $\mathcal{H}$ & $\Psi_0$ & $\Phi_a$ & $\mathcal{G}$ & \small{$\hspace{-3mm}\!\cd 1 \!\hspace{-3mm}$} & $\R$ & $\!\!Q_{\mathfrak{p}}\!\!$  \\ \hline \hline 

    \multirow{4.6}{*}{$B_4$} & \multirow{4.6}{*}{$1$} & \multirow{4.6}{*}{$1$} &  $\dynkin{B}{oooo}$ & $\varnothing$  & $\varnothing$ & $\dynkin{F}{I}$ & $\checkmark$ & $\checkmark$ & $~\checkmark$ \topstrut  \\ \cdashline{4-11} 
&&& $\dynkin{B}{ooo*}$ & $\widetilde{A_1}$ & $\widetilde{A_1}$ &&&& \topstrut \\ \cdashline{4-11} 
&&& $\dynkin{B}{oo**}$ & $B_2$ & $B_2$ &&&& \topstrut \\ \cdashline{4-11} 
&&& $\dynkin{B}{o***}$ & $B_3$ & $B_3$ & $\dynkin{F}{II}$ & $\xmark$ & $\checkmark$ & $~\xmark$ \topstrut \\ \hline 

\multirow{2.5}{*}{$D_4$} & \multirow{2.5}{*}{$S_3$} & \multirow{2.5}{*}{\makecell{$\Z_3 \textnormal{ or}\!$ \\ $S_3$}} & \begin{tikzpicture}[scale=0.8,baseline=-2.2ex] \dynkin[ply=3,rotate=270]{D}{oooo} \end{tikzpicture} & $\varnothing$ & $\widetilde{A_2}$  & $\dynkin{F}{I}$ & $\checkmark$ & $\xmark$ & $\checkmark$ \\ \cdashline{4-11} 
&&& \begin{tikzpicture}[scale=0.8,baseline=-2.2ex] \dynkin[ply=3,rotate=270,fold radius = 2mm]{D}{*o**} \end{tikzpicture} & $(A_1)^3$ & $C_3$  & $\dynkin{F}{I}$ & $\xmark$ & $\xmark$ & $\xmark$  \\ \hline 

\multirow{4.5}{*}{\makecell{$C_4$ \\ \small{$\!(p\!=\!2)\!$}}} & \multirow{4.5}{*}{$1$} & \multirow{4.5}{*}{$1$} & $\begin{tikzpicture}[baseline=-1.0ex]\dynkin[rotate=180]{C}{oooo}\end{tikzpicture}$ & $\varnothing$ & $\varnothing$ & $\dynkin{F}{I}$ & $\checkmark$ & $\xmark$ & $~\xmark$ \topstrut \\ \cdashline{4-11} 
&&& $\begin{tikzpicture}[baseline=-1.0ex]\dynkin[rotate=180]{C}{*o*o}\end{tikzpicture}$ & $(\widetilde{A_1})^2$ & $B_2$  & $\dynkin{F}{I}$ & $\xmark$ & $\xmark$ & $~\xmark$ \topstrut \\ \cdashline{4-11} 
&&& \multirow{2}{*}{$\begin{tikzpicture}[baseline=-1.0ex]\dynkin[rotate=180]{C}{*o**}\end{tikzpicture}$} & \multirow{1.9}{*}{$B_2\widetilde{A_1}$} & \multirow{2}{*}{$B_3$}  & $\dynkin{F}{I}$ & $\xmark$ & $\xmark$ & $~\xmark$ \topstrut \\  
&&&&& & $\dynkin{F}{II}$ & $\xmark$ & $\xmark$ & $\xmark$  \\ \hline

\multirow{3.3}{*}{\makecell{$\!\widetilde{D_4}$ \\ \small{$\!(p\!=\!2)\!$}}} & \multirow{3.3}{*}{$S_3$} & \multirow{3.3}{*}{\makecell{$\!\Z_3 \textnormal{ or}\!$ \\ $S_3$}} & $\begin{tikzpicture}[scale=0.8,baseline=-1.7ex] \dynkin[ply=3,rotate=270]{D}{oooo} \end{tikzpicture}\!\!\!\!\!\!\!\!\!\!\!\!\!\!\!\!\!\!\hspace{0.18mm}\shimcal{\begin{tikzpicture}[scale=0.8,baseline=-1.7ex] \dynkin[ply=3,rotate=270]{D}{oooo} \end{tikzpicture}}\!\!\!\!\hspace{-0.2mm}$ & \multirow{1.2}{*}{$\varnothing$} & \multirow{1.2}{*}{$A_2$}  & \multirow{1.2}{*}{$\dynkin{F}{I}$} & \multirow{1.2}{*}{$\checkmark$} & \multirow{1.2}{*}{$\xmark$} & \multirow{1.2}{*}{$\xmark$} \\ \cdashline{4-11} 
&&& \multirow{1.4}{*}{$\shimcal{\begin{tikzpicture}[scale=0.8,baseline=-1.7ex] \dynkin[ply=3,rotate=270,fold radius = 2mm]{D}{*o**} \end{tikzpicture}}$} & \multirow{1.9}{*}{$(\widetilde{A_1})^3$} & \multirow{2}{*}{$B_3$}  & $\dynkin{F}{I}$ & $\xmark$ & $\xmark$ & $~\xmark$  \topstrut \\ 
&&&&& & $\dynkin{F}{II}$ & $\xmark$ & $\xmark$ & $\xmark$  \\ \hline 

    \multirow{6.8}{*}{\makecell{$\!C_3A_1\!$ \\ \small{$\!(p\!\neq\!2)\!$}}} & \multirow{6.8}{*}{$1$} & \multirow{6.8}{*}{$1$} & $\begin{tikzpicture}[baseline=-1.0ex]\dynkin[rotate=180]{C}{ooo}\end{tikzpicture} \hspace{-0.3mm}\times\hspace{-0.3mm} \dynkin{A}{o}$ & $\varnothing$ & $\varnothing$ & $\dynkin{F}{I}$ & $\checkmark$ & $\checkmark$ & $~\checkmark$ \topstrut \\ \cdashline{4-11} 
&&& $\begin{tikzpicture}[baseline=-1.0ex]\dynkin[rotate=180]{C}{ooo}\end{tikzpicture} \hspace{-0.3mm}\times\hspace{-0.3mm} \dynkin{A}{*}$ & $A_1$ & $A_1$ &&&& \topstrut \\ \cdashline{4-11} 
&&& $\begin{tikzpicture}[baseline=-1.0ex]\dynkin[rotate=180]{C}{*o*}\end{tikzpicture} \hspace{-0.3mm}\times\hspace{-0.3mm} \dynkin{A}{o}$ & $A_1\widetilde{A_1}$ & $A_1\widetilde{A_1}$ &&&& \topstrut  \\ \cdashline{4-11} 
&&& $\begin{tikzpicture}[baseline=-1.0ex]\dynkin[rotate=180]{C}{***}\end{tikzpicture} \hspace{-0.3mm}\times\hspace{-0.3mm} \dynkin{A}{o}$ & $C_3$ & $C_3$ &&&& \topstrut  \\ \cdashline{4-11} 
&&& \multirow{2}{*}{$\begin{tikzpicture}[baseline=-1.0ex]\dynkin[rotate=180]{C}{*o*}\end{tikzpicture} \hspace{-0.3mm}\times\hspace{-0.3mm} \dynkin{A}{*}$} & \multirow{1.9}{*}{$(A_1)^2\widetilde{A_1}$} & \multirow{2}{*}{$B_3$} & $\dynkin{F}{I}$ & $\xmark$ & $\checkmark$ & $~?$  \topstrut  \\ 
&&&&& & $\dynkin{F}{II}$ & $\xmark$ & $\checkmark$ & $\xmark$   \\ \hline

    \multirow{7.75}{*}{$\!A_2\widetilde{A_2}\!$} & \multirow{7.75}{*}{$\Z_2$} & \multirow{3.45}{*}{$1$} & $\dynkin{A}{oo} \hspace{-0.3mm}\times\hspace{-0.3mm} \widetilde{\smash{\dynkin{A}{oo}}}$ & $\varnothing$ & $\varnothing$ & $\dynkin{F}{I}$ & $\checkmark$ & $\checkmark$ & $~\checkmark$ \topstrut \\ \cdashline{4-11} 
&&& $\dynkin{A}{oo} \hspace{-0.3mm}\times\hspace{-0.3mm} \widetilde{\smash{\dynkin{A}{**}}}$ & $\widetilde{A_2}$ & $\widetilde{A_2}$ &&&& \topstrut \\ \cdashline{4-11} 
&&& $\dynkin{A}{**} \hspace{-0.3mm}\times\hspace{-0.3mm} \widetilde{\smash{\dynkin{A}{oo}}}$ & $A_2$ & $A_2$ &&&& \topstrut \\ \cline{3-11} 

 & & \multirow{4.3}{*}{$\Z_2$} & $\begin{tikzpicture}[baseline=-0.3ex]\dynkin[ply=2,rotate=90,fold radius = 2mm]{A}{oo}\end{tikzpicture} \hspace{0.2mm}\times\hspace{-0.2mm} \! \simcal{\begin{tikzpicture}[baseline=-0.3ex]\dynkin[ply=2,rotate=90,fold radius = 2mm]{A}{oo}\end{tikzpicture}}\!\hspace{0.4mm}$ & $\varnothing$ & $\!A_1\widetilde{A_1}\!$ & $\dynkin{F}{I}$ & $\checkmark$ & $\checkmark$ & $~\checkmark$  \topstrut \\ \cdashline{4-11} 
&&& $\begin{tikzpicture}[baseline=-0.3ex]\dynkin[ply=2,rotate=90,fold radius = 2mm]{A}{oo}\end{tikzpicture} \hspace{0.2mm}\times\hspace{-0.2mm} \! \simcal{\begin{tikzpicture}[baseline=-0.3ex]\dynkin[ply=2,rotate=90,fold radius = 2mm]{A}{**}\end{tikzpicture}}\!\hspace{0.4mm}$ & $\widetilde{A_2}$ & $C_3$ & $\dynkin{F}{I}$ & $\xmark$ & $\checkmark$ & $~\xmark$ \topstrut \\ \cdashline{4-11} 
&&& \multirow{2.1}{*}{$\begin{tikzpicture}[baseline=-0.3ex]\dynkin[ply=2,rotate=90,fold radius = 2mm]{A}{**}\end{tikzpicture} \hspace{0.2mm}\times\hspace{-0.2mm} \! \simcal{\begin{tikzpicture}[baseline=-0.3ex]\dynkin[ply=2,rotate=90,fold radius = 2mm]{A}{oo}\end{tikzpicture}}\!\hspace{0.4mm}$} & \multirow{2.2}{*}{$A_2$} & \multirow{2.1}{*}{$B_3$} & $\dynkin{F}{I}$ & $\xmark$ & $\xmark$ & $~\xmark$  \topstrut \\ 
&&&&& & $\dynkin{F}{II}$ & $\xmark$ & $\checkmark$ & $\xmark$ \\ \hline

  \end{tabular}\hspace*{0pt}\end{center}\end{table}

\vspace{2mm}\noindent Recall from Table \ref{exceptlist} that $\Delta$ is one of $B_4$, $D_4$, $C_4$ (if $p=2$), $\smash{\widetilde{D_4}}$ (if $p=2$), $C_3A_1$ (if $p \neq 2$) or $\smash{A_2\widetilde{A_2}}$. In the fourth column of Table \ref{F_4} we list all indices $\mathcal{H}$ of type $\Delta$ that are maximal in $\Phi$. For each of those indices $\mathcal{H}$, we compute the $W$-conjugacy class of $\Phi_a$ in $\Phi$.

\vspace{2mm}\noindent If $\mathcal{H}=\begin{tikzpicture}[scale=0.8,baseline=-2.2ex] \dynkin[ply=3,rotate=270,fold radius = 2mm]{D}{*o**}\end{tikzpicture}$ then $\Phi_s=(\Psi_0)^{\perp}$. If $\mathcal{H}=\begin{tikzpicture}[baseline=-0.3ex]\dynkin[ply=2,rotate=90,fold radius = 2mm]{A}{oo}\end{tikzpicture} \hspace{0.2mm}\times\hspace{-0.2mm} \! \simcal{\begin{tikzpicture}[baseline=-0.3ex]\dynkin[ply=2,rotate=90,fold radius = 2mm]{A}{**}\end{tikzpicture}}\!\hspace{0.4mm}$ then $\Phi_s=\big\langle  \alpha_0^1 \big\rangle$. For both of these cases, $\Phi_s \cong A_1$ and $r(\Phi_s)=r_s(\mathcal{H}) =1$. So $\Phi_a =(\Phi_s)^{\perp} \cong C_3$ since $\Phi_s$ spans $E_s$.

\vspace{2mm}\noindent If $\mathcal{H}$ is either $\begin{tikzpicture}[baseline=-1.0ex]\dynkin[rotate=180]{C}{*o**}\end{tikzpicture}$, $\!\!\!\!\shimcal{\begin{tikzpicture}[scale=0.8,baseline=-1.7ex] \dynkin[ply=3,rotate=270,fold radius = 2mm]{D}{*o**} \end{tikzpicture}}\!\!\!\!$ or $\begin{tikzpicture}[baseline=-1.0ex]\dynkin[rotate=180]{C}{*o*}\end{tikzpicture} \hspace{-0.3mm}\times\hspace{-0.3mm} \dynkin{A}{*}$ then $\Phi_s=(\Psi_0)^{\perp}$. \textcolor{white}{$\!\!\!\Big($}If $\mathcal{H}=\begin{tikzpicture}[baseline=-0.3ex]\dynkin[ply=2,rotate=90,fold radius = 2mm]{A}{**}\end{tikzpicture} \hspace{0.2mm}\times\hspace{-0.2mm} \! \simcal{\begin{tikzpicture}[baseline=-0.3ex]\dynkin[ply=2,rotate=90,fold radius = 2mm]{A}{oo}\end{tikzpicture}}\!\hspace{0.4mm}$ then $\Phi_s=\big\langle  \alpha_0^2 \big\rangle$. For each of these cases, $\Phi_s \cong \smash{\widetilde{A_1}}$ and $r(\Phi_s)=r_s(\mathcal{H}) =1$. Hence $\Phi_a =(\Phi_s)^{\perp} \cong B_3$.

\vspace{2mm}\noindent If $\mathcal{H}=\begin{tikzpicture}[baseline=-0.3ex]\dynkin[ply=2,rotate=90,fold radius = 2mm]{A}{oo}\end{tikzpicture} \hspace{0.2mm}\times\hspace{-0.2mm} \! \simcal{\begin{tikzpicture}[baseline=-0.3ex]\dynkin[ply=2,rotate=90,fold radius = 2mm]{A}{oo}\end{tikzpicture}}\!\hspace{0.4mm}$ then $\Phi_s=\big\langle  \alpha_0^1,\alpha_0^2 \big\rangle \cong A_1\widetilde{A_1}$ and $r(\Phi_s)=r_s(\mathcal{H}) =2$. Hence $\Phi_a =(\Phi_s)^{\perp} \cong \smash{A_1\widetilde{A_1}}$.

\vspace{2mm}\noindent If $\mathcal{H}=\begin{tikzpicture}[scale=0.8,baseline=-2.2ex] \dynkin[ply=3,rotate=270]{D}{oooo} \end{tikzpicture}$ then $\Phi_s= \big\langle \alpha_0, \beta \big\rangle$, where $\beta=\smash{\small\makecell{$\small0$ \\ $\small0~1~0$}} \in \Delta$. Observe that $\Phi_s \cong A_2$ and $r(\Phi_s)=r_s(\mathcal{H}) =2$. Hence $\Phi_a =(\Phi_s)^{\perp} \cong \smash{\widetilde{A_2}}$.

\vspace{2mm}\noindent If $\mathcal{H}=\!\!\!\shimcal{\begin{tikzpicture}[scale=0.8,baseline=-1.7ex] \dynkin[ply=3,rotate=270,fold radius = 2mm]{D}{oooo} \end{tikzpicture}}\!\!\!$ then $\Phi_s= \big\langle \alpha_0, \gamma \big\rangle$, where $\gamma=\smash{\small\makecell{$\small0$ \\ $\small0~1~0$}} \in \Delta$. \textcolor{white}{$\!\!\!\Big($}Observe that $\Phi_s \cong \smash{\widetilde{A_2}}$ and $r(\Phi_s)=r_s(\mathcal{H}) =2$. Hence $\Phi_a =(\Phi_s)^{\perp} \cong A_2$.

\vspace{2mm}\noindent For all remaining indices $\mathcal{H}$ listed in the fourth column of Table \ref{F_4}, we observe that $\Psi_0$ is a parabolic subsystem of $\Phi$ and we use Lemma \ref{complemma} to check that $r(\Phi_a)=r(\Psi_0)$. So $\Phi_a=\Psi_0$ since $\Phi_a$ is a parabolic subsystem of $\Phi$ that contains $\Psi_0$. That is, $\mathcal{H}$ is independent in $\Phi$.

\subsection*{Appendix C. $\Phi=E_6 = \dynkin[scale=2,labels={\alpha_1,\alpha_...,\alpha_6}]{E}{6}$}\label{E6}

\noindent Let $\Phi=E_6$. We use the following \textit{standard} parametrisation in $\R^6$ of $\Phi$. Let $\alpha_1=\big(-\sum_{i=1}^5 e_i+\sqrt{3}e_6\big)/2$, $\alpha_2=e_4-e_5$, $\alpha_3=e_4+e_5$, $\alpha_4=e_3-e_4$, $\alpha_5=e_2-e_3$ and $\alpha_6=e_1-e_2$. The highest root of $\Phi$ is $\big(\sum_{i=1}^4 e_i-e_5+\sqrt{3}e_6\big)/2$. The Weyl group of $\Phi$ is $W \cong \PSU_4(2) \rtimes \Z_2$ and $\Iso(\Phi)=W \rtimes \langle \rho \rangle$, where $\rho$ acts on the Dynkin diagram of $\Phi$ as the non-trivial graph automorphism. The longest element of $W$ is $-\rho$.

\vspace{2mm}\noindent By Table II of \cite{T}, the isotropic indices of inner type $\Phi$ are $\begin{tikzpicture}[baseline=0.4ex]\dynkin{E}{I}\end{tikzpicture}$, $\begin{tikzpicture}[baseline=0.4ex]\dynkin{E}{*o*o**}\end{tikzpicture}$ and $\begin{tikzpicture}[baseline=0.4ex]\dynkin{E}{o****o}\end{tikzpicture}$, and those of outer type are $\begin{tikzpicture}[baseline=-0.6ex]\dynkin[fold radius=2mm]{E}{II}\end{tikzpicture}$, $\begin{tikzpicture}[baseline=-0.6ex]\dynkin[ply=2,fold radius=2mm]{E}{*o*o**}\end{tikzpicture}$, $\begin{tikzpicture}[baseline=-0.6ex]\dynkin[ply=2,fold radius=2mm]{E}{oo***o}\end{tikzpicture}$, $\begin{tikzpicture}[baseline=-0.6ex]\dynkin[ply=2,fold radius=2mm]{E}{o****o}\end{tikzpicture}$ and $\begin{tikzpicture}[baseline=-0.6ex]\dynkin[ply=2,fold radius=2mm]{E}{*o****}\end{tikzpicture}$.

\vspace{2mm}\noindent Recall from Table \ref{exceptlist} that $\Delta$ is one of $(A_2)^3$, $D_5$, $A_5A_1$ or $D_4$. In the fourth column of Table \ref{E_6} we list all indices $\mathcal{H}$ of type $\Delta$ that are maximal in $\Phi$. For each of those indices $\mathcal{H}$, we compute the $W$-conjugacy class of $\Phi_a$ in $\Phi$.

\vspace{2mm}\noindent If $\mathcal{H}$ is one of $\begin{tikzpicture}[baseline=-0.6ex,rotate=180]\dynkin[ply=2,fold radius=2mm]{D}{*o***}\end{tikzpicture} \hspace{-0.3mm}\times\hspace{-0.3mm} \mathcal{T}^1_0$, $\begin{tikzpicture}[baseline=-0.6ex]\dynkin[ply=2,fold radius=2mm]{A}{o***o}\end{tikzpicture} \hspace{-0.3mm}\times\hspace{-0.3mm} \dynkin{A}{*}$, $\smash{\begin{tikzpicture}[scale=0.8,baseline=-2.2ex] \dynkin[ply=3,rotate=270,fold radius = 2mm]{D}{*o**} \end{tikzpicture}}\hspace{-0.5mm}\times\hspace{-0.3mm} \mathcal{T}^2_0$ (with $\Pi=\Z_3$ or $S_3$) or $\smash{\begin{tikzpicture}[scale=0.8,baseline=-2.2ex] \dynkin[ply=3,rotate=270,fold radius = 2mm]{D}{*o**} \end{tikzpicture}}\hspace{-0.5mm}\times\hspace{-0.3mm} \mathcal{T}^2_0$ (with $\Pi=\Z_6$ or $\Z_2 \hspace{-0.5mm}\times\hspace{-0.5mm} S_3$) then $\Phi_s=(\Psi_0)^{\perp}$. If $\mathcal{H}$ is one of $\begin{tikzpicture}[baseline=0.3ex] 
\dynkin[name=upper]{A}{**}
\node (current) at ($(upper root 1)+(-0,3.1mm)$) {};
\dynkin[at=(current),name=lower]{A}{**}
\begin{scope}[on background layer]
\foreach \i in {1,2}%
{%
\draw[/Dynkin diagram/fold style]
($(upper root \i)$) -- ($(lower
root \i)$);%
}%
\end{scope}
\end{tikzpicture} \hspace{-0.3mm}\times\hspace{-0.3mm} \begin{tikzpicture}[baseline=-0.3ex]\dynkin[ply=2,fold radius = 2mm,rotate=90]{A}{oo}\end{tikzpicture}$, $(\begin{tikzpicture}[baseline=-0.3ex]\dynkin[ply=2,fold radius = 2mm,rotate=90]{A}{**}\end{tikzpicture})^2 \hspace{-0.3mm}\times\hspace{-0.3mm} \begin{tikzpicture}[baseline=-0.3ex]\dynkin[ply=2,fold radius = 2mm,rotate=90]{A}{oo}\end{tikzpicture}$ or $\smash{\begin{tikzpicture}[baseline=0.3ex] 
\dynkin[name=upper]{A}{**}
\node (current) at ($(upper root 1)+(-0,3.1mm)$) {};
\dynkin[at=(current),name=lower]{A}{**}
\begin{scope}[on background layer]
\foreach \i in {1,2}%
{%
\draw[/Dynkin diagram/fold style]
($(upper root \i)$) -- ($(lower
root \i)$);%
}%
\draw[/Dynkin diagram/fold style] ($(upper root 1)$) -- ($(lower root 2)$);
\draw[/Dynkin diagram/fold style] ($(upper root 2)$) -- ($(lower root 1)$);
\end{scope}
\end{tikzpicture}} \hspace{-0.3mm}\times\hspace{-0.3mm} \begin{tikzpicture}[baseline=-0.3ex]\dynkin[ply=2,fold radius = 2mm,rotate=90]{A}{oo}\end{tikzpicture}$ then $\Phi_s=\big\langle  \alpha_0^3 \big\rangle$. For each of these cases, $\Phi_s \cong A_1$ and $r(\Phi_s)=r_s(\mathcal{H}) =1$. Hence $\Phi_a =(\Phi_s)^{\perp} \cong A_5$ since $\Phi_s$ spans $E_s$.

\vspace{2mm}\noindent If $\mathcal{H}$ is either $\begin{tikzpicture}[scale=0.8,baseline=-2.2ex] \dynkin[ply=3,rotate=270,fold radius = 2mm]{D}{oooo} \end{tikzpicture}\hspace{-0.5mm}\times\hspace{-0.3mm} \mathcal{T}^2_0$ (with $\Pi=\Z_3$ or $S_3$) or $\begin{tikzpicture}[scale=0.8,baseline=-2.2ex] \dynkin[ply=3,rotate=270,fold radius = 2mm]{D}{oooo} \end{tikzpicture}\hspace{-0.5mm}\times\hspace{-0.3mm} \mathcal{T}^2_0$ (with $\Pi=\Z_6$ or $\Z_2 \times S_3$) then $\Phi_s= \big\langle \alpha_0, \beta \big\rangle$, where $\beta=\smash{\small\makecell{$\small0$ \\ $\small0~1~0$}} \in \Delta$. \textcolor{white}{$\!\!\!\Big($}For both of these cases, $\Phi_s \cong A_2$ and $r(\Phi_s)=r_s(\mathcal{H}) =2$. Hence $\Phi_a =(\Phi_s)^{\perp} \cong (A_2)^2$.

\vspace{2mm}\noindent If $\mathcal{H}=\smash{\begin{tikzpicture}[baseline=-0.6ex]\dynkin[ply=2,fold radius=2mm,rotate=180]{D}{*o*oo}\end{tikzpicture}} \hspace{-0.3mm}\times\hspace{-0.3mm} \mathcal{T}^1_0$ then $\Phi_s=\smash{(\Psi_0)^{\perp_{\Psi}}}$. If $\mathcal{H}=\smash{\begin{tikzpicture}[baseline=-0.6ex]\dynkin[ply=2,fold radius=2mm]{A}{oo*oo}\end{tikzpicture} \hspace{-0.3mm}\times\hspace{-0.3mm} \dynkin{A}{*}}$ then $\Phi_s=\big\langle \alpha_0^1,\gamma \big\rangle$, where $\gamma=0~1~1~1~0 \in \Psi^1$. If $\mathcal{H}=\begin{tikzpicture}[baseline=-0.3ex]\dynkin[ply=2,fold radius = 2mm,rotate=90]{A}{**}\end{tikzpicture} \hspace{-0.3mm}\times\hspace{-0.3mm} (\begin{tikzpicture}[baseline=-0.3ex]\dynkin[ply=2,fold radius = 2mm,rotate=90]{A}{oo}\end{tikzpicture})^2$ then $\Phi_s=\big\langle \alpha_0^2,\alpha_0^3 \big\rangle$. For each of these cases, $\Phi_s \cong (A_1)^2$ and $r(\Phi_s)=r_s(\mathcal{H}) =2$. Hence $\Phi_a =(\Phi_s)^{\perp} \cong A_3$.

\vspace{2mm}\noindent If $\mathcal{H}=(\begin{tikzpicture}[baseline=-0.3ex]\dynkin[ply=2,fold radius = 2mm,rotate=90]{A}{oo}\end{tikzpicture})^3$ then $\Phi_s=\big\langle \alpha_0^1,\alpha_0^2,\alpha_0^3 \big\rangle$. Observe that $\Phi_s \cong (A_1)^3$ and $r(\Phi_s)=r_s(\mathcal{H}) =3$. Hence $\Phi_a =(\Phi_s)^{\perp} \cong A_1$.

\vspace{2mm}\noindent If $\mathcal{H}$ is either $\dynkin{A}{*o*o*} \hspace{-0.5mm}\times\hspace{-0.5mm} \dynkin{A}{*}$ or $\smash{\begin{tikzpicture}[baseline=-0.6ex]\dynkin[ply=2,fold radius=2mm]{A}{*o*o*}\end{tikzpicture} \hspace{-0.3mm}\times\hspace{-0.3mm} \dynkin{A}{*}}$ then $\Phi_a \cong D_4$ as this is the only rank $4$ parabolic subsystem of $\Phi$ that contains $\Psi_0 \cong (A_1)^4$. 

\vspace{2mm}\noindent If $\mathcal{H}=\smash{\begin{tikzpicture}[baseline=0.3ex] 
\dynkin[name=upper]{A}{oo}
\node (current) at ($(upper root 1)+(-0,3.1mm)$) {};
\dynkin[at=(current),name=lower]{A}{oo}
\begin{scope}[on background layer]
\foreach \i in {1,2}%
{%
\draw[/Dynkin diagram/fold style]
($(upper root \i)$) -- ($(lower
root \i)$);%
}%
\end{scope}
\end{tikzpicture} \hspace{-0.3mm}\times\hspace{-0.3mm} \begin{tikzpicture}[baseline=-0.3ex]\dynkin[ply=2,fold radius = 2mm,rotate=90]{A}{**}\end{tikzpicture}}$ then, by Lemma \ref{complemma}, $\Phi_a$ is a rank $4$ parabolic subsystem of $\Phi$ that contains $\Psi_0 \cong A_2$ and $\Aut_{W_{\Phi_a}}\hspace{-0.5mm}\!\big((W_{\Phi_a})_0\big)$ acts non-trivially on $\Delta_0$. The only possibility for such a subsystem $\Phi_a$ is of type $D_4$.

\vspace{2mm}\noindent For the next five cases, without loss of generality, we choose a representative of the $W$-conjugacy class of $\Delta \cong (A_2)^3$ in $\Phi$ using the standard parametrisation in $\R^6$ of $\Phi$. We choose $\Delta = \big\{\big(-\sum_{i=1}^4 e_i+e_5-\sqrt{3}e_6\big)/2,\alpha_2,\alpha_1,\alpha_3,\alpha_5,\alpha_6\big\}$. We will need the following roots in $\Phi$. Let $\delta_1:=(-e_1-e_2+e_3+e_4+e_5-\sqrt{3}e_6)/2$ and $\delta_2:=(-e_1-e_2+e_3+e_4-e_5+\sqrt{3}e_6)/2$. Let $\epsilon_1:=(-e_1-e_2-e_3+e_4+e_5+\sqrt{3}e_6)/2$ and $\epsilon_2:=(-e_1-e_2-e_3+e_4-e_5-\sqrt{3}e_6)/2$. Let $\zeta:=(-e_1-e_2+e_3-e_4-e_5-\sqrt{3}e_6)/2$.

\vspace{2mm}\noindent If $\mathcal{H}=\begin{tikzpicture}[baseline=0.5ex,rotate=90] 
\dynkin[name=upper]{A}{oo}
\node (current) at ($(upper root 1)+(-0,3.4mm)$) {};
\dynkin[at=(current),name=middle]{A}{oo}
\node (currenta) at ($(upper root 1)+(-0,7.4mm)$) {};
\dynkin[at=(currenta),name=lower]{A}{oo}
\begin{scope}[on background layer]
\foreach \i in {1,2}%
{%
\draw[/Dynkin diagram/fold style]
($(upper root \i)$) -- ($(lower
root \i)$);%
}%
\end{scope}
\end{tikzpicture}$ then the roots $\alpha_4$ and $e_2+e_3$ of $\Phi$ are fixed by $\Pi$ and that $r_s(\mathcal{H}) =2$. So $\big\{\alpha_4,e_2+e_3\big\}$ is a basis for $E_s$. Hence $\Phi_a = \big\langle e_1-e_5,\delta_1,e_1+e_5,\delta_2 \big\rangle \cong (A_2)^2$.

\vspace{2mm}\noindent If $\mathcal{H}=\smash{\begin{tikzpicture}[baseline=0.35ex,rotate=90] 
\dynkin[name=upper]{A}{oo}
\node (current) at ($(upper root 1)+(-0,3.4mm)$) {};
\dynkin[at=(current),name=middle]{A}{oo}
\node (currenta) at ($(upper root 1)+(-0,7.4mm)$) {};
\dynkin[at=(currenta),name=lower]{A}{oo}
\begin{scope}[on background layer]
\draw[/Dynkin diagram/fold style]
($(upper root 1)$) -- ($(lower
root 2)$);
\draw[/Dynkin diagram/fold style]
($(upper root 2)$) -- ($(lower
root 1)$);
\draw[/Dynkin diagram/fold style]
($(upper root 1)$) -- ($(middle
root 2)$);
\draw[/Dynkin diagram/fold style]
($(upper root 2)$) -- ($(middle
root 1)$);
\draw[/Dynkin diagram/fold style]
($(middle root 1)$) -- ($(lower
root 2)$);
\draw[/Dynkin diagram/fold style]
($(middle root 2)$) -- ($(lower
root 1)$);
\end{scope}
\end{tikzpicture}}$ then the vector $-e_2-2e_3+e_4$ is fixed by $\Pi$ and $r_s(\mathcal{H}) =1$. So $E_s$ is generated by $-e_2-2e_3+e_4$. Hence $\Phi_a = \big\langle e_1-e_5,\delta_1,e_1+e_5,\delta_2,e_2+e_4 \big\rangle \cong (A_2)^2A_1$.

\vspace{2mm}\noindent If $\mathcal{H}=\smash{\begin{tikzpicture}[baseline=0.3ex] 
\dynkin[name=upper]{A}{oo}
\node (current) at ($(upper root 1)+(-0,3.1mm)$) {};
\dynkin[at=(current),name=lower]{A}{oo}
\begin{scope}[on background layer]
\foreach \i in {1,2}%
{%
\draw[/Dynkin diagram/fold style]
($(upper root \i)$) -- ($(lower
root \i)$);%
}%
\end{scope}
\end{tikzpicture} \hspace{-0.3mm}\times\hspace{-0.3mm} \begin{tikzpicture}[baseline=-0.3ex]\dynkin[ply=2,fold radius = 2mm,rotate=90]{A}{oo}\end{tikzpicture}}$ then $\epsilon_1$, $\epsilon_2$ and $e_1-e_3$ are fixed by $\Pi$ and $r_s(\mathcal{H}) =3$. So $\big\{\epsilon_1,\epsilon_2,e_1-e_3\big\}$ is a basis for $E_s$. Hence $\Phi_a = \big\langle e_3+e_4,\zeta,e_1+e_3 \big\rangle \cong (A_2)^2$.

\vspace{2mm}\noindent If $\mathcal{H}=\smash{\begin{tikzpicture}[baseline=0.3ex] 
\dynkin[name=upper]{A}{oo}
\node (current) at ($(upper root 1)+(-0,3.1mm)$) {};
\dynkin[at=(current),name=lower]{A}{oo}
\begin{scope}[on background layer]
\foreach \i in {1,2}%
{%
\draw[/Dynkin diagram/fold style]
($(upper root \i)$) -- ($(lower
root \i)$);%
}%
\draw[/Dynkin diagram/fold style] ($(upper root 1)$) -- ($(lower root 2)$);
\draw[/Dynkin diagram/fold style] ($(upper root 2)$) -- ($(lower root 1)$);
\end{scope}
\end{tikzpicture} \hspace{-0.3mm}\times\hspace{-0.3mm} \begin{tikzpicture}[baseline=-0.3ex]\dynkin[ply=2,fold radius = 2mm,rotate=90]{A}{**}\end{tikzpicture}}$ then $\epsilon_1+\epsilon_2=-e_1-e_2-e_3+e_4$ is fixed by $\Pi$ and $r_s(\mathcal{H}) =1$. So $E_s$ is generated by $-e_1-e_2-e_3+e_4$. Hence $\Phi_a = \big\langle e_3+e_4,e_2-e_3,\zeta,e_1+e_3 \big\rangle \cong D_4$.

\vspace{2mm}\noindent If $\mathcal{H}=\smash{\begin{tikzpicture}[baseline=0.3ex] 
\dynkin[name=upper]{A}{oo}
\node (current) at ($(upper root 1)+(-0,3.1mm)$) {};
\dynkin[at=(current),name=lower]{A}{oo}
\begin{scope}[on background layer]
\foreach \i in {1,2}%
{%
\draw[/Dynkin diagram/fold style]
($(upper root \i)$) -- ($(lower
root \i)$);%
}%
\draw[/Dynkin diagram/fold style] ($(upper root 1)$) -- ($(lower root 2)$);
\draw[/Dynkin diagram/fold style] ($(upper root 2)$) -- ($(lower root 1)$);
\end{scope}
\end{tikzpicture} \hspace{-0.3mm}\times\hspace{-0.3mm} \begin{tikzpicture}[baseline=-0.3ex]\dynkin[ply=2,fold radius = 2mm,rotate=90]{A}{oo}\end{tikzpicture}}$ then $-e_1-e_2-e_3+e_4$ and $e_1-e_3$ are fixed by $\Pi$ and $r_s(\mathcal{H}) =2$. So $\big\{-e_1-e_2-e_3+e_4,e_1-e_3\big\}$ is a basis for $E_s$. Hence $\Phi_a = \big\langle e_1+e_3,\zeta,e_3+e_4 \big\rangle \cong (A_1)^3$.

\vspace{2mm}\noindent For all remaining indices $\mathcal{H}$ listed in the fourth column of Table \ref{E_6}, we observe that $\Psi_0$ is a parabolic subsystem of $\Phi$ and we use Lemma \ref{complemma} to check that $r(\Phi_a)=r(\Psi_0)$. So $\Phi_a=\Psi_0$ since $\Phi_a$ is a parabolic subsystem of $\Phi$ that contains $\Psi_0$. That is, $\mathcal{H}$ is independent in $\Phi$.

\newpage\newgeometry{top=0.7cm,bottom=0.7cm,left=0cm, right=0cm,foot=4mm}
\begin{longtable}{| c | c | c | c | c | c | c | c | c | c | c |}\caption{Isotropic maximal embeddings of abstract indices of type $E_6$}\label{E_6} \\ \hline

\multicolumn{1}{|c|}{} & \multicolumn{1}{c|}{} & \multicolumn{1}{c|}{} & \multicolumn{1}{c|}{} & \multicolumn{1}{c|}{} & \multicolumn{1}{c|}{} & \multicolumn{1}{c|}{} & \multicolumn{3}{c|}{\small{Special fields}} \\

  $\Delta$ & \small{$\!\!\Stab_{\Iso(\Phi)}(\Delta)\!\!$} & $\Pi$ & $\mathcal{H}$ & $\Psi_0$ & $\Phi_a$ & $\mathcal{G}$ & \small{$\hspace{-0.8mm}\!\cd 1 \!\hspace{-1mm}$} & $\R$ & $\!\!Q_{\mathfrak{p}}\!\!$ \\ \hline \hline 
\endfirsthead

\multicolumn{11}{|c|}%
{{\bfseries \tablename\ \thetable{} -- continued from previous page}} \\
\hline
\multicolumn{1}{|c|}{} & \multicolumn{1}{c|}{} & \multicolumn{1}{c|}{} & \multicolumn{1}{c|}{} & \multicolumn{1}{c|}{} & \multicolumn{1}{c|}{} & \multicolumn{1}{c|}{} & \multicolumn{3}{c|}{\small{Special fields}} \\

  $\Delta$ & \small{$\!\!\Stab_{\Iso(\Phi)}(\Delta)\!\!$} & $\Pi$ & $\mathcal{H}$ & $\Psi_0$ & $\Phi_a$ & $\mathcal{G}$ & \small{$\hspace{-0.8mm}\!\cd 1 \!\hspace{-1mm}$} & $\R$ & $\!\!Q_{\mathfrak{p}}\!\!$ \\ \hline \hline \endhead 
\hline \multicolumn{11}{|c|}{{Continued on next page}} \\ \hline \endfoot
 \endlastfoot

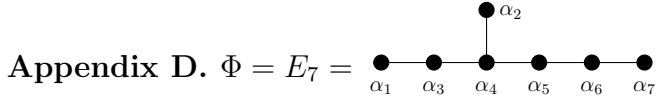
\multirow{38.5}{*}{$\!(A_2)^3\!$} & \multirow{38.5}{*}{\makecell{\small{$\!\!\big\langle  \psi,r,\rho \big\rangle \!\!$} \\  $\!\cong \Z_2 \hspace{-0.3mm}\times\hspace{-0.3mm} S_3$}} & \multirow{3.6}{*}{$1$} & $(\dynkin{A}{oo})^3$ & $\varnothing$ & $\varnothing$ & $\begin{tikzpicture}[baseline=0.4ex]\dynkin{E}{I}\end{tikzpicture}$ & $\checkmark$ & $\checkmark$ & $~\checkmark$ \topstrut \\ \cdashline{4-11} 
&&& $\!\dynkin{A}{**} \hspace{-0.3mm}\times\hspace{-0.3mm} (\dynkin{A}{oo})^2\!$ & $A_2$ & $A_2$ &&&& \topstrut \\ \cdashline{4-11} 
&&& $\!(\dynkin{A}{**})^2 \hspace{-0.3mm}\times\hspace{-0.3mm} \dynkin{A}{oo}\!$ & $(A_2)^2$ & $(A_2)^2$ & $\begin{tikzpicture}[baseline=0.4ex]\dynkin{E}{*o*o**}\end{tikzpicture}$ & $\xmark$ & $\xmark$ & $~\checkmark$ \topstrut \\ \cline{3-11} 

& & \multirow{6.1}{*}{\makecell{$\!\langle \psi\rho \rangle\!$ \\ $\!\cong \! \Z_2\!$}} & $\begin{tikzpicture}[baseline=0.3ex] 
\dynkin[name=upper]{A}{oo}
\node (current) at ($(upper root 1)+(-0,3.1mm)$) {};
\dynkin[at=(current),name=lower]{A}{oo}
\begin{scope}[on background layer]
\foreach \i in {1,2}%
{%
\draw[/Dynkin diagram/fold style]
($(upper root \i)$) -- ($(lower
root \i)$);%
}%
\end{scope}
\end{tikzpicture} \hspace{-0.3mm}\times\hspace{-0.3mm} \begin{tikzpicture}[baseline=-0.3ex]\dynkin[ply=2,fold radius = 2mm,rotate=90]{A}{oo}\end{tikzpicture}$ & $\varnothing$ & $(A_1)^3$ & $\begin{tikzpicture}[baseline=0.4ex]\dynkin{E}{I}\end{tikzpicture}$ & $\checkmark$ & $\checkmark$ & $\checkmark$ \\
\end{longtable}\restoregeometry

\subsection*{Appendix D. $\Phi=E_7 = \dynkin[scale=2,labels={\alpha_1,\alpha_...,\alpha_7}]{E}{7}$}\label{E7}

\noindent Let $\Phi=E_7$. We use the following \textit{standard} parametrisation in $\R^7$ of $\Phi$. Let $\alpha_1=\big(-\sum_{i=1}^6 e_i+\sqrt{2}e_7\big)/2$, $\alpha_2=e_5-e_6$, $\alpha_3=e_5+e_6$, $\alpha_4=e_4-e_5$, $\alpha_5=e_3-e_4$, $\alpha_6=e_2-e_3$ and $\alpha_7=e_1-e_2$. The highest root of $\Phi$ is $\sqrt{2}e_7$. The Weyl group of $\Phi$ is $W \cong \Z_2 \times \PSp_6(2)$ and the longest element of $W$ is $-1$.

\vspace{1mm}\noindent By Table II of \cite{T}, the isotropic indices of type $\Phi$ are $\begin{tikzpicture}[baseline=0.4ex]\dynkin{E}{V}\end{tikzpicture}$, $\begin{tikzpicture}[baseline=0.4ex]\dynkin{E}{VI}\end{tikzpicture}$, $\begin{tikzpicture}[baseline=0.4ex]\dynkin{E}{VII}\end{tikzpicture}$, $\begin{tikzpicture}[baseline=0.4ex]\dynkin{E}{o****o*}\end{tikzpicture}$, $\begin{tikzpicture}[baseline=0.4ex]\dynkin{E}{*****o*}\end{tikzpicture}$, $\begin{tikzpicture}[baseline=0.4ex]\dynkin{E}{o******}\end{tikzpicture}$ and $\begin{tikzpicture}[baseline=0.4ex]\dynkin{E}{******o}\end{tikzpicture}$. All of these indices are of inner type since $\Iso(\Phi)=W$.

\vspace{2mm}\noindent Recall from Table \ref{exceptlist} that $\Delta$ is one of $A_5A_2$, $D_6A_1$, $D_4(A_1)^3$, $A_7$, $E_6$ or $(A_1)^7$. In the fourth column of Table \ref{E_7} we list all indices $\mathcal{H}$ of type $\Delta$ that are maximal in $\Phi$. For each of those indices $\mathcal{H}$, we compute the $W$-conjugacy class of $\Phi_a$ in $\Phi$.

\vspace{2mm}\noindent If $\mathcal{H}$ is one of $\begin{tikzpicture}[baseline=-0.6ex]\dynkin[ply=2,fold radius = 2mm]{A}{o***o}\end{tikzpicture} \hspace{-0.5mm}\times\hspace{-0.5mm} \begin{tikzpicture}[baseline=-0.6ex]\dynkin[ply=2,fold radius = 2mm]{A}{**}\end{tikzpicture}$, $\begin{tikzpicture}[scale=.175,baseline=0.3ex]
    \foreach \y in {0,...,3}
    \draw[thin,xshift=\y cm] (\y cm,0) ++(.3 cm, 0) -- +(14 mm,0);
    \draw[thin,fill=black] (0 cm,0) circle (3 mm);
    \draw[thin,fill=black] (2 cm,0) circle (3 mm);
    \draw[thin,fill=black] (4 cm,0) circle (3 mm);
    \draw[thin] (6 cm,0) circle (3 mm);
    \draw[thin,fill=black] (8 cm,0) circle (3 mm);
    \draw[thin,fill=black] (2 cm,2 cm) circle (3 mm);
    \draw[thin, fill=black] (2 cm, 3mm) -- +(0, 1.4 cm);
  \end{tikzpicture} \times\hspace{-0.5mm} \dynkin{A}{1}$, $\smash{\begin{tikzpicture}[scale=0.8,baseline=-2.2ex] \dynkin[ply=3,rotate=270,fold radius = 2mm]{D}{*o**} \end{tikzpicture}} \hspace{-0.5mm}\times\hspace{-0.5mm} \begin{tikzpicture}[baseline=-0.5ex,rotate=90] 
\dynkin[name=upper]{A}{*}
\node (current) at ($(upper root 1)+(-0,3.4mm)$) {};
\dynkin[at=(current),name=middle]{A}{*}
\node (currenta) at ($(upper root 1)+(-0,7.4mm)$) {};
\dynkin[at=(currenta),name=lower]{A}{*}
\begin{scope}[on background layer]
\foreach \i in {1}%
{%
\draw[/Dynkin diagram/fold style]
($(upper root \i)$) -- ($(lower
root \i)$);%
}%
\end{scope}
\end{tikzpicture}$, $\dynkin{A}{***o***}$, $\begin{tikzpicture}[baseline=-0.6ex]\dynkin[ply=2,fold radius = 2mm]{A}{***o***}\end{tikzpicture}$, $\begin{tikzpicture}[baseline=-0.6ex]\dynkin[ply=2,fold radius = 2mm]{A}{o*****o}\end{tikzpicture}$ or $\smash{\begin{tikzpicture}[baseline=-0.6ex]\dynkin[ply=2,fold radius=2mm]{E}{*o****}\end{tikzpicture}} \hspace{-0.3mm}\times\hspace{-0.5mm} \mathcal{T}^1_0$ then $\Phi_s=(\Psi_0)^{\perp}$ (these last two indices satisfy $\Psi_0 \cong (A_5)''$). If $\mathcal{H}= \begin{tikzpicture}[baseline=-0.6ex]\dynkin[ply=2,fold radius = 2mm]{A}{*****}\end{tikzpicture} \hspace{-0.5mm}\times\hspace{-0.5mm} \begin{tikzpicture}[baseline=-0.6ex]\dynkin[ply=2,fold radius = 2mm]{A}{oo}\end{tikzpicture}$ then $\Phi_s=\big\langle  \alpha_0^2 \big\rangle$. For each of these cases, $\Phi_s \cong A_1$ and $r(\Phi_s)=r_s(\mathcal{H}) =1$. Hence $\Phi_a =(\Phi_s)^{\perp} \cong D_6$ since $\Phi_s$ spans $E_s$.

\vspace{2mm}\noindent If $\mathcal{H}=\begin{tikzpicture}[scale=.175,baseline=0.3ex]
    \foreach \y in {0,...,3}
    \draw[thin,xshift=\y cm] (\y cm,0) ++(.3 cm, 0) -- +(14 mm,0);
    \draw[thin,fill=black] (0 cm,0) circle (3 mm);
    \draw[thin] (2 cm,0) circle (3 mm);
    \draw[thin,fill=black] (4 cm,0) circle (3 mm);
    \draw[thin] (6 cm,0) circle (3 mm);
    \draw[thin,fill=black] (8 cm,0) circle (3 mm);
    \draw[thin,fill=black] (2 cm,2 cm) circle (3 mm);
    \draw[thin, fill=black] (2 cm, 3mm) -- +(0, 1.4 cm);
  \end{tikzpicture} \times\hspace{-0.5mm} \dynkin{A}{*}$ then $\Phi_s=(\Psi_0)^{\perp}$. If $\mathcal{H}=\begin{tikzpicture}[baseline=-0.6ex]\dynkin[ply=2,fold radius = 2mm]{E}{oo***o}\end{tikzpicture} \hspace{-0.3mm}\times\hspace{-0.5mm} \mathcal{T}^1_0$ then $\Phi_s=(\Psi_0)^{\perp_{\Psi}}$. If $\mathcal{H}=\begin{tikzpicture}[baseline=-0.6ex]\dynkin[ply=2,fold radius = 2mm]{A}{o***o}\end{tikzpicture} \hspace{-0.5mm}\times\hspace{-0.5mm} \begin{tikzpicture}[baseline=-0.6ex]\dynkin[ply=2,fold radius = 2mm]{A}{oo}\end{tikzpicture}$ then $\Phi_s= \big\langle \alpha_0^1, \alpha_0^2\big\rangle$. If $\mathcal{H}=\smash{\begin{tikzpicture}[baseline=-0.6ex]\dynkin[ply=2,fold radius = 2mm]{A}{oo*oo}\end{tikzpicture} \hspace{-0.5mm}\times\hspace{-0.5mm} \begin{tikzpicture}[baseline=-0.6ex]\dynkin[ply=2,fold radius = 2mm]{A}{**}\end{tikzpicture}}$ then $\Phi_s= \big\langle \alpha_0^1, \beta\big\rangle$, where $\beta = 0~1~1~1~0 \in \Psi^1$. If $\mathcal{H}=\smash{\begin{tikzpicture}[baseline=-0.6ex]\dynkin[ply=2,fold radius = 2mm]{A}{oo***oo}\end{tikzpicture}}$ then $\Phi_s= \big\langle \alpha_0, \gamma\big\rangle$, where $\gamma = 0~1~1~1~1~1~0 \in \Psi$. If $\mathcal{H}=\smash{\begin{tikzpicture}[baseline=-0.6ex]\dynkin[ply=2,fold radius = 2mm]{A}{*o*o*o*}\end{tikzpicture}}$ then $(\Psi_0)^{\perp} =\big\langle \delta_1,\delta_2,\delta_3\big\rangle\cong {\big((A_1)^3\big)}'$, where the non-trivial element of $\Pi \cong \Z_2$ fixes $\delta_1$ and $\delta_2$ pointwise and sends $\delta_3 \mapsto -\delta_3$, and so $\Phi_s= \big\langle \delta_1,\delta_2\big\rangle$. For each of these cases, $\Phi_s \cong (A_1)^2$ and $r(\Phi_s)=r_s(\mathcal{H}) =2$. Hence $\Phi_a =(\Phi_s)^{\perp} \cong D_4A_1$.

\vspace{2mm}\noindent If $\mathcal{H}=\begin{tikzpicture}[baseline=-0.6ex]\dynkin[ply=2,fold radius=2mm]{E}{*o*o**}\end{tikzpicture} \hspace{-0.3mm}\times\hspace{-0.5mm} \mathcal{T}^1_0$ then $\Phi_s=(\Psi_0)^{\perp}$. If $\mathcal{H}=\smash{\begin{tikzpicture}[scale=0.8,baseline=-2.2ex] \dynkin[ply=3,rotate=270,fold radius = 2mm]{D}{oooo} \end{tikzpicture}} \hspace{-0.5mm}\times\hspace{-0.5mm} \begin{tikzpicture}[baseline=-0.5ex,rotate=90] 
\dynkin[name=upper]{A}{*}
\node (current) at ($(upper root 1)+(-0,3.4mm)$) {};
\dynkin[at=(current),name=middle]{A}{*}
\node (currenta) at ($(upper root 1)+(-0,7.4mm)$) {};
\dynkin[at=(currenta),name=lower]{A}{*}
\begin{scope}[on background layer]
\foreach \i in {1}%
{%
\draw[/Dynkin diagram/fold style]
($(upper root \i)$) -- ($(lower
root \i)$);%
}%
\end{scope}
\end{tikzpicture}$ then $\Phi_s= \big\langle \alpha_0^1, \epsilon \big\rangle$, where $\epsilon=\smash{\small\makecell{$\small0$ \\ $\small0~1~0$}} \in \Psi^1$. \textcolor{white}{$\!\!\!\Big($}For both of these cases, $\Phi_s \cong A_2$ and $r(\Phi_s)=r_s(\mathcal{H}) =2$. Hence $\Phi_a =(\Phi_s)^{\perp} \cong (A_5)'$.

\vspace{2mm}\noindent If $\mathcal{H}$ is either $\smash{\begin{tikzpicture}[scale=.175,baseline=0.3ex]
    \foreach \y in {0,...,3}
    \draw[thin,xshift=\y cm] (\y cm,0) ++(.3 cm, 0) -- +(14 mm,0);
    \draw[thin,fill=black] (0 cm,0) circle (3 mm);
    \draw[thin] (2 cm,0) circle (3 mm);
    \draw[thin,fill=black] (4 cm,0) circle (3 mm);
    \draw[thin] (6 cm,0) circle (3 mm);
    \draw[thin,fill=black] (8 cm,0) circle (3 mm);
    \draw[thin] (2 cm,2 cm) circle (3 mm);
    \draw[thin] (2 cm, 3mm) -- +(0, 1.4 cm);
  \end{tikzpicture} \times\hspace{-0.5mm} \dynkin{A}{*}}$ or $\dynkin{A}{*o*o*o*}$ then $\Phi_s=(\Psi_0)^{\perp}$. If $\mathcal{H}=\smash{\begin{tikzpicture}[baseline=-0.6ex]\dynkin[ply=2,fold radius = 2mm]{A}{ooooo}\end{tikzpicture} \hspace{-0.5mm}\times\hspace{-0.5mm} \begin{tikzpicture}[baseline=-0.6ex]\dynkin[ply=2,fold radius = 2mm]{A}{**}\end{tikzpicture}}$ then $\Phi_s= \big\langle \alpha_0^1, \zeta_1,\zeta_2\big\rangle$, where $\zeta_1 = 0~1~1~1~0 \in \Psi^1$ and $\zeta_2 = 0~0~1~0~0 \in \Psi^1$. For each of these cases, $\Phi_s \cong {\big((A_1)^3\big)}'$ and $r(\Phi_s)=r_s(\mathcal{H}) =3$. Hence $\Phi_a =(\Phi_s)^{\perp} \cong D_4$.

\vspace{2mm}\noindent If $\mathcal{H}=\smash{\begin{tikzpicture}[baseline=-0.6ex]\dynkin[ply=2,fold radius = 2mm]{A}{oo*oo}\end{tikzpicture} \hspace{-0.5mm}\times\hspace{-0.5mm} \begin{tikzpicture}[baseline=-0.6ex]\dynkin[ply=2,fold radius = 2mm]{A}{oo}\end{tikzpicture}}$ then $\Phi_s= \big\langle \alpha_0^1, \alpha_0^2,\eta \big\rangle$, where $\eta = 0~1~1~1~0 \in \Psi^1$. If $\mathcal{H}=\smash{\begin{tikzpicture}[baseline=-0.6ex]\dynkin[ply=2,fold radius = 2mm]{A}{ooo*ooo}\end{tikzpicture}}$ then $\Phi_s= \big\langle \alpha_0,\kappa_1,\kappa_2 \big\rangle$, where $\kappa_1 = 0~1~1~1~1~1~0 \in \Psi$ and $\kappa_2 = 0~0~1~1~1~0~0 \in \Psi$. For both of these cases, $\Phi_s \cong {\big((A_1)^3\big)}''$ and $r(\Phi_s)=r_s(\mathcal{H}) =3$. Hence $\Phi_a =(\Phi_s)^{\perp} \cong {\big((A_1)^4\big)}''$.

\vspace{2mm}\noindent If $\mathcal{H}=\begin{tikzpicture}[baseline=-0.6ex]\dynkin[ply=2,fold radius = 2mm]{A}{ooooo}\end{tikzpicture} \hspace{-0.5mm}\times\hspace{-0.5mm} \begin{tikzpicture}[baseline=-0.6ex]\dynkin[ply=2,fold radius = 2mm]{A}{oo}\end{tikzpicture}$ then $\Phi_s= \big\langle \alpha_0^1, \alpha_0^2,\lambda_1,\lambda_2 \big\rangle$, where $\lambda_1 = 0~1~1~1~0 \in \Psi^1$ and $\lambda_2 = 0~0~1~0~0 \in \Psi^1$. Observe that $\Phi_s \cong {\big((A_1)^4\big)}''$ and $r(\Phi_s)=r_s(\mathcal{H}) =4$. Hence $\Phi_a =(\Phi_s)^{\perp} \cong {\big((A_1)^3\big)}''$.

\vspace{2mm}\noindent If $\mathcal{H}=\smash{\begin{tikzpicture}[baseline=-0.6ex]\dynkin[ply=2,fold radius = 2mm]{A}{ooooooo}\end{tikzpicture}}$ then $\Phi_s= \big\langle \alpha_0,\mu_1,\mu_2,\mu_3 \big\rangle \cong {\big((A_1)^4\big)}'$, where $\mu_1 = 0~1~1~1~1~1~0 \in \Psi$, $\mu_1 = 0~0~1~1~1~0~0 \in \Psi$ and $\mu_3 = 0~0~0~1~0~0~0 \in \Psi$. If $\mathcal{H}=\begin{tikzpicture}[baseline=-0.6ex]\dynkin[ply=2,fold radius = 2mm]{E}{oooooo}\end{tikzpicture} \hspace{-0.3mm}\times\hspace{-0.5mm} \mathcal{T}^1_0$ then $\Phi_s= \big\langle \alpha_0,\nu_1,\nu_2,\nu_3 \big\rangle \cong D_4$, \textcolor{white}{$\!\!\!\Big($}where $\nu_1=\smash{\small\makecell{$\small0$ \\ $\small0~0~1~0~0$}} \in \Psi$, $\nu_2=\smash{\small\makecell{$\small1$ \\ $\small0~0~1~0~0$}} \in \Psi$ and $\nu_3=\smash{\small\makecell{$\small0$ \\ $\small1~1~1~1~1$}} \in \Psi$. For both of these cases, $r(\Phi_s)=r_s(\mathcal{H}) =4$. Hence $\Phi_a =(\Phi_s)^{\perp} \cong \smash{{\big((A_1)^3\big)}'}$.

\vspace{2mm}\noindent If $\mathcal{H}$ is either $\dynkin{A}{**o**} \hspace{-0.5mm}\times\hspace{-0.5mm} \dynkin{A}{**}$ or $\smash{\begin{tikzpicture}[baseline=-0.6ex]\dynkin[ply=2,fold radius = 2mm]{A}{**o**}\end{tikzpicture} \hspace{-0.5mm}\times\hspace{-0.5mm} \begin{tikzpicture}[baseline=-0.6ex]\dynkin[ply=2,fold radius = 2mm]{A}{**}\end{tikzpicture}}$ then $\Phi_a \cong E_6$ as this is the only parabolic subsystem of $\Phi$ that contains $\Psi_0 \cong (A_2)^3$. Similarly, if $\mathcal{H}=\begin{tikzpicture}[scale=.175,baseline=0.3ex]
    \foreach \y in {0,...,3}
    \draw[thin,xshift=\y cm] (\y cm,0) ++(.3 cm, 0) -- +(14 mm,0);
    \draw[thin,fill=black] (0 cm,0) circle (3 mm);
    \draw[thin] (2 cm,0) circle (3 mm);
    \draw[thin,fill=black] (4 cm,0) circle (3 mm);
    \draw[thin,fill=black] (6 cm,0) circle (3 mm);
    \draw[thin,fill=black] (8 cm,0) circle (3 mm);
    \draw[thin,fill=black] (2 cm,2 cm) circle (3 mm);
    \draw[thin, fill=black] (2 cm, 3mm) -- +(0, 1.4 cm);
  \end{tikzpicture} \times\hspace{-0.5mm} \dynkin{A}{1}$ then $\Phi_a \cong D_5A_1$ as this is the only parabolic subsystem of $\Phi$ that contains $\Psi_0 \cong A_3(A_1)^3$. 

\vspace{2mm}\noindent If $\mathcal{H}=\begin{tikzpicture}[baseline=-0.5ex,rotate=90] 
\dynkin[name=upper]{A}{o}
\node (current) at ($(upper root 1)+(-0,3.4mm)$) {};
\dynkin[at=(current),name=first]{A}{o}
\node (currenta) at ($(upper root 1)+(-0,7.4mm)$) {};
\dynkin[at=(currenta),name=second]{A}{o}
\node (currentaa) at ($(upper root 1)+(-0,11.4mm)$) {};
\dynkin[at=(currentaa),name=third]{A}{o}
\node (currentaaa) at ($(upper root 1)+(-0,15.4mm)$) {};
\dynkin[at=(currentaaa),name=fourth]{A}{o}
\node (currentaaaa) at ($(upper root 1)+(-0,19.4mm)$) {};
\dynkin[at=(currentaaaa),name=fifth]{A}{o}
\node (currentaaaaa) at ($(upper root 1)+(-0,23.4mm)$) {};
\dynkin[at=(currentaaaaa),name=lower]{A}{o}
\begin{scope}[on background layer]
\foreach \i in {1}%
{%
\draw[/Dynkin diagram/fold style]
($(upper root \i)$) -- ($(lower
root \i)$);%
}%
\end{scope}
\end{tikzpicture}$ then $\Pi$ contains an element of order $7$. By Lemma \ref{complemma}, $W_{\Phi_a}$ contains an element of order $7$. The only possibility for such a subsystem $\Phi_a$ is of type $A_6$.

\vspace{2mm}\noindent If $\mathcal{H}=\smash{\begin{tikzpicture}[scale=0.8,baseline=-2.2ex] \dynkin[ply=3,rotate=270,fold radius = 2mm]{D}{*o**} \end{tikzpicture} \hspace{-0.5mm}\times\hspace{-0.5mm} \begin{tikzpicture}[baseline=-0.5ex,rotate=90] 
\dynkin[name=upper]{A}{o}
\node (current) at ($(upper root 1)+(-0,3.4mm)$) {};
\dynkin[at=(current),name=middle]{A}{o}
\node (currenta) at ($(upper root 1)+(-0,7.4mm)$) {};
\dynkin[at=(currenta),name=lower]{A}{o}
\begin{scope}[on background layer]
\foreach \i in {1}%
{%
\draw[/Dynkin diagram/fold style]
($(upper root \i)$) -- ($(lower
root \i)$);%
}%
\end{scope}
\end{tikzpicture}}$ then, by Lemma \ref{complemma}, $\Phi_a$ is a rank $5$ parabolic subsystem of $\Phi$ that contains $\Psi_0 \cong \big((A_1)^3\big)''$ and $\Aut_{W_{\Phi_a}}(W_{\Delta_0})$ acts transitively on $\Delta_0$. The only possibility for such a subsystem $\Phi_a$ is of type $(A_5)''$.

\vspace{2mm}\noindent If $\mathcal{H}=\begin{tikzpicture}[scale=0.8,baseline=-2.2ex] \dynkin[ply=3,rotate=270,fold radius = 2mm]{D}{****} \end{tikzpicture} \hspace{-0.5mm}\times\hspace{-0.5mm} \begin{tikzpicture}[baseline=-0.5ex,rotate=90] 
\dynkin[name=upper]{A}{o}
\node (current) at ($(upper root 1)+(-0,3.4mm)$) {};
\dynkin[at=(current),name=middle]{A}{o}
\node (currenta) at ($(upper root 1)+(-0,7.4mm)$) {};
\dynkin[at=(currenta),name=lower]{A}{o}
\begin{scope}[on background layer]
\foreach \i in {1}%
{%
\draw[/Dynkin diagram/fold style]
($(upper root \i)$) -- ($(lower
root \i)$);%
}%
\end{scope}
\end{tikzpicture}$ then, by Lemma \ref{complemma}, $\Phi_a$ is a rank $6$ parabolic subsystem of $\Phi$ that contains $\Psi_0 \cong D_4$ and $\Aut_{W_{\Phi_a}}(W_{\Delta_0})$ acts on $\Delta_0$ by triality. The only possibility for such a subsystem $\Phi_a$ is of type $E_6$.

\vspace{2mm}\noindent For the next cases, without loss of generality, we choose a representative of the $W$-conjugacy class of $\Delta$ in $\Phi$ using the standard parametrisation in $\R^7$ of $\Phi$. We will need the following roots in $\Phi$. Let $\xi_1:=(-e_1-e_2+e_3-e_4-e_5+e_6+\sqrt{2}e_7)/2$ and $\xi_2:=(-e_1-e_2-e_3+e_4+e_5-e_6+\sqrt{2}e_7)/2$. Let $\tau:=(e_1+e_2+e_3+e_4+e_5+e_6+\sqrt{2}e_7)/2$, $\upsilon:=(e_1-e_2-e_3-e_4+e_5-e_6+\sqrt{2}e_7)/2$ and $\psi:=(-e_1+e_2-e_3-e_4-e_5-e_6-\sqrt{2}e_7)/2$.

\vspace{2mm}\noindent If $\mathcal{H}=\begin{tikzpicture}[scale=0.8,baseline=-2.2ex] \dynkin[ply=3,rotate=270,fold radius = 2mm]{D}{oooo} \end{tikzpicture} \hspace{-0.5mm}\times\hspace{-0.5mm} \begin{tikzpicture}[baseline=-0.5ex,rotate=90] 
\dynkin[name=upper]{A}{o}
\node (current) at ($(upper root 1)+(-0,3.4mm)$) {};
\dynkin[at=(current),name=middle]{A}{o}
\node (currenta) at ($(upper root 1)+(-0,7.4mm)$) {};
\dynkin[at=(currenta),name=lower]{A}{o}
\begin{scope}[on background layer]
\foreach \i in {1}%
{%
\draw[/Dynkin diagram/fold style]
($(upper root \i)$) -- ($(lower
root \i)$);%
}%
\end{scope}
\end{tikzpicture}$ then let $\Delta = \big\{\alpha_2,\alpha_3,\alpha_4,\alpha_5,\alpha_7,\sqrt{2}e_7,e_1+e_2\big\}$. Observe that the vectors $e_4-e_5$, $e_3+e_4$ and $2e_1+\sqrt{2}e_7$ are fixed by $\Pi$ and that $r_s(\mathcal{H}) =3$. So $\big\{e_4-e_5,e_3+e_4,2e_1+\sqrt{2}e_7\big\}$ is a basis for $E_s$. Hence $\Phi_a = \big\langle e_2-e_6,\xi_1,e_2+e_6,\xi_2 \big\rangle \cong (A_2)^2$.

\vspace{2mm}\noindent If $\mathcal{H}=\smash{\begin{tikzpicture}[baseline=-0.6ex]\dynkin[ply=2,fold radius = 2mm]{A}{*o*o*}\end{tikzpicture} \hspace{-0.5mm}\times\hspace{-0.5mm} \begin{tikzpicture}[baseline=-0.6ex]\dynkin[ply=2,fold radius = 2mm]{A}{**}\end{tikzpicture}}$ then let $\Delta^1 = \big\{ \alpha_2,\alpha_4,\alpha_5,\alpha_6,\alpha_7 \big\}$ and $\Delta^2 = \big\{\alpha_1,\tau\}$. Observe that $e_1+e_2-e_5-e_6$ is fixed by $\Pi$ and that $r_s(\mathcal{H}) =1$. So $E_s$ is generated by $e_1+e_2-e_5-e_6$. Hence $\Phi_a = \big\langle e_1-e_2,e_2+e_5,e_5-e_6,\alpha_1,e_3+e_4,e_3-e_4 \big\rangle \cong D_5A_1$.

\vspace{2mm}\noindent If $\mathcal{H}=\smash{\begin{tikzpicture}[baseline=-0.6ex]\dynkin[ply=2,fold radius = 2mm]{A}{*o*o*}\end{tikzpicture} \hspace{-0.5mm}\times\hspace{-0.5mm} \begin{tikzpicture}[baseline=-0.6ex]\dynkin[ply=2,fold radius = 2mm]{A}{oo}\end{tikzpicture}}$ then again let $\Delta^1 = \big\{ \alpha_2,\alpha_4,\alpha_5,\alpha_6,\alpha_7 \big\}$ and $\Delta^2 = \big\{\alpha_1,\tau \big\}$. Observe that $e_7$ and $e_1+e_2-e_5-e_6$ are both fixed by $\Pi$ and that $r_s(\mathcal{H}) =2$. So $\big\{e_7,e_1+e_2-e_5-e_6\big\}$ is a basis for $E_s$. Hence $\Phi_a = \big\langle e_1-e_2,e_2+e_5,e_5-e_6,e_3+e_4,e_3-e_4 \big\rangle \cong {\big(A_3(A_1)^2\big)}'$.

\vspace{2mm}\noindent If $\mathcal{H}=\smash{\begin{tikzpicture}[baseline=-0.6ex]\dynkin[ply=2,fold radius = 2mm]{A}{**o**}\end{tikzpicture} \hspace{-0.5mm}\times\hspace{-0.5mm} \begin{tikzpicture}[baseline=-0.6ex]\dynkin[ply=2,fold radius = 2mm]{A}{oo}\end{tikzpicture}}$ then again let $\Delta^1 = \big\{ \alpha_2,\alpha_4,\alpha_5,\alpha_6,\alpha_7 \big\}$ and $\Delta^2 = \big\{\alpha_1,\tau \big\}$. Observe that $e_7$ and $e_1+e_2+e_3-e_4-e_5-e_6$ are both fixed by $\Pi$ and that $r_s(\mathcal{H}) =2$. So $\big\{e_7,e_1+e_2+e_3-e_4-e_5-e_6\big\}$ is a basis for $E_s$. Hence $\Phi_a = \big\langle \alpha_2,\alpha_4,e_3+e_4,\alpha_6,\alpha_7 \big\rangle \cong (A_5)''$.

\vspace{2mm}\noindent If $\mathcal{H}=\smash{\begin{tikzpicture}[baseline=-0.6ex]\dynkin[ply=2,fold radius = 2mm]{A}{*o***o*}\end{tikzpicture}}$ then let $\Delta = \big\{ -\sqrt{2}e_7,\alpha_1,\alpha_3,\alpha_4,\alpha_5,\alpha_6,\alpha_7 \big\}$. Observe that $e_3+e_4+e_5-e_6$ is fixed by $\Pi$ and that $r_s(\mathcal{H}) =1$. So $E_s$ is generated by $e_3+e_4+e_5-e_6$. Hence $\Phi_a = \big\langle -\sqrt{2}e_7,\upsilon,\alpha_4,\alpha_3,\alpha_5,\alpha_7 \big\rangle \cong D_5A_1$.

\vspace{2mm}\noindent If $\mathcal{H}=\smash{\begin{tikzpicture}[baseline=-0.6ex]\dynkin[ply=2,fold radius=2mm]{E}{o****o}\end{tikzpicture} \hspace{-0.3mm}\times\hspace{-0.5mm} \mathcal{T}^1_0}$ then let $\Delta = \big\{ \alpha_1,\alpha_2,\alpha_3,\alpha_4,\alpha_5,\alpha_6 \big\}$. Observe that $-\sqrt{2}e_7+e_1-e_2$ is fixed by $\Pi$ and that $r_s(\mathcal{H}) =1$. So $E_s$ is generated by $-\sqrt{2}e_7+e_1-e_2$. Hence $\Phi_a = \big\langle \psi,e_5+e_6,e_4-e_5,e_5-e_6,e_3-e_4,e_1+e_2 \big\rangle \cong D_5A_1$.

\vspace{2mm}\noindent For all remaining indices $\mathcal{H}$ listed in the fourth column of Table \ref{E_7}, we observe that $\Psi_0$ is a parabolic subsystem of $\Phi$ and we use Lemma \ref{complemma} to check that $r(\Phi_a)=r(\Psi_0)$. So $\Phi_a=\Psi_0$ since $\Phi_a$ is a parabolic subsystem of $\Phi$ that contains $\Psi_0$. That is, $\mathcal{H}$ is independent in $\Phi$.

\newpage\newgeometry{top=1.3cm,bottom=1.3cm,left=0cm, right=0cm,foot=5mm}
\begin{longtable}{| c | c | c | c | c | c | c | c | c | c | c |}\caption{Isotropic maximal embeddings of abstract indices of type $E_7$}\label{E_7} \\ \hline
\multicolumn{1}{|c|}{} & \multicolumn{1}{c|}{} & \multicolumn{1}{c|}{} & \multicolumn{1}{c|}{} & \multicolumn{1}{c|}{} & \multicolumn{1}{c|}{} & \multicolumn{1}{c|}{} & \multicolumn{3}{c|}{\small{Special fields}} \\

  $\textcolor{white}{ai}\Delta\textcolor{white}{ia}$ & \small{$\hspace{-0.5mm}\Stab_W(\Delta)\hspace{-0.5mm}$} & $\Pi$ & $\textcolor{white}{iaaa} \mathcal{H} \textcolor{white}{aaai}$ & $\Psi_0$ & $\Phi_a$ & $\mathcal{G}$ & \small{$\hspace{-0.8mm}\!\cd 1 \!\hspace{-1mm}$} & $\!\R\!$ & $\!\!Q_{\mathfrak{p}}\!\!$  \\ \hline \hline 
\endfirsthead

\multicolumn{11}{|c|}%
{{\bfseries \tablename\ \thetable{} -- continued from previous page}} \\
\hline
\multicolumn{1}{|c|}{} & \multicolumn{1}{c|}{} & \multicolumn{1}{c|}{} & \multicolumn{1}{c|}{} & \multicolumn{1}{c|}{} & \multicolumn{1}{c|}{} & \multicolumn{1}{c|}{} & \multicolumn{3}{c|}{\small{Special fields}} \\

  $\Delta$ & \small{$\hspace{-0.5mm}\Stab_W(\Delta)\hspace{-0.5mm}$} & $\Pi$ & $\mathcal{H}$ & $\Psi_0$ & $\Phi_a$ & $\mathcal{G}$ & \small{$\hspace{-0.8mm}\!\cd 1 \!\hspace{-1mm}$} & $\!\R\!$ & $\!\!Q_{\mathfrak{p}}\!\!$ \\ \hline \hline \endhead
\hline \multicolumn{11}{|c|}{{Continued on next page}} \\ \hline \endfoot 
 \endlastfoot


\multirow{48}{*}{$\!A_5A_2\!$} & \multirow{48}{*}{$\Z_2$} & \multirow{11.25}{*}{$1$} & $\!\!\dynkin{A}{ooooo} \hspace{-0.5mm}\times\hspace{-0.5mm} \dynkin{A}{oo}\!\!$ & $\varnothing$ & $\varnothing$ & $\begin{tikzpicture}[baseline=0.4ex]\dynkin{E}{V}\end{tikzpicture}$ & $\checkmark$ & $\checkmark$ & $\checkmark$ \\ \cdashline{4-11} 
&&& $\!\!\dynkin{A}{ooooo} \hspace{-0.5mm}\times\hspace{-0.5mm} \dynkin{A}{**}\!\!$ & $A_2$ & $A_2$ &&&& \topstrut \\ \cdashline{4-11} 
&&& $\!\!\dynkin{A}{*o*o*} \hspace{-0.5mm}\times\hspace{-0.5mm} \dynkin{A}{oo}\!\!$ & $\big((A_1)^3\big)'$ & $\big((A_1)^3\big)'$ & $\begin{tikzpicture}[baseline=0.4ex]\dynkin{E}{o*oo*o*}\end{tikzpicture}$ & $\xmark$ & $\checkmark$ & $~\checkmark$ \topstrut \\ \cdashline{4-11} 
&&& $\!\!\dynkin{A}{*o*o*} \hspace{-0.5mm}\times\hspace{-0.5mm} \dynkin{A}{**}\!\!$ & $(A_1)^3A_2$ & $(A_1)^3A_2$ &&&& \topstrut \\ \cdashline{4-11} 
&&& $\!\!\dynkin{A}{**o**} \hspace{-0.5mm}\times\hspace{-0.5mm} \dynkin{A}{oo}\!\!$ & $(A_2)^2$ & $(A_2)^2$ &&&& \topstrut \\ \cdashline{4-11} 
&&& \multirow{3.6}{*}{$\!\!\dynkin{A}{**o**} \hspace{-0.5mm}\times\hspace{-0.5mm} \dynkin{A}{**}\!\!$} & \multirow{3.6}{*}{$(A_2)^3$} & \multirow{3.6}{*}{$E_6$} & $\begin{tikzpicture}[baseline=0.4ex]\dynkin{E}{V}\end{tikzpicture}$ & $\xmark$ & $\xmark$ & $?$ \\ 
&&&&&& $\begin{tikzpicture}[baseline=0.4ex]\dynkin{E}{o****oo}\end{tikzpicture}$ & $\xmark$ & $\xmark$ & $\xmark$ \\ 
&&&&&& $\begin{tikzpicture}[baseline=0.4ex]\dynkin{E}{******o}\end{tikzpicture}$ & $\xmark$ & $\xmark$ & $\xmark$ \\
\end{longtable}\restoregeometry

\subsection*{Appendix E. $\Phi=E_8 = \dynkin[scale=2,labels={\alpha_1,\alpha_...,\alpha_8}]{E}{8}$}\label{E8}

\noindent Let $\Phi=E_8$. We use the following \textit{standard} parametrisation in $\R^8$ of $\Phi$. Let $\alpha_1=-(\sum_{i=1}^8 e_i)/2$, $\alpha_2=e_6+e_7$, $\alpha_3=e_6-e_7$, $\alpha_4=e_5-e_6$, $\alpha_5=e_4-e_5$, $\alpha_6=e_3-e_4$, $\alpha_7=e_2-e_3$ and $\alpha_8=e_1-e_2$. The highest root of $\Phi$ is $e_1-e_8$. The Weyl group of $\Phi$ is $W \cong \Z_2 \hspace{0.5mm}.\hspace{0.5mm}\mathrm{PSO}_8^+ (2) \hspace{0.5mm}.\hspace{0.5mm} \Z_2$ and the longest element of $W$ is $-1$.

\vspace{1mm}\noindent By Table II of \cite{T}, the isotropic indices of type $\Phi$ consist of $\begin{tikzpicture}[baseline=0.3ex]\dynkin{E}{VIII}\end{tikzpicture}$, $\begin{tikzpicture}[baseline=0.3ex]\dynkin{E}{IX}\end{tikzpicture}$, $\begin{tikzpicture}[baseline=0.3ex]\dynkin{E}{o******o}\end{tikzpicture}$, $\begin{tikzpicture}[baseline=0.3ex]\dynkin{E}{******oo}\end{tikzpicture}$, $\begin{tikzpicture}[baseline=0.3ex]\dynkin{E}{o*******}\end{tikzpicture}$ and $\begin{tikzpicture}[baseline=0.3ex]\dynkin{E}{*******o}\end{tikzpicture}$. All of these indices are of inner type since $\Iso(\Phi)=W$.

\vspace{2mm}\noindent Recall from Table \ref{exceptlist} that $\Delta$ is one of $(A_4)^2$, $A_8$, $E_6A_2$, $D_8$, $(A_2)^4$, $E_7A_1$, $(D_4)^2$ or $(A_1)^8$. In the fourth column of Table \ref{E_8} we list all indices $\mathcal{H}$ of type $\Delta$ that are maximal in $\Phi$. For each of those indices $\mathcal{H}$, we compute the $W$-conjugacy class of $\Phi_a$ in $\Phi$.

\vspace{2mm}\noindent If $\mathcal{H}$ is one of $\begin{tikzpicture}[baseline=-0.6ex]\dynkin[ply=2,fold radius = 2mm]{A}{o**o}\end{tikzpicture} \times\hspace{-0.5mm} \begin{tikzpicture}[baseline=-0.6ex]\dynkin[ply=2,fold radius = 2mm]{A}{****}\end{tikzpicture}$, $\begin{tikzpicture}[baseline=-0.6ex]\dynkin[ply=2,fold radius = 2mm]{A}{o******o}\end{tikzpicture}$, $\begin{tikzpicture}[baseline=-0.6ex]\dynkin[ply=2,fold radius = 2mm]{E}{*o****}\end{tikzpicture} \hspace{-0.3mm}\times\hspace{-0.5mm} \begin{tikzpicture}[baseline=-0.6ex]\dynkin[ply=2,fold radius = 2mm]{A}{**}\end{tikzpicture}$, $\begin{tikzpicture}[baseline=0.4ex]\dynkin{E}{o******}\end{tikzpicture} \!\times\! \dynkin{A}{1}$, $\begin{tikzpicture}[scale=0.8,baseline=-2.2ex] \dynkin[ply=3,rotate=270,fold radius = 2mm]{D}{*o**} \end{tikzpicture} \!\times\! \begin{tikzpicture}[scale=0.8,baseline=-2.2ex] \dynkin[ply=3,rotate=270,fold radius = 2mm]{D}{****} \end{tikzpicture}$, $\begin{tikzpicture}[scale=.175,baseline=0.3ex]
    \foreach \y in {0,...,5}
    \draw[thin,xshift=\y cm] (\y cm,0) ++(.3 cm, 0) -- +(14 mm,0);
    \draw[thin,fill=black] (0 cm,0) circle (3 mm);
    \draw[thin,fill=black] (2 cm,0) circle (3 mm);
    \draw[thin,fill=black] (4 cm,0) circle (3 mm);
    \draw[thin,fill=black] (6 cm,0) circle (3 mm);
    \draw[thin,fill=black] (8 cm,0) circle (3 mm);
    \draw[thin] (10 cm,0) circle (3 mm);
    \draw[thin,fill=black] (12 cm,0) circle (3 mm);
    \draw[thin,fill=black] (2 cm,2 cm) circle (3 mm);
    \draw[thin,fill=black] (2 cm, 3mm) -- +(0, 1.4 cm);
  \end{tikzpicture}$ or $\begin{tikzpicture}[scale=.175,baseline=0.3ex]
    \foreach \y in {0,...,5}
    \draw[thin,xshift=\y cm] (\y cm,0) ++(.3 cm, 0) -- +(14 mm,0);
    \draw[thin] (0 cm,0) circle (3 mm);
    \draw[thin,fill=black] (2 cm,0) circle (3 mm);
    \draw[thin,fill=black] (4 cm,0) circle (3 mm);
    \draw[thin,fill=black] (6 cm,0) circle (3 mm);
    \draw[thin,fill=black] (8 cm,0) circle (3 mm);
    \draw[thin,fill=black] (10 cm,0) circle (3 mm);
    \draw[thin,fill=black] (12 cm,0) circle (3 mm);
    \draw[thin,fill=black] (2 cm,2 cm) circle (3 mm);
    \draw[thin,fill=black] (2 cm, 3mm) -- +(0, 1.4 cm);
  \end{tikzpicture}$ then $\Phi_s=(\Psi_0)^{\perp}$ (this last index is oriented such that $\Psi_0 \cong (A_7)'$, that is, $\Psi_0$ is a copy of $A_7$ in $E_8$ with a non-trivial perp). If $\mathcal{H}= \begin{tikzpicture}[baseline=-0.6ex]\dynkin[ply=2,fold radius = 2mm]{E}{******}\end{tikzpicture} \hspace{-0.3mm}\times\hspace{-0.5mm} \begin{tikzpicture}[baseline=-0.6ex]\dynkin[ply=2,fold radius = 2mm]{A}{oo}\end{tikzpicture}$ then $\Phi_s=\big\langle  \alpha_0^2 \big\rangle$. For each of these cases, $\Phi_s \cong A_1$ and $r(\Phi_s)=r_s(\mathcal{H}) =1$. Hence $\Phi_a =(\Phi_s)^{\perp} \cong E_7$ since $\Phi_s$ spans $E_s$.

\vspace{2mm}\noindent If $\mathcal{H}$ is one of $\begin{tikzpicture}[baseline=-0.6ex]\dynkin[ply=2,fold radius = 2mm]{E}{oo***o}\end{tikzpicture} \hspace{-0.3mm}\times\hspace{-0.5mm} \begin{tikzpicture}[baseline=-0.6ex]\dynkin[ply=2,fold radius = 2mm]{A}{**}\end{tikzpicture}$, $\begin{tikzpicture}[baseline=0.4ex]\dynkin{E}{o****o*}\end{tikzpicture} \!\times\! \dynkin{A}{1}$, $\begin{tikzpicture}[scale=0.8,baseline=-2.2ex] \dynkin[ply=3,rotate=270,fold radius = 2mm]{D}{*o**} \end{tikzpicture} \!\times\! \begin{tikzpicture}[scale=0.8,baseline=-2.2ex] \dynkin[ply=3,rotate=270,fold radius = 2mm]{D}{*o**} \end{tikzpicture}$, $\begin{tikzpicture}[scale=.175,baseline=0.3ex]
    \foreach \y in {0,...,5}
    \draw[thin,xshift=\y cm] (\y cm,0) ++(.3 cm, 0) -- +(14 mm,0);
    \draw[thin,fill=black] (0 cm,0) circle (3 mm);
    \draw[thin,fill=black] (2 cm,0) circle (3 mm);
    \draw[thin,fill=black] (4 cm,0) circle (3 mm);
    \draw[thin] (6 cm,0) circle (3 mm);
    \draw[thin,fill=black] (8 cm,0) circle (3 mm);
    \draw[thin] (10 cm,0) circle (3 mm);
    \draw[thin,fill=black] (12 cm,0) circle (3 mm);
    \draw[thin,fill=black] (2 cm,2 cm) circle (3 mm);
    \draw[thin,fill=black] (2 cm, 3mm) -- +(0, 1.4 cm);
  \end{tikzpicture}$ or $\begin{tikzpicture}[scale=.175,baseline=0.3ex]
    \foreach \y in {0,...,5}
    \draw[thin,xshift=\y cm] (\y cm,0) ++(.3 cm, 0) -- +(14 mm,0);
    \draw[thin] (0 cm,0) circle (3 mm);
    \draw[thin,fill=black] (2 cm,0) circle (3 mm);
    \draw[thin,fill=black] (4 cm,0) circle (3 mm);
    \draw[thin] (6 cm,0) circle (3 mm);
    \draw[thin,fill=black] (8 cm,0) circle (3 mm);
    \draw[thin,fill=black] (10 cm,0) circle (3 mm);
    \draw[thin,fill=black] (12 cm,0) circle (3 mm);
    \draw[thin,fill=black] (2 cm,2 cm) circle (3 mm);
    \draw[thin,fill=black] (2 cm, 3mm) -- +(0, 1.4 cm);
  \end{tikzpicture}$ then $\Phi_s=(\Psi_0)^{\perp}$ (this last index is oriented such that $\Psi_0 \cong \big((A_3)^2\big)'$, that is, $\Psi_0$ is a copy of $(A_3)^2$ in $E_8$ with a non-trivial perp). If $\mathcal{H}=\begin{tikzpicture}[baseline=-0.6ex]\dynkin[ply=2,fold radius = 2mm]{A}{o**o}\end{tikzpicture} \times\hspace{-0.5mm} \begin{tikzpicture}[baseline=-0.6ex]\dynkin[ply=2,fold radius = 2mm]{A}{o**o}\end{tikzpicture}$ then $\Phi_s=(\Psi_0)^{\perp_{\Psi}}$. If $\mathcal{H}=\begin{tikzpicture}[baseline=-0.6ex]\dynkin[ply=2,fold radius = 2mm]{E}{*o****}\end{tikzpicture} \hspace{-0.3mm}\times\hspace{-0.5mm} \begin{tikzpicture}[baseline=-0.6ex]\dynkin[ply=2,fold radius = 2mm]{A}{oo}\end{tikzpicture}$ then $\Phi_s = \big\langle \alpha_0^1, \alpha_0^2\big\rangle$. If $\mathcal{H}=\begin{tikzpicture}[baseline=-0.6ex]\dynkin[ply=2,fold radius = 2mm]{A}{oooo}\end{tikzpicture} \times\hspace{-0.5mm} \begin{tikzpicture}[baseline=-0.6ex]\dynkin[ply=2,fold radius = 2mm]{A}{****}\end{tikzpicture}$ then $\Phi_s= \big\langle \alpha_0^1, \beta\big\rangle$, where $\beta = 0~1~1~0 \in \Psi^1$. If $\mathcal{H}=\begin{tikzpicture}[baseline=-0.6ex]\dynkin[ply=2,fold radius = 2mm]{A}{oo****oo}\end{tikzpicture}$ then $\Phi_s = \big\langle \alpha_0, \gamma \big\rangle$, where $\gamma=0~1~1~1~1~1~1~0 \in \Psi$. For each of these cases, $\Phi_s \cong (A_1)^2$ and $r(\Phi_s)=r_s(\mathcal{H}) =2$. Hence $\Phi_a =(\Phi_s)^{\perp} \cong D_6$.

\vspace{2mm}\noindent If $\mathcal{H}$ is one of $\begin{tikzpicture}[baseline=-0.6ex]\dynkin[ply=2,fold radius = 2mm]{E}{*o*o**}\end{tikzpicture} \hspace{-0.3mm}\times\hspace{-0.5mm} \begin{tikzpicture}[baseline=-0.6ex]\dynkin[ply=2,fold radius = 2mm]{A}{**}\end{tikzpicture}$, $\dynkin{A}{**o**o**}$ or $\begin{tikzpicture}[baseline=0.4ex]\dynkin{E}{*o*o**}\end{tikzpicture} \hspace{-0.5mm}\times\hspace{-0.5mm} \dynkin{A}{**}$ then $\Phi_s=(\Psi_0)^{\perp}$. If $\mathcal{H}=\begin{tikzpicture}[scale=0.8,baseline=-2.2ex] \dynkin[ply=3,rotate=270,fold radius = 2mm]{D}{oooo} \end{tikzpicture} \!\times\! \begin{tikzpicture}[scale=0.8,baseline=-2.2ex] \dynkin[ply=3,rotate=270,fold radius = 2mm]{D}{****} \end{tikzpicture}$ then $\Phi_s = \big\langle \alpha_0^1,\delta \big\rangle$, where $\delta=\smash{\small\makecell{$\small0$ \\ $\small0~1~0$}} \in \Psi^1$. \textcolor{white}{$\!\!\!\Big($}For each of these cases, $\Phi_s \cong A_2$ and $r(\Phi_s)=r_s(\mathcal{H}) =2$. Hence $\Phi_a =(\Phi_s)^{\perp} \cong E_6$.

\vspace{2mm}\noindent If $\mathcal{H}=\begin{tikzpicture}[baseline=-0.6ex]\dynkin[ply=2,fold radius = 2mm]{A}{oooo}\end{tikzpicture} \times\hspace{-0.5mm} \begin{tikzpicture}[baseline=-0.6ex]\dynkin[ply=2,fold radius = 2mm]{A}{o**o}\end{tikzpicture}$ then $\Phi_s = \big\langle \alpha_0^1,\alpha_0^2, \epsilon \big\rangle$, where $\epsilon =0~1~1~0 \in \Psi^1$. If $\mathcal{H}=\begin{tikzpicture}[baseline=-0.6ex]\dynkin[ply=2,fold radius = 2mm]{A}{ooo**ooo}\end{tikzpicture}$ then $\Phi_s = \big\langle \alpha_0,\zeta_1,\zeta_2 \big\rangle$, where $\zeta_1 =0~1~1~1~1~1~1~0 \in \Psi$ and $\zeta_2=0~0~1~1~1~1~0~0 \in \Psi$. If $\mathcal{H}=\begin{tikzpicture}[baseline=-0.6ex]\dynkin[ply=2,fold radius = 2mm]{E}{oo***o}\end{tikzpicture} \hspace{-0.3mm}\times\hspace{-0.5mm} \begin{tikzpicture}[baseline=-0.6ex]\dynkin[ply=2,fold radius = 2mm]{A}{oo}\end{tikzpicture}$ then $\Phi_s$ is generated by $\alpha_0^2$ and the perp of $\Psi_0$ in $\Psi^1$. For each of these cases, $\Phi_s \cong (A_1)^3$ and $r(\Phi_s)=r_s(\mathcal{H}) =3$. Hence $\Phi_a =(\Phi_s)^{\perp} \cong D_4A_1$.

\vspace{2mm}\noindent If $\mathcal{H}=\smash{\begin{tikzpicture}[baseline=-0.6ex]\dynkin[ply=2,fold radius = 2mm]{E}{*o*o**}\end{tikzpicture} \hspace{-0.3mm}\times\hspace{-0.5mm} \begin{tikzpicture}[baseline=-0.6ex]\dynkin[ply=2,fold radius = 2mm]{A}{oo}\end{tikzpicture}}$ then $\Phi_s$ is generated by $\alpha_0^2$ and the perp of $\Psi_0$ in $\Psi^1$. If $\mathcal{H}=\begin{tikzpicture}[scale=0.8,baseline=-2.2ex] \dynkin[ply=3,rotate=270,fold radius = 2mm]{D}{oooo} \end{tikzpicture} \!\times\! \begin{tikzpicture}[scale=0.8,baseline=-2.2ex] \dynkin[ply=3,rotate=270,fold radius = 2mm]{D}{*o**} \end{tikzpicture}$ then $\Phi_s =\big\langle \alpha_0^1, \alpha_0^2, \eta \big\rangle$, where $\eta=\smash{\small\makecell{$\small0$ \\ $\small0~1~0$}} \in \Psi^1$. \textcolor{white}{$\!\!\!\Big($}For both of these cases, $\Phi_s \cong A_2A_1$ and $r(\Phi_s)=r_s(\mathcal{H}) =3$. Hence $\Phi_a =(\Phi_s)^{\perp} \cong A_5$.

\vspace{3mm}\noindent If $\mathcal{H}=\smash{\begin{tikzpicture}[baseline=-0.6ex]\dynkin[ply=2,fold radius = 2mm]{A}{oooo}\end{tikzpicture} \times\hspace{-0.5mm} \begin{tikzpicture}[baseline=-0.6ex]\dynkin[ply=2,fold radius = 2mm]{A}{oooo}\end{tikzpicture}}$ then $\Phi_s = \big\langle \alpha_0^1,\alpha_0^2, \kappa_1, \kappa_2 \big\rangle$, where $\kappa_i =0~1~1~0 \in \Psi^i$ for $i=1,2$. If $\mathcal{H}=\begin{tikzpicture}[baseline=-0.6ex]\dynkin[ply=2,fold radius = 2mm]{A}{oooooooo}\end{tikzpicture}$ then $\Phi_s = \big\langle \alpha_0,\lambda_1,\lambda_2,\lambda_3 \big\rangle$, where $\lambda_1=0~1~1~1~1~1~1~0 \in \Psi$, $\lambda_2=0~0~1~1~1~1~0~0 \in \Psi$ and $\lambda_3=0~0~0~1~1~0~0~0 \in \Psi$. For both of these cases, $\Phi_s \cong \smash{\big((A_1)^4\big)''}$ and $r(\Phi_s)=r_s(\mathcal{H}) =4$. That is, $\Phi_s$ is a copy of $(A_1)^4$ in $E_8$ that is parabolic and is not contained in some $D_4$ subsystem. Hence $\Phi_a =(\Phi_s)^{\perp} \cong \smash{\big((A_1)^4\big)''}$.

\vspace{2mm}\noindent If $\mathcal{H}=\begin{tikzpicture}[scale=0.8,baseline=-2.2ex] \dynkin[ply=3,rotate=270,fold radius = 2mm]{D}{oooo} \end{tikzpicture} \!\times\! \begin{tikzpicture}[scale=0.8,baseline=-2.2ex] \dynkin[ply=3,rotate=270,fold radius = 2mm]{D}{oooo} \end{tikzpicture}$ then $\Phi_s =\big\langle \alpha_0^1, \alpha_0^2, \mu_1, \mu_2 \big\rangle$, where $\mu_i=\smash{\small\makecell{$\small0$ \\ $\small0~1~0$}} \in \Psi^i$ for $i=1,2$. Observe that $\Phi_s \cong (A_2)^2$ and $r(\Phi_s)=r_s(\mathcal{H}) =4$. Hence $\Phi_a =(\Phi_s)^{\perp} \cong (A_2)^2$.

\vspace{2mm}\noindent If $\mathcal{H}$ is either $\begin{tikzpicture}[baseline=0.4ex]\dynkin{E}{VI}\end{tikzpicture} \!\times\! \dynkin{A}{1}$ or $\begin{tikzpicture}[scale=.175,baseline=0.3ex]
    \foreach \y in {0,...,5}
    \draw[thin,xshift=\y cm] (\y cm,0) ++(.3 cm, 0) -- +(14 mm,0);
    \draw[thin] (0 cm,0) circle (3 mm);
    \draw[thin] (2 cm,0) circle (3 mm);
    \draw[thin,fill=black] (4 cm,0) circle (3 mm);
    \draw[thin] (6 cm,0) circle (3 mm);
    \draw[thin,fill=black] (8 cm,0) circle (3 mm);
    \draw[thin] (10 cm,0) circle (3 mm);
    \draw[thin,fill=black] (12 cm,0) circle (3 mm);
    \draw[thin,fill=black] (2 cm,2 cm) circle (3 mm);
    \draw[thin,fill=black] (2 cm, 3mm) -- +(0, 1.4 cm);
  \end{tikzpicture}$ then $\Phi_s=(\Psi_0)^{\perp}$ (note that $\Psi_0 \cong \smash{\big((A_1)^4\big)'}$, that is, $\Psi_0$ is a copy of $(A_1)^4$ in $E_8$ that is not parabolic). If $\mathcal{H}=\smash{\begin{tikzpicture}[baseline=-0.6ex]\dynkin[ply=2,fold radius = 2mm]{E}{oooooo}\end{tikzpicture} \hspace{-0.3mm}\times\hspace{-0.5mm} \begin{tikzpicture}[baseline=-0.6ex]\dynkin[ply=2,fold radius = 2mm]{A}{**}\end{tikzpicture}}$ then $\Phi_s =\big\langle \alpha_0^1, \nu_1,\nu_2,\nu_3 \big\rangle$, where \textcolor{white}{$\!\!\!\Big($}$\nu_1=\smash{\small\makecell{$\small0$ \\ $\small0~0~1~0~0$}} \in \Psi^1$, $\nu_2=\smash{\small\makecell{$\small1$ \\ $\small0~0~1~0~0$}} \in \Psi^1$ and $\nu_3=\smash{\small\makecell{$\small0$ \\ $\small1~1~1~1~1$}} \in \Psi^1$. For each of these cases, $\Phi_s \cong D_4$ and $r(\Phi_s)=r_s(\mathcal{H}) =4$. Hence $\Phi_a =(\Phi_s)^{\perp} \cong D_4$.

\vspace{2mm}\noindent If \textcolor{white}{$\!\!\!\Big($}$\mathcal{H}=\begin{tikzpicture}[baseline=-0.6ex]\dynkin[ply=2,fold radius = 2mm]{E}{oooooo}\end{tikzpicture} \hspace{-0.3mm}\times\hspace{-0.5mm} \begin{tikzpicture}[baseline=-0.6ex]\dynkin[ply=2,fold radius = 2mm]{A}{oo}\end{tikzpicture}$ then let $\xi_1=\smash{\small\makecell{$\small0$ \\ $\small0~0~1~0~0$}} \in \Psi^1$, $\xi_2=\smash{\small\makecell{$\small1$ \\ $\small0~0~1~0~0$}} \in \Psi^1$ and $\xi_3=\smash{\small\makecell{$\small0$ \\ $\small1~1~1~1~1$}} \in \Psi^1$. Then $\Phi_s =\big\langle \alpha_0^1, \alpha_0^2, \xi_1,\xi_2,\xi_3 \big\rangle \cong D_4A_1$ and $r(\Phi_s)=r_s(\mathcal{H}) =5$. Hence $\Phi_a =(\Phi_s)^{\perp} \cong (A_1)^3$.

 \vspace{2mm}\noindent If $\mathcal{H}=\smash{\begin{tikzpicture}[baseline=0.4ex]\dynkin{E}{*****o*}\end{tikzpicture}} \!\times\! \dynkin{A}{1}$ then $\Phi_a \cong D_7$ as this is the only parabolic subsystem of $\Phi$ that contains $\Psi_0 \cong D_5(A_1)^2$. Similarly, if $\mathcal{H}=\begin{tikzpicture}[scale=.175,baseline=0.3ex]
    \foreach \y in {0,...,5}
    \draw[thin,xshift=\y cm] (\y cm,0) ++(.3 cm, 0) -- +(14 mm,0);
    \draw[thin,fill=black] (0 cm,0) circle (3 mm);
    \draw[thin,fill=black] (2 cm,0) circle (3 mm);
    \draw[thin,fill=black] (4 cm,0) circle (3 mm);
    \draw[thin] (6 cm,0) circle (3 mm);
    \draw[thin,fill=black] (8 cm,0) circle (3 mm);
    \draw[thin,fill=black] (10 cm,0) circle (3 mm);
    \draw[thin,fill=black] (12 cm,0) circle (3 mm);
    \draw[thin,fill=black] (2 cm,2 cm) circle (3 mm);
    \draw[thin,fill=black] (2 cm, 3mm) -- +(0, 1.4 cm);
  \end{tikzpicture}~$ then $\Phi_a \cong D_7$ as this is the only parabolic subsystem of $\Phi$ that contains $\Psi_0 \cong D_4A_3$. 

\vspace{2mm}\noindent Consider the case where $\mathcal{H}=\smash{\begin{tikzpicture}[baseline=0.5ex,rotate=90] 
\dynkin[name=upper]{A}{oo}
\node (current) at ($(upper root 1)+(-0,3.4mm)$) {};
\dynkin[at=(current),name=middle]{A}{oo}
\node (currenta) at ($(upper root 1)+(-0,7.4mm)$) {};
\dynkin[at=(currenta),name=lower]{A}{oo}
\node (currentaa) at ($(upper root 1)+(-0,11.4mm)$) {};
\dynkin[at=(currentaa),name=muchlower]{A}{oo}
\begin{scope}[on background layer]
\draw[/Dynkin diagram/fold style]
($(upper root 1)$) -- ($(muchlower root 2)$);
\draw[/Dynkin diagram/fold style]
($(upper root 2)$) -- ($(muchlower root 1)$);
\draw[/Dynkin diagram/fold style]
($(upper root 1)$) -- ($(muchlower root 1)$);
\draw[/Dynkin diagram/fold style]
($(upper root 2)$) -- ($(muchlower root 2)$); 
\end{scope}
 \end{tikzpicture}}$ (resp. $\smash{\begin{tikzpicture}[scale=1.05,baseline=0.5ex,rotate=90] 
\dynkin[name=upper]{A}{oo}
\node (current) at ($(upper root 1)+(-0,3.4mm)$) {};
\dynkin[at=(current),name=middle]{A}{oo}
\begin{scope}[on background layer]
\draw[/Dynkin diagram/fold style]
($(upper root 1)$) -- ($(middle
root 1)$);
\draw[/Dynkin diagram/fold style]
($(upper root 2)$) -- ($(middle
root 2)$);
\draw[/Dynkin diagram/fold style]
($(upper root 1)$) -- ($(middle
root 2)$);
\draw[/Dynkin diagram/fold style]
($(upper root 2)$) -- ($(middle
root 1)$);
\end{scope}
\end{tikzpicture} \times \begin{tikzpicture}[scale=1.05,baseline=0.5ex,rotate=90] 
\dynkin[name=upper]{A}{oo}
\node (current) at ($(upper root 1)+(-0,3.4mm)$) {};
\dynkin[at=(current),name=middle]{A}{oo}
\begin{scope}[on background layer]
\draw[/Dynkin diagram/fold style]
($(upper root 1)$) -- ($(middle
root 1)$);
\draw[/Dynkin diagram/fold style]
($(upper root 2)$) -- ($(middle
root 2)$);
\draw[/Dynkin diagram/fold style]
($(upper root 1)$) -- ($(middle
root 2)$);
\draw[/Dynkin diagram/fold style]
($(upper root 2)$) -- ($(middle
root 1)$);
\end{scope}
\end{tikzpicture}}$). Let $\Theta :=\big\langle \alpha_0^1,\alpha_0^2,\alpha_0^3,\alpha_0^4\big\rangle$. Observe that $\Theta \cong \Theta^{\perp} \cong \smash{\big((A_1)^4\big)''}$ (that is, $\Theta$ and $\Theta^{\perp}$ are parabolic subsystems of $E_8$ that are not contained in some $D_4$ subsystem). Observe that $\Pi$ acts transitively (resp. with two orbits of equal size) on the factors of $\Theta$. By inspection of the action of $\Aut_W(W_{\Theta \Theta^{\perp}}) \cong \AGL_3(2)$ on the factors of $\Theta \Theta^{\perp}$, we deduce that $\Pi$ also acts transitively (resp. with two orbits of equal size) on the factors of $\Theta^{\perp}$. Let $\Theta_{\R}$ denote the subspace of $E$ that is generated by $\Theta$. Note that $E_s$ is a subspace of $\Theta_{\R}$ of dimension $1$ (resp. $2$). Hence, by Lemma \ref{complemma}, $\Phi_a$ is a rank $7$ (resp. $6$) parabolic subsystem of $\Phi$ that contains $\Theta^{\perp}$ and $\Aut_{W_{\Phi_a}}(W_{\Theta^{\perp}})$ acts transitively (resp. with two orbits of equal size) on the factors of $\Theta^{\perp}$. The only possibility for such a subsystem $\Phi_a$ is of type $(A_7)''$ (resp. $\smash{\big((A_3)^2\big)''}$).

\vspace{2mm}\noindent If $\mathcal{H}=\smash{\begin{tikzpicture}[baseline=-0.6ex]\dynkin[ply=2,fold radius = 2mm]{A}{**o**o**}\end{tikzpicture}}$ then $\Pi=\langle \rho \rangle \cong \Z_2$, where $\rho$ acts on $\Delta$ as the non-trivial graph automorphism. Let $\Omega:=\smash{(\Psi_0)^{\perp}} \cong A_2$. Observe that $\rho$ acts on $\Omega$ as an element of order $2$ in $W_{\Omega}$. In particular, $\rho$ does not fix pointwise any root in $\Omega$. By Lemma \ref{complemma}, $\Phi_a$ is a rank $7$ parabolic subsystem of $\Phi$ that contains $\Psi_0 \cong (A_2)^3$. There are two possibilities, either $\Phi_a$ is of type $E_7$ or of type $E_6A_1$. Assume that $\Phi_a \cong E_7$. Then $\Pi$ acts trivially on $(\Phi_a)^{\perp} \cong A_1$ by Lemma \ref{complemma}. But $(\Phi_a)^{\perp} \subset \Omega$ and so $\rho$ fixes some root in $\Omega$. This is a contradiction. Hence $\Phi_a \cong E_6A_1$.

\vspace{2mm}\noindent For the next cases, without loss of generality, we choose a representative of the $W$-conjugacy class of $\Delta$ in $\Phi$ using the standard parametrisation in $\R^8$ of $\Phi$. We will need the following roots in $\Phi$. Let $\omega_1:=(-e_1+e_2+e_3+e_4+e_5+e_6-e_7+e_8)/2$, $\omega_2:=(-e_1-e_2+e_3-e_4-e_5-e_6-e_7+e_8)/2$, $\omega_3:=(-e_1-e_2-e_3-e_4+e_5+e_6+e_7+e_8)/2$ and $\omega_4:=(e_1+e_2+e_3-e_4-e_5+e_6+e_7+e_8)/2$. Let $\tau_1:=(e_1+e_2+e_3+e_4+e_5+e_6-e_7-e_8)/2$, $\tau_2:=(e_1+e_2+e_3-e_4-e_5-e_6+e_7-e_8)/2$ and $\tau_3:=(-e_1-e_2-e_3+e_4+e_5+e_6+e_7-e_8)/2$. Let $\upsilon:=(e_1+e_2+e_3-e_4-e_5+e_6-e_7-e_8)/2$, $\phi:=(-e_1+e_2-e_3+e_4-e_5+e_6+e_7-e_8)/2$ and $\psi:=(-e_1+e_2+e_3-e_4+e_5+e_6+e_7+e_8)/2$.

\vspace{2mm}\noindent If $\mathcal{H}=\smash{\begin{tikzpicture}[scale=1,baseline=0.3ex] 
\dynkin[name=upper]{A}{o**o}
\node (current) at ($(upper root 1)+(-0,3.1mm)$) {};
\dynkin[at=(current),name=middle]{A}{o**o}
\begin{scope}[on background layer]
\draw[/Dynkin diagram/fold style]
($(upper root 1)$) -- ($(middle
root 1)$);
\draw[/Dynkin diagram/fold style]
($(upper root 2)$) -- ($(middle
root 2)$);
\draw[/Dynkin diagram/fold style]
($(upper root 3)$) -- ($(middle
root 3)$);
\draw[/Dynkin diagram/fold style]
($(upper root 4)$) -- ($(middle
root 4)$);
\draw[/Dynkin diagram/fold style]
($(upper root 1)$) -- ($(middle
root 4)$);
\draw[/Dynkin diagram/fold style]
($(upper root 2)$) -- ($(middle
root 3)$);
\draw[/Dynkin diagram/fold style]
($(upper root 3)$) -- ($(middle
root 2)$);
\draw[/Dynkin diagram/fold style]
($(upper root 4)$) -- ($(middle
root 1)$);
\end{scope}
\end{tikzpicture}}$ then let $\Delta^1 = \big\{\omega_1,\alpha_8,\alpha_7,\omega_2\big\}$ and $\Delta^2 = \big\{ \omega_3,\alpha_5,\alpha_4, \omega_4\big\}$. If $\mathcal{H}=\smash{\begin{tikzpicture}[scale=1.05,baseline=0.5ex,rotate=90] 
\dynkin[name=upper]{A}{oo}
\node (current) at ($(upper root 1)+(-0,3.4mm)$) {};
\dynkin[at=(current),name=middle]{A}{oo}
\begin{scope}[on background layer]
\draw[/Dynkin diagram/fold style]
($(upper root 1)$) -- ($(middle
root 1)$);
\draw[/Dynkin diagram/fold style]
($(upper root 2)$) -- ($(middle
root 2)$);
\draw[/Dynkin diagram/fold style]
($(upper root 1)$) -- ($(middle
root 2)$);
\draw[/Dynkin diagram/fold style]
($(upper root 2)$) -- ($(middle
root 1)$);
\end{scope}
\end{tikzpicture} \times \begin{tikzpicture}[scale=1.05,baseline=0.5ex,rotate=90] 
\dynkin[name=upper]{A}{**}
\node (current) at ($(upper root 1)+(-0,3.4mm)$) {};
\dynkin[at=(current),name=middle]{A}{**}
\begin{scope}[on background layer]
\draw[/Dynkin diagram/fold style]
($(upper root 1)$) -- ($(middle
root 1)$);
\draw[/Dynkin diagram/fold style]
($(upper root 2)$) -- ($(middle
root 2)$);
\draw[/Dynkin diagram/fold style]
($(upper root 1)$) -- ($(middle
root 2)$);
\draw[/Dynkin diagram/fold style]
($(upper root 2)$) -- ($(middle
root 1)$);
\end{scope}
\end{tikzpicture}}$ then let $\Delta = \big\{\alpha_1,\tau_1,\tau_2,\tau_3, \alpha_4, \alpha_5,\alpha_7,\alpha_8 \big\}$ where $\Delta_0 = \big\{ \alpha_4, \alpha_5,\alpha_7,\alpha_8  \big\}$. If $\mathcal{H}=\begin{tikzpicture}[baseline=-0.6ex]\dynkin[ply=2,fold radius = 2mm]{E}{o****o}\end{tikzpicture} \smash{\times} \begin{tikzpicture}[baseline=-0.6ex]\dynkin[ply=2,fold radius = 2mm]{A}{**}\end{tikzpicture}$ then let $\Delta^1 = \big\{ \alpha_1,\alpha_2,\alpha_3,\alpha_4,\alpha_5,\upsilon\big\}$ and $\Delta^2 = \big\{\alpha_7,\alpha_8 \big\}$. If $\mathcal{H}=\smash{\begin{tikzpicture}[scale=1,baseline=-0.3ex] 
\dynkin[name=upper]{D}{*o**}
\node (current) at ($(upper root 1)+(4.25mm,-3.45mm)$) {};
\dynkin[at=(current),name=middle,rotate=60]{D}{*o**}
\begin{scope}[on background layer]
\draw[/Dynkin diagram/fold style]
($(upper root 1)$) -- ($(middle
root 1)$);
\draw[/Dynkin diagram/fold style]
($(upper root 2)$) -- ($(middle
root 2)$);
\draw[/Dynkin diagram/fold style]
($(upper root 3)$) -- ($(middle
root 3)$);
\draw[/Dynkin diagram/fold style]
($(upper root 4)$) -- ($(middle
root 4)$);
\draw[/Dynkin diagram/fold style]
($(upper root 1)$) -- ($(middle
root 3)$);
\draw[/Dynkin diagram/fold style]
($(upper root 3)$) -- ($(middle
root 4)$);
\draw[/Dynkin diagram/fold style]
($(upper root 4)$) -- ($(middle
root 1)$);
\end{scope}
\end{tikzpicture}}$ then let $\Delta^1=\big\{\alpha_8,\alpha_6,\alpha_4,\phi \big\}$ and $\Delta^2=\big\{e_1+e_2,e_3+e_4,e_5+e_6,\alpha_1 \big\}$. For each of these cases, observe that $e_8$ is fixed by $\Pi$ and $r_s(\mathcal{H}) =1$. So $E_s$ is generated by $e_8$. Hence $\Phi_a = \big\langle \alpha_2, \alpha_3, \alpha_4, \alpha_5, \alpha_6, \alpha_7, \alpha_8 \big\rangle \cong D_7$.

\vspace{2mm}\noindent If $\mathcal{H}=\begin{tikzpicture}[baseline=-0.6ex,rotate=90] 
\dynkin[name=upper]{A}{o}
\node (current) at ($(upper root 1)+(-0,3.1mm)$) {};
\dynkin[at=(current),name=first]{A}{o}
\node (currenta) at ($(upper root 1)+(-0,6.9mm)$) {};
\dynkin[at=(currenta),name=second]{A}{o}
\node (currentaa) at ($(upper root 1)+(-0,10.7mm)$) {};
\dynkin[at=(currentaa),name=third]{A}{o}
\node (currentaaa) at ($(upper root 1)+(-0,14.5mm)$) {};
\dynkin[at=(currentaaa),name=fourth]{A}{o}
\node (currentaaaa) at ($(upper root 1)+(-0,18.3mm)$) {};
\dynkin[at=(currentaaaa),name=fifth]{A}{o}
\node (currentaaaaa) at ($(upper root 1)+(-0,22.1mm)$) {};
\dynkin[at=(currentaaaaa),name=sixth]{A}{o}
\node (currentaaaaaa) at ($(upper root 1)+(-0,25.9mm)$) {};
\dynkin[at=(currentaaaaaa),name=lower]{A}{o}
\begin{scope}[on background layer]
\foreach \i in {1}%
{%
\draw[/Dynkin diagram/fold style]
($(upper root \i)$) -- ($(lower
root \i)$);%
}%
\end{scope}
\end{tikzpicture}$ then let $\Delta = \big\{ e_1-e_8,-e_1-e_8,e_2-e_3,e_2+e_3,e_4-e_5,e_4+e_5,e_6-e_7,e_6+e_7 \big\}$. Observe that $e_2+e_4+e_6-e_8$ is fixed by $\Pi$ and that $r_s(\mathcal{H}) =1$. So $E_s$ is generated by $e_2+e_4+e_6-e_8$. Hence $\Phi_a = \big\langle -e_4-e_8,e_4-e_2,\psi, e_1-e_3, e_3-e_5,e_5-e_7,-e_1-e_3 \big\rangle \cong D_7$.

\vspace{2mm}\noindent If $\mathcal{H}=\smash{\begin{tikzpicture}[scale=1,baseline=0.3ex] 
\dynkin[name=upper]{A}{oooo}
\node (current) at ($(upper root 1)+(-0,3.1mm)$) {};
\dynkin[at=(current),name=middle]{A}{oooo}
\begin{scope}[on background layer]
\draw[/Dynkin diagram/fold style]
($(upper root 1)$) -- ($(middle
root 1)$);
\draw[/Dynkin diagram/fold style]
($(upper root 2)$) -- ($(middle
root 2)$);
\draw[/Dynkin diagram/fold style]
($(upper root 3)$) -- ($(middle
root 3)$);
\draw[/Dynkin diagram/fold style]
($(upper root 4)$) -- ($(middle
root 4)$);
\draw[/Dynkin diagram/fold style]
($(upper root 1)$) -- ($(middle
root 4)$);
\draw[/Dynkin diagram/fold style]
($(upper root 2)$) -- ($(middle
root 3)$);
\draw[/Dynkin diagram/fold style]
($(upper root 3)$) -- ($(middle
root 2)$);
\draw[/Dynkin diagram/fold style]
($(upper root 4)$) -- ($(middle
root 1)$);
\end{scope}
\end{tikzpicture}}$ then again let $\Delta^1 = \big\{\omega_1,\alpha_8,\alpha_7,\omega_2\big\}$ and $\Delta^2 = \big\{ \omega_3,\alpha_5,\alpha_4, \omega_4\big\}$. Observe that $e_8$ and $e_1-e_3+e_4-e_6$ are both fixed by $\Pi$ and that $r_s(\mathcal{H}) =2$. So $\big\{e_8,e_1-e_3+e_4-e_6\big\}$ is a basis for $E_s$. Hence $\Phi_a = \big\langle e_1+e_3, -e_3-e_4, e_4+e_6 \big\rangle \cup \big\langle e_2-e_5,e_5-e_7,-e_2-e_5 \big\rangle \cong \smash{\big((A_3)^2\big)''}$.

\vspace{2mm}\noindent If $\mathcal{H}=\smash{\begin{tikzpicture}[baseline=-0.6ex]\dynkin[ply=2,fold radius = 2mm]{E}{o****o}\end{tikzpicture} \hspace{-0.3mm}\times\hspace{-0.5mm} \begin{tikzpicture}[baseline=-0.6ex]\dynkin[ply=2,fold radius = 2mm]{A}{oo}\end{tikzpicture}}$ then let $\Delta^1 = \big\{ \alpha_1,\alpha_2,\alpha_3,\alpha_4,\alpha_5,\upsilon\big\}$ and $\Delta^2 = \big\{\alpha_7,\alpha_8 \big\}$. Observe that $e_8$ and $e_1-e_3$ are both fixed by $\Pi$ and that $r_s(\mathcal{H}) =2$. So $\big\{e_8,e_1-e_3\big\}$ is a basis for $E_s$. Hence $\Phi_a = \big\langle \alpha_2, \alpha_3, \alpha_4, \alpha_5, e_2-e_4,e_1+e_3 \big\rangle \cong D_5A_1$.

\vspace{2mm}\noindent If $\mathcal{H}=\begin{tikzpicture}[scale=1,baseline=-0.3ex] 
\dynkin[name=upper]{D}{oooo}
\node (current) at ($(upper root 1)+(4.2mm,-3.45mm)$) {};
\dynkin[at=(current),name=middle,rotate=60]{D}{oooo}
\begin{scope}[on background layer]
\draw[/Dynkin diagram/fold style]
($(upper root 1)$) -- ($(middle
root 1)$);
\draw[/Dynkin diagram/fold style]
($(upper root 2)$) -- ($(middle
root 2)$);
\draw[/Dynkin diagram/fold style]
($(upper root 3)$) -- ($(middle
root 3)$);
\draw[/Dynkin diagram/fold style]
($(upper root 4)$) -- ($(middle
root 4)$);
\draw[/Dynkin diagram/fold style]
($(upper root 1)$) -- ($(middle
root 3)$);
\draw[/Dynkin diagram/fold style]
($(upper root 3)$) -- ($(middle
root 4)$);
\draw[/Dynkin diagram/fold style]
($(upper root 4)$) -- ($(middle
root 1)$);
\end{scope}
\end{tikzpicture}$ then let $\Delta^1=\big\{\alpha_8,\alpha_6,\alpha_4,\phi \big\}$ and $\Delta^2=\big\{e_1+e_2,e_3+e_4,e_5+e_6,\alpha_1 \big\}$. Observe that $e_8$ and $e_1+e_3+e_5$ are both fixed by $\Pi$ and that $r_s(\mathcal{H}) =2$. So $\big\{e_8,e_1+e_3+e_5\big\}$ is a basis for $E_s$. Hence $\Phi_a = \big\langle e_2-e_4,e_4-e_6,e_6-e_7,e_6+e_7 \big\rangle \cup \big\langle e_1-e_3,e_3-e_5 \big\rangle\cong D_4A_2$.

\vspace{2mm}\noindent For all remaining indices $\mathcal{H}$ listed in the fourth column of Table \ref{E_8}, we observe that $\Psi_0$ is a parabolic subsystem of $\Phi$ and we use Lemma \ref{complemma} to check that $r(\Phi_a)=r(\Psi_0)$. So $\Phi_a=\Psi_0$ since $\Phi_a$ is a parabolic subsystem of $\Phi$ that contains $\Psi_0$. That is, $\mathcal{H}$ is independent in $\Phi$.

\newpage\newgeometry{top=0.7cm,bottom=0.7cm,left=0cm, right=0cm,foot=3mm}
\begin{longtable}{| c | c | c | c | c | c | c | c | c | c | c |}\caption{Isotropic maximal embeddings of abstract indices of type $E_8$}\label{E_8} \\ \hline

\multicolumn{1}{|c|}{} & \multicolumn{1}{c|}{} & \multicolumn{1}{c|}{} & \multicolumn{1}{c|}{} & \multicolumn{1}{c|}{} & \multicolumn{1}{c|}{} & \multicolumn{1}{c|}{} & \multicolumn{3}{c|}{\small{Special fields}} \\

  $\Delta$ & \small{$\hspace{-0.5mm}\Stab_W(\Delta)\hspace{-0.5mm}$} & $\Pi$ & $\textcolor{white}{iiaaa} \mathcal{H} \textcolor{white}{aaaii}$ & $\Psi_0$ & $\Phi_a$ & $\mathcal{G}$ & \small{$\hspace{-0.8mm}\!\cd 1 \!\hspace{-1mm}$} & $\!\R\!$ & $\!\!Q_{\mathfrak{p}}\!\!$ \\ \hline \hline 
\endfirsthead

\multicolumn{11}{|c|}%
{{\bfseries \tablename\ \thetable{} -- continued from previous page}} \\
\hline
\multicolumn{1}{|c|}{} & \multicolumn{1}{c|}{} & \multicolumn{1}{c|}{} & \multicolumn{1}{c|}{} & \multicolumn{1}{c|}{} & \multicolumn{1}{c|}{} & \multicolumn{1}{c|}{} & \multicolumn{3}{c|}{\small{Special fields}} \\

  $\Delta$ & \small{$\hspace{-0.5mm}\Stab_W(\Delta)\hspace{-0.5mm}$} & $\Pi$ & $\mathcal{H}$ & $\Psi_0$ & $\Phi_a$ & $\mathcal{G}$ & \small{$\hspace{-0.8mm}\!\cd 1 \!\hspace{-1mm}$} & $\!\R\!$ & $\!\!Q_{\mathfrak{p}}\!\!$  \\ \hline \hline \endhead
\hline \multicolumn{11}{|c|}{{Continued on next page}} \\ \hline \endfoot 
 \endlastfoot

\multirow{26.5}{*}{$\!\!(A_4)^2\!\!$} & \multirow{26.5}{*}{$\Z_4$} & \multirow{2.5}{*}{$1$} & $(\dynkin{A}{oooo})^2$ & $\varnothing$ & $\varnothing$ & $\begin{tikzpicture}[baseline=0.4ex]\dynkin{E}{VIII}\end{tikzpicture}$ & $\checkmark$ & $\checkmark$ & $\checkmark$  \\ \cdashline{4-11} 
&&& $\!\!\!\dynkin{A}{4} \!\times\! \dynkin{A}{oooo}\!\!\!$ & $A_4$ & $A_4$ &&&&  \topstrut \\ \cline{3-11} 

& & \multirow{17.7}{*}{$\Z_2$} & $\begin{tikzpicture}[baseline=-0.6ex]\dynkin[ply=2,fold radius = 2mm]{A}{oooo}\end{tikzpicture} \times\hspace{-0.5mm} \begin{tikzpicture}[baseline=-0.6ex]\dynkin[ply=2,fold radius = 2mm]{A}{oooo}\end{tikzpicture}$ & $\varnothing$ & $\!\big((A_1)^4\big)''\!$ & $\begin{tikzpicture}[baseline=0.4ex]\dynkin{E}{VIII}\end{tikzpicture}$ & $\checkmark$ & $\checkmark$ & $\checkmark$ \\   \cdashline{4-11} 
&&& \multirow{2.4}{*}{$\begin{tikzpicture}[baseline=-0.6ex]\dynkin[ply=2,fold radius = 2mm]{A}{oooo}\end{tikzpicture} \times\hspace{-0.5mm} \begin{tikzpicture}[baseline=-0.6ex]\dynkin[ply=2,fold radius = 2mm]{A}{o**o}\end{tikzpicture}$} & \multirow{2.4}{*}{$A_2$} & \multirow{2.4}{*}{$D_4A_1$} & $\begin{tikzpicture}[baseline=0.4ex]\dynkin{E}{VIII}\end{tikzpicture}$ & $\xmark$ & $\xmark$ & $\xmark$ \\ 
&&&&&& $\begin{tikzpicture}[baseline=0.4ex]\dynkin{E}{IX}\end{tikzpicture}$ & $\xmark$ & $\checkmark$ & $\xmark$ \\  \cdashline{4-11}
&&& \multirow{3.6}{*}{$\begin{tikzpicture}[baseline=-0.6ex]\dynkin[ply=2,fold radius = 2mm]{A}{oooo}\end{tikzpicture} \times\hspace{-0.5mm} \begin{tikzpicture}[baseline=-0.6ex]\dynkin[ply=2,fold radius = 2mm]{A}{****}\end{tikzpicture}$} & \multirow{3.6}{*}{$A_4$} & \multirow{3.6}{*}{$D_6$} & $\begin{tikzpicture}[baseline=0.4ex]\dynkin{E}{VIII}\end{tikzpicture}$ & $\xmark$ & $\xmark$ & $\xmark$ \\ 
&&&&&& $\begin{tikzpicture}[baseline=0.4ex]\dynkin{E}{IX}\end{tikzpicture}$ & $\xmark$ & $\checkmark$ & $\xmark$ \\
\end{longtable}\restoregeometry

\section*{Acknowledgements}

\noindent The author would like to thank his PhD supervisor Prof. Martin Liebeck for all the help and direction he has given with this project. In addition, the author wishes to thank Prof. Philippe Gille for some help with Galois cohomology, and to Profs. Nikolay
Nikolov, David Evans and Gerhard Röhrle for some useful suggestions.

\end{document}